\documentclass[11pt,reqno]{amsart}
\usepackage[margin=1.1in]{geometry}

\usepackage{amsmath,amsfonts,amsthm,amssymb,enumerate,graphicx}
\usepackage{url,mathtools,tabularx,multirow,esint,hyperref,booktabs}

\numberwithin{equation}{section}

\newcommand{\R}{\mathbb{R}}
\newcommand{\Z}{\mathbb{Z}}
\DeclareMathOperator{\re}{Re}
\DeclareMathOperator{\im}{Im}
\DeclareMathOperator{\dv}{div}
\DeclareMathOperator{\diag}{diag}

\theoremstyle{plain}
\newtheorem{theorem}{Theorem}[section]
\newtheorem{lem}[theorem]{Lemma}
\newtheorem{cor}[theorem]{Corollary}
\newtheorem{proposition}[theorem]{Proposition}

\theoremstyle{definition}

\theoremstyle{remark}
\newtheorem{remark}[theorem]{Remark}

\title{A shrinking torus without rotational symmetry}

\author{Stan Palasek}
\address{Department of Mathematics, Princeton University, Princeton, NJ 08544, USA}
\email{spalasek@princeton.edu}

\begin{document}

\begin{abstract}
We construct a compact embedded self-shrinker for mean curvature flow in $\R^3$ with genus one that is not a surface of revolution. This answers negatively the question whether the Angenent torus is the unique shrinking torus, and thus shows that Brendle's uniqueness theorem in genus zero does not extend to genus one. We show that the new surface has symmetry group $D_3\times\Z_2$ and is locally unique up to rotations. In conformal coordinates, the shrinker equation reduces to a semilinear elliptic system with cubic nonlinearity, which we solve by a Newton--Kantorovich argument. Using rigorous spectral bounds for a Schrödinger operator on the flat torus, we prove that the surface has Morse index 10 and that its only Jacobi fields come from rotations.
\end{abstract}

\maketitle

\section{Introduction}
In this paper, we consider self-similar shrinking solutions of mean curvature flow in $\R^3$. We say a surface $\Sigma\subset\R^3$ is a \emph{self-shrinker} if its mean curvature $H$ satisfies
\begin{align}\label{eq:shrinker}
H=\frac12\langle X,n\rangle,
\end{align}
where $X$ is the position vector and $n$ is a unit normal, with the convention that a sphere of radius $r$ has $H=2/r$ for the outward normal. Equivalently, for $\Sigma$ a self-shrinker, $\sqrt{-t}\,\Sigma$ evolves by mean curvature flow for $t<0$.

By Huisken's monotonicity formula \cite{Hu}, self-shrinkers model the singularities of mean curvature flow. For a flow starting from a closed embedded surface in $\R^3$, Ilmanen \cite{Il} and Bamler and Kleiner \cite{BaK} showed that the tangent flows at the first singular time are given by smooth embedded self-shrinkers of multiplicity one, whose genus is at most that of the initial surface. In particular, the first singularities of an embedded torus are modeled by shrinkers of genus zero or one.

In the genus zero case, Brendle \cite{Br} proved that the only properly embedded self-shrinkers are the plane, the round cylinder, and the round sphere. In genus one, Angenent \cite{An} constructed an embedded shrinking torus of revolution. The uniqueness of the Angenent torus has remained open, both in the class of compact embedded shrinking tori (see, e.g., \cite[Open Problem~1]{DLN}) and in the rotationally symmetric class \cite{KM,DN,Mr,Pe}.

The genus-one question also has a counterpart in minimal surface theory. Indeed, self-shrinkers are minimal surfaces for the conformally flat metric $e^{-|X|^2/4}\delta$ on $\R^3$ \cite{An,CM}. For minimal surfaces in the round $S^3$, Brendle's proof of the Lawson conjecture \cite{Br2} establishes that the Clifford torus is the unique embedded minimal torus up to ambient isometry. Andrews and Li \cite{AL} extended this to constant mean curvature, proving that every embedded torus in $S^3$ with constant mean curvature is a surface of revolution. It is natural to ask whether a similar rigidity holds for compact embedded shrinking tori.

The purpose of this work is to resolve this question in the negative. Specifically, we construct a compact embedded shrinking torus with finite symmetry group, which we prove to be locally unique modulo rotations. We also show that its entropy is strictly larger than that of a particular shrinking torus of revolution, so it does not minimize entropy among compact embedded shrinking tori.

\subsection{The main result}

Self-shrinkers are critical points of the Gaussian area
\begin{align*}
F(\Sigma)=\frac{1}{4\pi}\int_\Sigma e^{-|X|^2/4}\,dA.
\end{align*}
Colding and Minicozzi \cite{CM} introduced the entropy
\begin{align*}
\lambda(\Sigma)=\sup_{x_0\in\R^3,\ t_0>0}\ \frac{1}{4\pi t_0}\int_\Sigma e^{-|X-x_0|^2/(4t_0)}\,dA,
\end{align*}
which is non-increasing along mean curvature flow, and equals $F(\Sigma)$ for closed self-shrinkers.

Let $G\subset O(3)$ be the group of order $12$ generated by
\begin{equation}\label{eq:G}
\rho=\text{rotation by } 2\pi/3 \text{ about the } X_3\text{-axis},\qquad \iota_2=\diag(1,-1,1),\qquad \iota_3=\diag(1,1,-1),
\end{equation}
which is the symmetry group of a right prism over an equilateral triangle and is isomorphic to $D_3\times\Z_2$. Our main result is the following.

\begin{theorem}\label{thm:main}
There exists a real-analytic closed embedded surface $\Sigma\subset\R^3$ of genus one which satisfies \eqref{eq:shrinker} and has the following properties.
\begin{enumerate}
\item The group of isometries of $\R^3$ which map $\Sigma$ onto itself is $G$ as in \eqref{eq:G}. In particular, $\Sigma$ is not a surface of revolution about any line of $\R^3$.
\item The entropy of $\Sigma$ satisfies $|\lambda(\Sigma)-1.8525734|\leq 6.2\cdot10^{-6}$.
\item The surface $\Sigma$ is conformally equivalent to the flat torus $\R^2/(2\pi\Z\times2\pi c\,\Z)$ with $|c-0.93006727552|\leq 8\cdot 10^{-11}$.
\end{enumerate}
\end{theorem}

\begin{cor}\label{cor:angenent}
The surface $\Sigma$ of Theorem~\ref{thm:main} is not congruent to any rotationally symmetric self-shrinker. In particular, the Angenent torus is not the only closed embedded self-shrinker of genus one in $\R^3$.
\end{cor}

\begin{figure}[ht]
\centering
\includegraphics[width=0.95\textwidth]{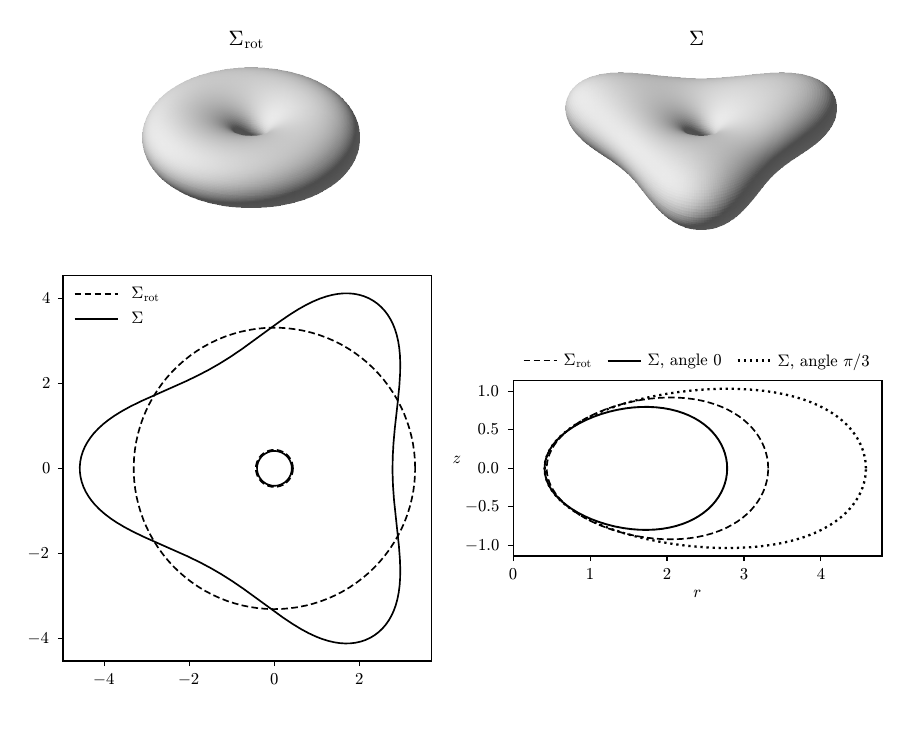}
\caption{The rotationally symmetric torus $\Sigma_{\mathrm{rot}}$ of Theorem~\ref{thm:rot} (top left) and the new surface $\Sigma$ of Theorem~\ref{thm:main} (top right). The bottom left panel shows the intersections with the plane $X_3=0$ of $\Sigma_{\mathrm{rot}}$ (dashed) and $\Sigma$ (solid). The bottom right panel shows meridian sections in the half-plane with coordinates $r=(X_1^2+X_2^2)^{1/2}$ and $z=X_3$. The dashed curve is the profile curve of $\Sigma_{\mathrm{rot}}$, and the solid and dotted curves are the sections of $\Sigma$ sliced in the directions $(1,0,0)$ and $(\cos\frac\pi3,\sin\frac\pi3,0)$, respectively.}
\label{fig:comparison}
\end{figure}

We refer to Figure~\ref{fig:comparison} for a comparison of $\Sigma$ with the torus of revolution constructed in Theorem~\ref{thm:rot}. The existence proof consists of a Newton--Kantorovich argument in the vicinity of an explicit trigonometric polynomial. From there, we obtain quantitative estimates which imply the other claims. See \S\ref{sec:computation} for details on the numerical methods.

\subsection{Index and nullity}\label{sec:intro-index}

For a closed self-shrinker $\Sigma$, we write
\begin{equation}\label{eq:jacobi}
Lf=-\Delta_\Sigma f+\frac12\langle X^{T},\nabla_\Sigma f\rangle-\Big(|A|^2+\frac12\Big)f
\end{equation}
for its Jacobi operator, the linearization of $H-\frac12\langle X,n\rangle$ under perturbations normal to the surface. As an operator on $L^2(\Sigma,e^{-|X|^2/4}dA)$, $L$ is self-adjoint and has discrete spectrum. Then the Morse index of $\Sigma$ is the number of negative eigenvalues of $L$ counted with multiplicity. On every closed self-shrinker, dilation and translation account for four negative directions. Indeed, Colding and Minicozzi \cite{CM} showed that $LH=-H$ and $L\langle v,n\rangle=-\frac12\langle v,n\rangle$ for $v\in\R^3$. Up to constant factors, these eigenfunctions are the normal components of the variations of $\Sigma$ by dilation and by translations. Since rotations about the origin map self-shrinkers to self-shrinkers, the normal components of the rotation fields satisfy $L\langle v\times X,n\rangle=0$.

\begin{theorem}\label{thm:index}
With $\Sigma$ the surface of Theorem~\ref{thm:main} and $L$ its Jacobi operator \eqref{eq:jacobi}, we have the following.
\begin{enumerate}
\item The operator $L$ has exactly $10$ negative eigenvalues, counted with multiplicity.
\item The kernel of $L$ is three-dimensional and is spanned by the functions $\langle e_a\times X,n\rangle$, $a\in\{1,2,3\}$.
\item On the space of $G$-invariant functions, $L$ has exactly three negative eigenvalues and trivial kernel.
\end{enumerate}
\end{theorem}

Of the ten negative directions, four correspond to dilation and translation, and two more are $G$-invariant. The nullity is exactly the dimension of the rotation orbit, so $\Sigma$ is a non-degenerate critical point of the Gaussian area modulo $SO(3)$.

The same computation gives a quantitative spectral gap.

\begin{cor}\label{cor:gap}
Every nonzero eigenvalue of $L$ has modulus at least $1.07\cdot10^{-3}$.
\end{cor}

The description of the kernel also gives local uniqueness, even outside the symmetry class $G$.

\begin{cor}\label{cor:rigidity}
Let $0<\alpha<1$. There exists $\epsilon_0>0$ with the following property: if $u\in C^{2,\alpha}(\Sigma)$ satisfies $\|u\|_{C^{2,\alpha}}<\epsilon_0$ and the surface $\{p+u(p)n(p) : p\in\Sigma\}$ is a self-shrinker, then this surface is the image of $\Sigma$ under a rotation of $\R^3$ about a line through the origin.
\end{cor}

\subsection{Entropy}

Our third result compares the entropy of $\Sigma$ with that of a torus of revolution, which we construct numerically using an analogous existence argument. 

\begin{theorem}\label{thm:rot}
There exists a real-analytic closed embedded surface $\Sigma_{\mathrm{rot}}\subset\R^3$ of genus one which satisfies \eqref{eq:shrinker} and is a surface of revolution about the $X_3$-axis invariant under $\iota_3$. Its entropy satisfies
\begin{align*}
1.85121667167\leq\lambda(\Sigma_{\mathrm{rot}})\leq1.85121667169.
\end{align*}
\end{theorem}

In particular, comparing to the surface $\Sigma$ of Theorem~\ref{thm:main} by means of \eqref{eq:ball-enclosures}, we have
\begin{align*}
\lambda(\Sigma)-\lambda(\Sigma_{\mathrm{rot}})\geq1.35064\cdot10^{-3}.
\end{align*}
This gives the following.

\begin{cor}\label{cor:least}
The surface $\Sigma$ of Theorem~\ref{thm:main} does not globally minimize entropy among the closed embedded self-shrinkers of genus one.
\end{cor}

We remark that numerically, the profile of $\Sigma_{\mathrm{rot}}$ agrees with the Angenent torus, and the enclosure for $\lambda(\Sigma_{\mathrm{rot}})$ agrees with the $1.85122$ obtained by Berchenko-Kogan \cite{BK1}. Whether $\Sigma_{\mathrm{rot}}$ \emph{is} the Angenent torus (or, relatedly, the Drugan--Nguyen torus~\cite{DN}) would be settled by uniqueness of compact embedded rotationally symmetric shrinkers, which is not known.

\subsection{Relation to previous work}

First we point out that the compactness and embeddedness assumptions are essential for the genus-one uniqueness question; see for instance \cite{DK}. In the class of surfaces of revolution in $\R^3$, Kleene and M\o ller \cite{KM} showed that the only complete embedded self-shrinkers are the plane, the round cylinder, the round sphere, or a torus. In 1992, Angenent constructed a self-shrinking embedded torus which is a surface of revolution by a shooting argument \cite{An}. Subsequently, Drugan and Nguyen gave a new construction of such a torus by a variational approach \cite{DN}. It is conjectured that the tori of \cite{An} and \cite{DN} coincide, although this remains open.

To our knowledge, the surface $\Sigma$ constructed here is the first compact embedded self-shrinker of genus one that is not a surface of revolution. Compactness distinguishes it from the genus-one example of Buzano, Nguyen, and Schulz \cite{BNS}, which is noncompact and hence nonrotational by \cite{KM}. Ketover \cite{Ke2} constructed embedded shrinkers of every genus $g\geq1$ with entropy less than $2$ and symmetry group $D_{g+1}\times\Z_2$, for which non-compactness is known for $g$ large and expected based on numerics when $g=1$.

In \cite{CS}, Chu and Sun proved that there is an embedded genus-one self-shrinker with entropy below that of the torus constructed by Drugan and Nguyen, though its symmetry group and compactness are not known. They also proved that a least-entropy complete embedded genus-one shrinker is either noncompact or has Morse index $5$, where the index includes dilation and translations. The existence of such a minimizer was proved by Sun and Wang \cite{SW}. Our compact index-$10$ torus therefore does not minimize entropy among complete embedded genus-one shrinkers. Theorem~\ref{thm:rot} gives the stronger conclusion that it does not minimize entropy even among compact embedded genus-one shrinkers.

At higher genus, M{\o}ller \cite{Mo} constructed compact embedded shrinkers by desingularizing the union of the sphere and the Angenent torus. More recently, Kapouleas and McGrath \cite{KMc} obtained compact examples by a doubling construction. Following numerical work of Chopp \cite{Ch}, Ketover \cite{Ke} and Hoffman--Mart\'in--White \cite{HMW} constructed shrinkers with Platonic and dihedral symmetries, whose compactness is not known in general. In particular, both \cite{Ke2} and \cite{HMW} give a genus-two shrinker with the same symmetries as our surface $\Sigma$. Additional noncompact examples with high genus were constructed by Kapouleas--Kleene--M{\o}ller \cite{KKM} and by Nguyen \cite{Ng}.

Colding and Minicozzi showed that the Morse index of a closed embedded self-shrinker in $\R^3$ must be at least $5$ except in the case of the round sphere \cite{CM}. Impera, Rimoldi, and Savo \cite{IRS} proved the lower bound $\frac23g+3$ for a compact shrinker of genus $g$. For rotational shrinking tori\footnote{Note that Liu and Berchenko-Kogan count the index with the convention to exclude the four directions from dilation and translations, whereas we include them throughout (including when we cite their results).}, Liu \cite{Li} proved the index must be at least $7$; see also the upper and lower bounds proved by Berchenko-Kogan in \cite{BK3}. Berchenko-Kogan's numerical computation for the Angenent torus \cite{BK2} gives full Morse index $9$. For our new surface $\Sigma$, Theorem~\ref{thm:index} rigorously establishes an index of 10.

\subsection{Strategy of the proofs}

We parameterize the surface in conformal coordinates, in which the shrinker equation takes on an explicit cubic nonlinearity. A conformal immersion of a flat torus into $\R^3$ gives a self-shrinker exactly when it is a harmonic map into $(\R^3,e^{-|X|^2/4}\delta)$. For the domain, we take the torus to be $\R^2/(2\pi\Z\times2\pi c\,\Z)$, where $c>0$ is a parameter to be found numerically. Rescaling the second variable, we can regard the immersion as a map $X(x,t)$ on the fixed torus $(\R/2\pi\Z)^2$. With $\sigma=c^{-2}$, the harmonic map equation becomes
\begin{equation}\label{eq:N-intro}
X_{xx}+\sigma X_{tt}-\frac12\big((X\cdot X_x)X_x+\sigma(X\cdot X_t)X_t\big)+\frac14(|X_x|^2+\sigma|X_t|^2)X=0.
\end{equation}
To enforce that $X$ is conformal, we use that the Hopf differential of a harmonic map from a torus is constant. For a $G$-equivariant map this constant is real, and its vanishing is equivalent to the scalar equation
\begin{equation}\label{eq:eta-intro}
\frac{1}{4\pi^2}\int_{(\R/2\pi\Z)^2}(|X_x|^2-\sigma|X_t|^2) dxdt=0.
\end{equation}
Equations \eqref{eq:N-intro}--\eqref{eq:eta-intro} together form a square polynomial system for the Fourier coefficients of $X$ and the parameter $\sigma$. By a continuation argument described below in \S\ref{sec:discovery}, we obtain an approximate solution in the form of an explicit trigonometric polynomial with $45813$ real unknowns. We proceed to construct an exact solution of the full system \eqref{eq:N-intro}--\eqref{eq:eta-intro} using a Newton--Kantorovich argument of radii-polynomial type \cite{Ya,DLM,BL} in a weighted Wiener-type space; see also \cite{NPW,GS}. Because our numerical solution is a polynomial, we can compute exact expressions for the residual and the coefficients of the linearization. The high frequency tail of the linearization can be estimated as a first-order perturbation of the flat Laplacian.

The other geometric conclusions can be inferred from the explicit numerical solution for sufficiently small perturbations. For instance, to show that the surface is an immersion, it suffices to prove a lower bound for $|X_x|$. Similarly, by \cite{CIMW}, embeddedness follows from an upper bound on the Gaussian area. To rule out additional symmetries beyond $G$, we need only show that a third angular moment of the Gaussian area is non-vanishing.

Writing the normal variation as $f=e^{|X|^2/8}\varphi$ and using conformal coordinates, we transform the quadratic form of $L$ into that of the Schr\"odinger operator $-\partial_x^2-\sigma\partial_t^2+V$ on the flat torus. This change preserves the index and nullity, though not the individual eigenvalues. We approximate $V$ by an explicit trigonometric polynomial with a small error in the Wiener norm. The existence bound controls the first derivatives, and bootstrapping with the equation allows us to bound the second derivatives entering $V$. Using a Neumann series, we may also control $|X_x|^{-2}$.

We then count negative eigenvalues separately in each symmetry sector. A Rayleigh--Ritz matrix gives a lower bound, and a Schur complement against the positive high-frequency block gives an upper bound. In both cases we are left with a finite symmetric matrix with rational entries, whose inertia we compute exactly. By comparing with shifted operators, we show that the eigenvalues near zero are exactly the zero eigenvalues coming from rotations.

For Theorem~\ref{thm:rot} in which we compare $\Sigma$ to a rotationally symmetric torus $\Sigma_{\mathrm{rot}}$, we carry out the same Newton--Kantorovich argument, now to maps of the form $X_1+iX_2=e^{ix}a(t)$, $X_3=b(t)$ to construct $\Sigma_{\mathrm{rot}}$.

\subsection{Numerical discovery of the second torus}\label{sec:discovery}

In this section, we describe the method by which we found an explicit approximation of the new torus $\Sigma$. In short, we carry out numerical continuation from a torus of revolution\footnote{The Angenent torus, presumably.}, which we perturb by a threefold rotationally symmetric deformation with amplitude $a\in\mathbb R$. The shrinker equation is balanced by a force which vanishes at $a=0$ by construction, and again at $a\approx 0.3197$, giving the numerical candidate for $\Sigma$. The nearest precedents for this method, to our knowledge, are \cite{MR2943603} and \cite{MR4494205}.

More specifically, we search for $\Sigma$ in the space of $G$-invariant graphs over a round torus, described by functions $q(u,v)$ that are even in $u$ and $v$ and $2\pi/3$-periodic in $u$. Such a graph is parameterized as
\begin{align*}
(u,v)&\mapsto\big((R_0+e^{q}\cos v)\cos u,(R_0+e^{q}\cos v)\sin u,e^{q}\sin v\big),\\
q(u,v)&=\sum_{k\in3\Z_{\geq0}}\sum_{j\geq0}c_{k,j}\cos(ku)\cos(jv),
\end{align*}
where the radius $R_0$ is an additional unknown, and we set $c_{0,1}=0$. We write $S(q)$ for the function $H-\frac12\langle X,n\rangle$ of this surface.

We initiated the search at a rotationally symmetric torus $\Sigma_{\mathrm{rot}}$ which solves $S(q_0)=0$ with $q_0=q_0(v)$. The computed $G$-invariant spectrum suggested that $\Sigma_{\rm rot}$ is isolated among the $G$-invariant self-shrinkers. This rules out a continuation method on a curve of \emph{exact} shrinkers, so we instead follow a curve of forced solutions using a Lyapunov--Schmidt-type method. We perturb in the direction of the non-rotationally symmetric mode (i.e., $k\neq0$) with lowest $G$-invariant eigenvalue $\approx 0.113$ which occurs at angular frequency $k=3$ (see Table~\ref{tab:spectrum}).

With $\phi$ this lowest mode, a fixed function $g$, and a real parameter $a$, we numerically solved
\begin{align*}
S(q)+\Lambda g=0,\qquad \langle q-q_0,\phi\rangle=a.
\end{align*}
for $q$, $R_0$, and $\Lambda$. Here, the second equation prescribes the amplitude of the $\phi$ mode, and $\Lambda g$ is the force needed to hold the surface at this amplitude. The self-shrinkers on the path are the solutions with $\Lambda=0$, so the search reduces to finding the zeros of the scalar function $\Lambda(a)$. We proceeded by a continuation in $a$ and found that $\Lambda$ vanishes at $a=0$ with nonzero slope, as the positive eigenvalue predicts, and vanishes again at $a\approx0.3197$. Two different choices of $g$ gave the same second zero.

Along this path, the Gaussian area increases from $\lambda(\Sigma_{\mathrm{rot}})$ until it attains its maximum at the second zero, then decreases. This agrees with Theorem~\ref{thm:rot}, and suggests that the direction of the path is the additional negative direction of $\Sigma$. Indeed, the computed index is $9$ for $\Sigma_{\mathrm{rot}}$, in agreement with Berchenko-Kogan \cite{BK2}, and $10$ for $\Sigma$, with $2$ and $3$ negative directions among the $G$-invariant functions (Appendix~\ref{sec:numerics}).

At the second zero, we set the forcing term to zero and refined the discretization, which yielded an approximate solution with residual of order $10^{-14}$ in double precision. We recomputed this solution in the conformal parameterization of \S\ref{sec:conformal} with a residual of size $10^{-36}$ in multiple-precision arithmetic. From there, the approximate solution $\bar U$ of \S\ref{sec:ubar} was obtained by Newton's method in double precision at the truncation $\Theta=250$.

In addition to the successful continuation at $k=3$, we carried out the same procedure with $k\in m\Z_{\geq0}$ and the lowest mode of angular frequency $m\in\{2,4,5,6,7,8\}$. In these cases the multiplier $\Lambda$ did not vanish for $a\ne0$ in the range of amplitudes that the computation could resolve. Whether there exist other solutions at larger amplitudes or on other branches remains an interesting open question.

We emphasize that once the approximation of $\Sigma$ is found, the details of this method play no role in the rigorous existence proof. Nonetheless, we believe that these details could have wider applicability to other geometric problems and are therefore of interest.

\subsection{Organization}

\S\ref{sec:conformal} develops the conformal formulation and symmetry reduction. \S\ref{sec:functional} and \S\ref{sec:computation} give the analytic bounds and certified computation, from which \S\ref{sec:geometry} deduces Theorem~\ref{thm:main}. \S\ref{sec:index} proves the index, nullity, and local rigidity results, and \S\ref{sec:rot} proves the entropy comparison. Appendix~\ref{sec:numerics} records some other interesting computational observations, which are not used in the proofs and not certified.

\subsection*{Acknowledgments}
We are grateful to Camillo De Lellis for reading a draft of the manuscript. We also thank the Princeton Department of Mathematics for access to computing resources.

The author made substantial use of Claude and GPT models in searching for the shrinker, developing the computer-assisted proof, and preparing a preliminary draft. The author checked all facets of the proof and substantially rewrote the manuscript.

\section{Conformal formulation}\label{sec:conformal}

\subsection{Conventions}

For an immersed surface in $\R^3$ with local parameterization $X$ and unit normal $n$ we write the metric $\hat g_{\alpha\beta}=X_\alpha\cdot X_\beta$, second fundamental form $h_{\alpha\beta}=X_{\alpha\beta}\cdot n$, and mean curvature $H=-\hat g^{\alpha\beta}h_{\alpha\beta}$. Then $\Delta_\Sigma X=-Hn$. We put
\begin{align*}
S(X)=H-\frac12\langle X,n\rangle
\end{align*}
so the shrinker equation \eqref{eq:shrinker} becomes $S=0$, or equivalently $\Delta_\Sigma X=-\frac12X^\perp$.

The parameter domain is the torus $\mathbb{T}^2=(\R/2\pi\Z)^2$ whose coordinates we denote by $(x,t)$. We write $\langle f\rangle$ for the mean of a function $f$ over $\mathbb{T}^2$. For $\sigma>0$, we put $y=t/\sqrt\sigma$, so that $(x,y)$ are coordinates on $\R^2/(2\pi\Z\times2\pi\sigma^{-1/2}\Z)$, and
\begin{align*}
\Delta_\sigma=\partial_x^2+\sigma\partial_t^2=\partial_x^2+\partial_y^2.
\end{align*}
For a map $X\colon\mathbb{T}^2\to\R^3$ and $\sigma>0$ we define
\begin{equation}
w=e^{-|X|^2/4},\qquad P=|X_x|^2-\sigma|X_t|^2,\qquad Q=X_x\cdot X_t,\qquad \mu=\frac12(|X_x|^2+\sigma|X_t|^2).
\end{equation}
Since $X_y=\sqrt\sigma X_t$, we have $P=|X_x|^2-|X_y|^2$ and $\sqrt\sigma\,Q=X_x\cdot X_y$. The map $X$ is weakly conformal with respect to $dx^2+dy^2$ exactly when $P=Q=0$, and in this case $\mu=|X_x|^2=|X_y|^2$. We define
\begin{equation}\label{eq:N}
\mathcal{N}(X,\sigma)=X_{xx}+\sigma X_{tt}-\frac12\big[(X\cdot X_x)X_x+\sigma(X\cdot X_t)X_t\big]+\frac14(|X_x|^2+\sigma|X_t|^2)X.
\end{equation}
The Christoffel symbols of the metric $e^{-|X|^2/4}\delta$ are $\Gamma^k_{ij}=-\frac14(\delta_{ik}X_j+\delta_{jk}X_i-\delta_{ij}X_k)$, and therefore $\mathcal{N}(X,\sigma)$ is the tension field of the map $X$ from $(\mathbb{T}^2,dx^2+dy^2)$ to $(\R^3,e^{-|X|^2/4}\delta)$.

\subsection{Identities for arbitrary maps}

The harmonic-map equation alone does not impose conformality. We first derive identities that recover the Hopf differential from the residual and relate its normal component to the geometric shrinker equation. They will also be used to identify the kernel of the linearized system.

For an immersion $X$ we write $\hat g=\mu(\delta+\chi)$ in the coordinates $(x,y)$, where $\chi$ is symmetric and trace-free with
\begin{align*}
\chi_{11}=-\chi_{22}=\frac{P}{2\mu},\qquad \chi_{12}=\frac{\sqrt\sigma\,Q}{\mu},\qquad |\chi|^2=\chi_{\alpha\beta}\chi_{\alpha\beta}.
\end{align*}

\begin{proposition}
For every map $X\in C^2(\mathbb{T}^2;\R^3)$ and every $\sigma>0$,
\begin{align}
\partial_x(wP)+2\sigma\partial_t(wQ)&=2w\mathcal{N}(X,\sigma)\cdot X_x, \label{eq:I1}\\
\partial_t(wP)-2\partial_x(wQ)&=-2w\mathcal{N}(X,\sigma)\cdot X_t. \label{eq:I2}
\end{align}
If $X$ is an immersion with unit normal $n$, and $h_{\alpha\beta}$ is computed in the coordinates $(x,y)$, then
\begin{equation}\label{eq:I3}
\mathcal{N}(X,\sigma)\cdot n=-\mu S(X)+\frac12\mu|\chi|^2H+\chi_{11}(h_{11}-h_{22})+2\chi_{12}h_{12}.
\end{equation}
\end{proposition}

\begin{proof}
We have $\partial_xw=-\frac12(X\cdot X_x)w$ and $\partial_tw=-\frac12(X\cdot X_t)w$. Moreover $\partial_xP+2\sigma\partial_tQ=2X_x\cdot\Delta_\sigma X$, since the terms containing $X_t\cdot X_{xt}$ cancel. Hence
\begin{align*}
\partial_x(wP)+2\sigma\partial_t(wQ)=w\big[2X_x\cdot\Delta_\sigma X-\frac12(X\cdot X_x)P-\sigma(X\cdot X_t)Q\big].
\end{align*}
On the other hand \eqref{eq:N} gives
\begin{align*}
2\mathcal{N}\cdot X_x=2X_x\cdot\Delta_\sigma X-(X\cdot X_x)|X_x|^2-\sigma(X\cdot X_t)Q+\frac12(|X_x|^2+\sigma|X_t|^2)(X\cdot X_x),
\end{align*}
which is the same bracket. This proves \eqref{eq:I1}. The identity \eqref{eq:I2} follows in the same way from $\partial_tP-2\partial_xQ=-2X_t\cdot\Delta_\sigma X$.

For \eqref{eq:I3} we note that the terms of $\mathcal{N}$ which contain $X_x$ and $X_t$ are tangential, so that $\mathcal{N}\cdot n=h_{11}+h_{22}+\frac12\mu\langle X,n\rangle$. A symmetric trace-free $2\times2$ matrix satisfies $\chi^2=\frac12|\chi|^2\delta$. Hence $\hat g^{-1}=\mu^{-1}(1-|\chi|^2/2)^{-1}(\delta-\chi)$ and
\begin{align*}
h_{11}+h_{22}=\mu(1-|\chi|^2/2)\hat g^{\alpha\beta}h_{\alpha\beta}+\chi_{\alpha\beta}h_{\alpha\beta}=-\mu(1-|\chi|^2/2)H+\chi_{\alpha\beta}h_{\alpha\beta},
\end{align*}
which gives \eqref{eq:I3}.
\end{proof}

In the complex coordinate $z=x+iy$ the identities \eqref{eq:I1} and \eqref{eq:I2} read $\partial_{\bar z}\big(w(P-2i\sqrt\sigma\,Q)\big)=2w\mathcal{N}\cdot X_z$, where $w(P-2i\sqrt\sigma\,Q)=4wX_z\cdot X_z$ is the coefficient of the Hopf differential of $X$ with respect to the metric $w\delta$.

\begin{cor}\label{cor:equation}
Let $X\colon\mathbb{T}^2\to\R^3$ be a $C^2$ immersion which is conformal with respect to $dx^2+dy^2$. Then
\begin{align*}
\mathcal{N}(X,\sigma)=-\mu S(X)n.
\end{align*}
In particular $X$ satisfies \eqref{eq:shrinker} if and only if $\mathcal{N}(X,\sigma)=0$.
\end{cor}

\begin{proof}
Since $P=Q=0$, the identities \eqref{eq:I1} and \eqref{eq:I2} give $\mathcal{N}\cdot X_x=\mathcal{N}\cdot X_t=0$, and \eqref{eq:I3} with $\chi=0$ gives $\mathcal{N}\cdot n=-\mu S$.
\end{proof}

\begin{cor}\label{cor:hopf}
Let $X\in C^2(\mathbb{T}^2;\R^3)$ and $\sigma>0$ satisfy $\mathcal{N}(X,\sigma)=0$. Then $wP$ and $wQ$ are constant on $\mathbb{T}^2$.
\end{cor}

\begin{proof}
By \eqref{eq:I1} and \eqref{eq:I2} we have $\Delta_\sigma(wP)=0$ in the sense of distributions, so that $wP$ is constant. Then \eqref{eq:I1} and \eqref{eq:I2} give $\partial_t(wQ)=\partial_x(wQ)=0$.
\end{proof}

\subsection{The symmetry class}

The group $G$ of \eqref{eq:G} acts on $\mathbb{T}^2$ through
\begin{align*}
\rho\colon(x,t)\mapsto(x+2\pi/3,t),\qquad \iota_2\colon(x,t)\mapsto(-x,t),\qquad \iota_3\colon(x,t)\mapsto(x,-t).
\end{align*}
We let $C_G$ be the set of $G$-equivariant maps, in the sense of distributions if necessary,
\begin{equation}
C_G=\Big\{X\colon\mathbb{T}^2\to\R^3 : X(x+2\pi/3,t)=\rho X(x,t),\ \ X(-x,t)=\iota_2X(x,t),\ \ X(x,-t)=\iota_3X(x,t)\Big\}.
\end{equation}
A function on $\mathbb{T}^2$ is called $G$-invariant if it is $\frac{2\pi}{3}$-periodic in $x$, even in $x$, and even in $t$.

\begin{lem}\label{lem:fourier}
A map $X$ belongs to $C_G$ if and only if
\begin{equation}\label{eq:ansatz}
\begin{aligned}
X_1+iX_2&=\sum_{m\in\Z}\sum_{j\geq0}z_{m,j}e^{i(3m+1)x}\cos(jt),\\
X_3&=\sum_{m\geq0}\sum_{j\geq1}\zeta_{m,j}\cos(3mx)\sin(jt),
\end{aligned}
\end{equation}
with real coefficients $z_{m,j}$ and $\zeta_{m,j}$.
\end{lem}

\begin{proof}
Let $W=X_1+iX_2$. The equivariance under $\rho$ means that $W(x+\frac{2\pi}{3},t)=e^{2\pi i/3}W(x,t)$ and that $X_3$ is $\frac{2\pi}{3}$-periodic in $x$. Hence the Fourier coefficients $\hat W_{k}$, $k=(k_1,k_2)$, vanish unless $k_1\equiv1$ modulo $3$, and those of $X_3$ vanish unless $k_1\equiv0$ modulo $3$. The equivariance under $\iota_2$ means that $W(-x,t)=\overline{W(x,t)}$ and that $X_3$ is even in $x$, and the equivariance under $\iota_3$ means that $W$ is even and $X_3$ is odd in $t$. For $W$ these two conditions read $\hat W_{k_1,k_2}=\overline{\hat W_{k_1,-k_2}}$ and $\hat W_{k_1,k_2}=\hat W_{k_1,-k_2}$, so that the coefficients are real and even in $k_2$. This is \eqref{eq:ansatz}.
\end{proof}

No component of a map in $C_G$ has a nonzero mean. The operator $\mathcal{N}(\cdot,\sigma)$ commutes with $O(3)$ acting on the target and with the isometries of the flat torus acting on the domain. Hence $\mathcal{N}(X,\sigma)\in C_G$ for $X\in C_G$. We record the parities of the scalar quantities attached to $X\in C_G$. Since
\begin{equation}\label{eq:parity-derivatives}
\begin{aligned}
X_x(-x,t)&=-\iota_2X_x(x,t), & X_t(-x,t)&=\iota_2X_t(x,t),\\
X_x(x,-t)&=\iota_3X_x(x,t), & X_t(x,-t)&=-\iota_3X_t(x,t),
\end{aligned}
\end{equation}
the functions $|X|^2$, $w$, $\mu$, and $P$ are $G$-invariant, and $Q$ is odd in $x$ and odd in $t$. If $X$ is an immersion we always use the unit normal
\begin{equation}\label{eq:normal}
n=\frac{X_x\times X_t}{|X_x\times X_t|}.
\end{equation}
The identity $(RA)\times(RB)=\det(R)R(A\times B)$ for $R\in O(3)$ and \eqref{eq:parity-derivatives} show that $n$ is $G$-equivariant, that is, $n\in C_G$. Consequently the product $fn$ lies in $C_G$ for every $G$-invariant function $f$, and $\xi\cdot n$ is $G$-invariant for every $\xi\in C_G$.

The infinitesimal symmetries of the equation do not belong to $C_G$. Indeed, \eqref{eq:parity-derivatives} shows that $X_x$ and $X_t$ are not in $C_G$ unless they vanish, the field $e_3\times X$ satisfies $(e_3\times X)(-x,t)=-\iota_2(e_3\times X)(x,t)$, and the fields $e_1\times X$ and $e_2\times X$ are not equivariant under $\rho$.

\subsection{The system}

For $\sigma\in\R$ and $X\in C_G$ we define
\begin{equation}\label{eq:F}
\eta(X,\sigma)=\langle |X_x|^2-\sigma|X_t|^2\rangle,\qquad \mathcal{F}(\sigma,X)=\big(\eta(X,\sigma),\ \mathcal{N}(X,\sigma)\big).
\end{equation}
The map $\mathcal{N}$ is cubic in $X$ and affine in $\sigma$, and $\eta$ is quadratic in $X$ and affine in $\sigma$. In terms of the coefficients in \eqref{eq:ansatz}, the unknowns are $\sigma$, $z_{m,j}$, and $\zeta_{m,j}$, and the equations are $\eta=0$ and the vanishing of the coefficients of $\mathcal{N}(X,\sigma)\in C_G$ in the same expansion.

\begin{proposition}\label{prop:zeros}
Let $\sigma>0$ and let $X\in C_G\cap C^1(\mathbb{T}^2;\R^3)$ satisfy $\eta(X,\sigma)=0$ and $\mathcal{N}(X,\sigma)=0$, the second equation in the sense of distributions. Then $X$ is real-analytic and $P=Q=0$ on $\mathbb{T}^2$. If moreover $X_x$ has no zero, then $X$ is an immersion of $\mathbb{T}^2$ which is conformal with respect to $dx^2+\sigma^{-1}dt^2$ and satisfies \eqref{eq:shrinker}.
\end{proposition}

\begin{proof}
The equation reads
\begin{align*}
    \Delta_\sigma X=\mathcal{W}(X,\nabla X),
\end{align*}
where $\mathcal{W}$ is a polynomial and the right side is continuous. Elliptic regularity gives $X\in W^{2,p}$ for every $p<\infty$, hence $X\in C^{1,\alpha}$ for all $\alpha<1$. Thus the right side is H\"older continuous, so Schauder estimates give $X\in C^{2,\alpha}$. Iterating, $X$ is smooth, and Morrey's theorem \cite{Mor} gives analyticity. By Corollary~\ref{cor:hopf}, there are constants $\kappa_1$ and $\kappa_2$ with $wP=\kappa_1$ and $wQ=\kappa_2$. The function $wQ$ is odd in $x$, so that $\kappa_2=0$ and $Q=0$. Next $P=\kappa_1e^{|X|^2/4}$, and therefore
\begin{align*}
0=\eta(X,\sigma)=\langle P\rangle=\kappa_1\langle e^{|X|^2/4}\rangle.
\end{align*}
The mean on the right is positive, so that $\kappa_1=0$ and $P=0$. If $X_x$ has no zero, then $|X_y|=|X_x|>0$ and $X_x\cdot X_y=0$ everywhere, and $X$ is a conformal immersion. Corollary~\ref{cor:equation} then gives $S(X)=0$.
\end{proof}

The unweighted constraint $\eta=0$ is equivalent to vanishing of the remaining Hopf constant at a zero of $\mathcal{N}$, and has the advantage of being polynomial in $X$ and $\sigma$.

\subsection{The linearization and the Jacobi operator}

Let $(\sigma,X)$ be a zero of $\mathcal{F}$ with $\sigma>0$ and $X\in C_G$ a smooth immersion. By Proposition~\ref{prop:zeros} the map $X$ is conformal and $\mu=|X_x|^2=\sigma|X_t|^2$. With the normal \eqref{eq:normal} we put
\begin{align*}
h_{xx}=X_{xx}\cdot n,\qquad h_{xt}=X_{xt}\cdot n,\qquad h_{tt}=X_{tt}\cdot n.
\end{align*}
The operator $L$ of \eqref{eq:jacobi} is the linearization of $S$, in the sense that $\frac{d}{d\varepsilon}S(X+\varepsilon fn)\big|_{\varepsilon=0}=Lf$. Let $s\in\R$ and let $\xi\colon\mathbb{T}^2\to\R^3$ be of class $C^2$. We decompose
\begin{equation}\label{eq:decomposition}
\xi=fn+\phi X_x+\psi X_t,
\end{equation}
so that $f=\xi\cdot n$, $\phi=\xi\cdot X_x/\mu$, and $\psi=\sigma\xi\cdot X_t/\mu$, and we denote by $P'$, $Q'$, and $\mathcal{N}'$ the derivatives of $P$, $Q$, and $\mathcal{N}$ at $(X,\sigma)$ in the direction $(\xi,s)$.

\begin{lem}\label{lem:linearization}
The derivatives $P'$, $Q'$, and $\mathcal{N}'$ are given by
\begin{align}
P'&=2\mu(\phi_x-\psi_t)-\frac{\mu}{\sigma}s-2f\,(h_{xx}-\sigma h_{tt}), \label{eq:L1}\\
Q'&=\mu\Big(\phi_t+\frac{\psi_x}{\sigma}\Big)-2fh_{xt}, \label{eq:L2}\\
\mathcal{N}'\cdot n&=-\mu Lf+\frac{(h_{xx}-\sigma h_{tt})P'+4\sigma h_{xt}Q'}{2\mu}, \label{eq:L3}\\
\mathcal{N}'\cdot X_x&=\frac{1}{2w}\big[\partial_x(wP')+2\sigma\partial_t(wQ')\big], \label{eq:L4}\\
\mathcal{N}'\cdot X_t&=-\frac{1}{2w}\big[\partial_t(wP')-2\partial_x(wQ')\big]. \label{eq:L5}
\end{align}
Moreover the derivative of $\eta$ in the direction $(\xi,s)$ equals $\langle P'\rangle$.
\end{lem}

\begin{proof}
We have $P'=2X_x\cdot\xi_x-2\sigma X_t\cdot\xi_t-s|X_t|^2$ and $Q'=\xi_x\cdot X_t+X_x\cdot\xi_t$. Inserting \eqref{eq:decomposition} and using $X_x\cdot n_x=-h_{xx}$, $X_t\cdot n_t=-h_{tt}$, $X_x\cdot n_t=X_t\cdot n_x=-h_{xt}$, $X_x\cdot X_t=0$, and $|X_x|^2=\sigma|X_t|^2=\mu$, we obtain
\begin{align*}
P'=2\mu(\phi_x-\psi_t)-\frac{\mu}{\sigma}s-2f(h_{xx}-\sigma h_{tt})+\phi\partial_xP+\psi\partial_tP
\end{align*}
and the analogous formula for $Q'$ with the additional terms $\phi\partial_xQ+\psi\partial_tQ$. Since $P=Q=0$ on $\mathbb{T}^2$, this gives \eqref{eq:L1} and \eqref{eq:L2}. The identities \eqref{eq:L4} and \eqref{eq:L5} are the derivatives of \eqref{eq:I1} and \eqref{eq:I2} at a point where $\mathcal{N}=0$ and $P=Q=0$. For \eqref{eq:L3} we differentiate \eqref{eq:I3}, which is valid along $X+\varepsilon\xi$ for small $\varepsilon$. Since $\mathcal{N}=0$ at the base point, the derivative of the left side is $\mathcal{N}'\cdot n$. The term $\frac12\mu|\chi|^2H$ is quadratic in $\chi$ and does not contribute. The derivative of $S$ is $Lf+\phi\partial_xS+\psi\partial_tS=Lf$, because $S$ is invariant under reparameterizations and vanishes on the base surface. Finally $\chi_{11}'=P'/(2\mu)$, $\chi_{12}'=\sqrt\sigma\,Q'/\mu$, $h_{11}-h_{22}=h_{xx}-\sigma h_{tt}$, and $h_{12}=\sqrt\sigma\,h_{xt}$. The last assertion is immediate from \eqref{eq:F}.
\end{proof}

\begin{proposition}\label{prop:jacobi}
Let $(\sigma,X)$ be a zero of $\mathcal{F}$ with $\sigma>0$ and $X\in C_G$ a smooth immersion.
\begin{enumerate}
\item Let $s\in\R$ and $\xi\in C_G\cap C^2$ satisfy $D\mathcal{F}(\sigma,X)(s,\xi)=0$. Then $f=\xi\cdot n$ is $G$-invariant and $Lf=0$. If $f=0$, then $s=0$ and $\xi=0$.
\item Let $f$ be a smooth $G$-invariant function on $\mathbb{T}^2$ with $Lf=0$. Then there exist $s\in\R$ and a smooth $\xi\in C_G$ with $\xi\cdot n=f$ and $D\mathcal{F}(\sigma,X)(s,\xi)=0$.
\end{enumerate}
\end{proposition}

\begin{proof}
Let $(s,\xi)$ be as in the first assertion. The function $f$ is $G$-invariant because $n$ and $\xi$ are $G$-equivariant. By \eqref{eq:L4} and \eqref{eq:L5} we have $\Delta_\sigma(wP')=0$, so that $wP'$ is a constant $\kappa'$, and then $wQ'$ is constant. Since $\xi\in C_G$, the function $Q'$ is odd in $x$, and hence $Q'=0$. Since the derivative of $\eta$ vanishes, $0=\langle P'\rangle=\kappa'\langle 1/w\rangle$ and $P'=0$. Now \eqref{eq:L3} gives $Lf=0$. Assume that $f=0$. Then \eqref{eq:L1} and \eqref{eq:L2} read
\begin{align*}
\phi_x-\psi_t=\frac{s}{2\sigma},\qquad \phi_t+\frac{\psi_x}{\sigma}=0.
\end{align*}
The mean of the first equation gives $s=0$. Then $\Delta_\sigma\phi=\Delta_\sigma\psi=0$, so that $\phi$ and $\psi$ are constant. By \eqref{eq:parity-derivatives} the function $\phi$ is odd in $x$ and $\psi$ is odd in $t$. Hence $\phi=\psi=0$ and $\xi=0$.

Conversely let $f$ be smooth and $G$-invariant with $Lf=0$. We put
\begin{align*}
R_1=\frac{f\,(h_{xx}-\sigma h_{tt})}{\mu},\qquad R_2=\frac{2fh_{xt}}{\mu},\qquad s=-2\sigma\langle R_1\rangle.
\end{align*}
By \eqref{eq:parity-derivatives} and equivariance of $n$, we see that $R_1$ is $G$-invariant, while $R_2$ is odd in $x$ and odd in $t$. Let $\phi$ and $\psi$ be the zero-average solutions of
\begin{align*}
\Delta_\sigma\phi=\partial_xR_1+\sigma\partial_tR_2,\qquad \Delta_\sigma\psi=\sigma(\partial_xR_2-\partial_tR_1).
\end{align*}
Then $\phi$ is odd in $x$ and even in $t$, $\psi$ is even in $x$ and odd in $t$, and both are $\frac{2\pi}{3}$-periodic in $x$. Hence $\xi=fn+\phi X_x+\psi X_t$ is a smooth element of $C_G$. Defining the functions $E_1=\phi_x-\psi_t-R_1+\langle R_1\rangle$ and $E_2=\phi_t+\psi_x/\sigma-R_2$, we see that they satisfy
\begin{align*}
\partial_xE_1+\sigma\partial_tE_2=0,\qquad \partial_xE_2-\partial_tE_1=0.
\end{align*}
and are therefore harmonic, hence constant. Since they have zero mean, we see that $E_1=E_2=0$. Thus \eqref{eq:L1} and \eqref{eq:L2} give $P'=Q'=0$. Combining \eqref{eq:L3}, \eqref{eq:L4}, \eqref{eq:L5}, and $Lf=0$ gives $\mathcal{N}'=0$. Finally, we see that the derivative of $\eta$ is $\langle P'\rangle=0$.
\end{proof}

\section{Functional setting}\label{sec:functional}

\subsection{Norms and spaces}

We fix a number $\bar\sigma\geq1$, which will be the parameter $\sigma$ for the approximate solution, and define
\begin{align*}
\theta_k=(k_1^2+\bar\sigma k_2^2)^{1/2}\quad\text{for } k=(k_1,k_2)\in\Z^2,\qquad \tau=\bar\sigma^{-1/2},\qquad \bar\Delta=\partial_x^2+\bar\sigma\partial_t^2.
\end{align*}
For a scalar function
\begin{align*}
    V=\sum_{k\in\mathbb Z^2}\hat V_ke^{i(k_1x+k_2t)}
\end{align*}
on $\mathbb{T}^2$, we define the Wiener-type norms
\begin{align*}
|V|_{1}=\sum_{k}|\hat V_k|,\qquad \|V\|_{\theta}=\sum_k\theta_k|\hat V_k|.
\end{align*}
For vector-valued functions, the norm is taken as the sum of the norms of the cartesian components. For $3\times3$ matrix-valued functions, we define
\begin{align*}
    \|M\|=\max_b\sum_a|M^{ab}|_{1}.
\end{align*}

We use the following basis of $C_G$, which corresponds to the expansion \eqref{eq:ansatz}. Let $\mathcal{I}$ be the set of indices
\begin{align*}
    \mathcal I=\big\{i=(W;m,j):m\in\Z,\,j\geq0\big\}\cup\big\{i=(3;m,j):m\geq0,\,j\geq1\big\}.
\end{align*}
Then for $i\in\mathcal I$, we define the basis
\begin{align*}
e_i=\begin{cases}
    (\cos((3m+1)x)\cos(jt),\,\sin((3m+1)x)\cos(jt),\,0),&i=(W;m,j)\\
    (0,0,\cos(3mx)\sin(jt)),&i=(3;m,j).
\end{cases}
\end{align*}
In these two cases we define the frequency as $k(i)=(3m+1,j)$ and $k(i)=(3m,j)$, respectively, and set $\theta_i=\theta_{k(i)}$. Each component of $e_i$ is supported on frequencies $k$ with $\theta_k=\theta_i$, and for distinct indices the supports are disjoint in each component. A direct count of the exponentials gives
\begin{align*}
|e_i|_{1}=2\text{ for } i=(W;m,j),\qquad |e_i|_{1}=1\text{ for } i=(3;m,j),\qquad \|e_i\|_{\theta}=\theta_i|e_i|_{1}.
\end{align*}
Consequently, for $X=\sum_{i\in\mathcal{I}}x_ie_i$ in $C_G$,
\begin{equation}\label{eq:norms-basis}
|X|_{1}=\sum_{i\in\mathcal{I}}|x_i||e_i|_{1},\qquad \|X\|_{\theta}=\sum_{i\in\mathcal{I}}\omega_i|x_i|
\end{equation}
where we have set $\omega_i=\theta_i|e_i|_1$. We fix a weight $\vartheta>0$ and define the Banach spaces
\begin{align*}
\mathcal{X}&=\big\{U=(\sigma,X) : \sigma\in\R,\ X\in C_G,\ \|X\|_{\theta}<\infty\big\}, &\|U\|_{\mathcal{X}}&=\vartheta|\sigma|+\|X\|_{\theta},\\
\mathcal{Y}&=\Big\{(y,Y) : y\in\R,\ Y\in C_G,\ \sum_k\theta_k^{-1}|\hat Y_k|<\infty\Big\}, &\|(y,Y)\|_{\mathcal{Y}}&=|y|+\sum_k\theta_k^{-1}|\hat Y_k|.
\end{align*}
By \eqref{eq:norms-basis}, the space $\mathcal{X}$ is a weighted $\ell^1$ space over the index set $\{\sigma\}\cup\mathcal{I}$ with the weights $\vartheta$ and $\omega_i$. Hence the norm of a bounded operator $\mathcal{L}$ on $\mathcal{X}$ with matrix $(\mathcal{L}_{i'i})$ is
\begin{equation}\label{eq:opnorm}
\|\mathcal{L}\|=\sup_{i}\ \frac{1}{\omega_i}\sum_{i'}\omega_{i'}|\mathcal{L}_{i'i}|,
\end{equation}
where $i$ and $i'$ run over $\{\sigma\}\cup\mathcal{I}$ and $\omega_\sigma=\vartheta$.

\begin{lem}\label{lem:algebra}
Let $V$, $V'$ be scalar functions, $Y$, $Y'$ vector fields, and $M$ a matrix of functions on $\mathbb{T}^2$ with finite norms $|\cdot|_{1}$ and $\|M\|$. Then
\begin{equation}\label{eq:N1}
|VV'|_{1}\leq|V|_{1}|V'|_{1},\qquad |Y\cdot Y'|_{1}\leq|Y|_{1}|Y'|_{1},\qquad |MY|_{1}\leq\|M\||Y|_{1},\qquad |\langle V\rangle|\leq|V|_{1}.
\end{equation}
For $X\in C_G$ with $\|X\|_{\theta}<\infty$,
\begin{equation}\label{eq:N2}
|X|_{1}\leq\|X\|_{\theta},\qquad |X_x|_{1}\leq\|X\|_{\theta},\qquad |X_t|_{1}\leq\tau\|X\|_{\theta},
\end{equation}
and the Euclidean norms of $X$, $X_x$, and $X_t$ are bounded pointwise by $|X|_{1}$, $|X_x|_{1}$, and $|X_t|_{1}$. In particular $X$ is of class $C^1$. Finally $\|\bar\Delta^{-1}Y\|_{\theta}=\sum_k\theta_k^{-1}|\hat Y_k|$ for $Y\in C_G$.
\end{lem}

\begin{proof}
The inequalities \eqref{eq:N1} express that the Wiener algebra is a Banach algebra. For \eqref{eq:N2} we use that $|k_1|\leq\theta_k$ and $|k_2|\leq\tau\theta_k$, and that $\theta_k\geq1$ on the frequencies of $C_G$. Indeed, the first two components have $|k_1|\geq1$, and the third component has $|k_2|\geq1$ and $\bar\sigma\geq1$. The last assertion holds because $\bar\Delta$ has the symbol $-\theta_k^2$.
\end{proof}

We write
\begin{equation}\label{eq:N12}
\mathcal{N}(X,\sigma)=X_{xx}+\sigma X_{tt}-\mathcal{N}_1(X)-\sigma\mathcal{N}_2(X),\qquad
\begin{aligned}
\mathcal{N}_1(X)&=\frac12(X\cdot X_x)X_x-\frac14|X_x|^2X,\\
\mathcal{N}_2(X)&=\frac12(X\cdot X_t)X_t-\frac14|X_t|^2X.
\end{aligned}
\end{equation}
With the trilinear form
\begin{equation}\label{eq:cubic}
\mathcal{C}(a,b,c)=\frac12(a\cdot b)c-\frac14(b\cdot c)a,\qquad |\mathcal{C}(a,b,c)|_{1}\leq\frac34|a|_{1}|b|_{1}|c|_{1},
\end{equation}
we have $\mathcal{N}_1(X)=\mathcal{C}(X,X_x,X_x)$ and $\mathcal{N}_2(X)=\mathcal{C}(X,X_t,X_t)$.

\begin{lem}\label{lem:C1}
The map $\mathcal{F}$ of \eqref{eq:F} is cubic in $X$ and affine in $\sigma$, hence a polynomial of degree four in $(\sigma,X)$. Each term is a bounded multilinear map from $\mathcal{X}$ to $\mathcal{Y}$, and therefore $\mathcal{F}\in C^1(\mathcal{X},\mathcal{Y})$.
\end{lem}

\begin{proof}
For $(\sigma,X)\in\mathcal{X}$, the distribution $X_{xx}+\sigma X_{tt}$ has the coefficients $-(k_1^2+\sigma k_2^2)\hat X_k$ which are bounded by $k_1^2+|\sigma|k_2^2\leq\max(1,|\sigma|/\bar\sigma)\theta_k^2$. By \eqref{eq:N1}, \eqref{eq:N2}, and \eqref{eq:cubic}, the cubic terms $\mathcal{N}_1(X)$ and $\mathcal{N}_2(X)$ have finite $|\cdot|_1$ norm. Since $\theta_k\geq1$, this dominates the norm of $\mathcal{Y}$. The $\eta$ term is bounded by \eqref{eq:N1} and \eqref{eq:N2}.
\end{proof}

\subsection{The Newton--Kantorovich theorem}

We use the following contraction criterion, a standard form of the Newton--Kantorovich argument used in validated numerics \cite{Ya,DLM,BL}. We include the proof to specify the bounds needed below.

\begin{theorem}\label{thm:radii}
Let $\mathcal{X}$ and $\mathcal{Y}$ be Banach spaces, $\mathcal{F}\in C^1(\mathcal{X},\mathcal{Y})$ and $\bar U\in\mathcal{X}$. Let $\mathcal{A}\colon\mathcal{Y}\to\mathcal{X}$ and $\mathcal{A}^\dagger\colon\mathcal{X}\to\mathcal{Y}$ be bounded linear maps, and assume that $\mathcal{A}$ is injective. Let $Y_0,Z_0,Z_1\geq0$ and let $Z_2\colon[0,\infty)\to[0,\infty)$ be non-decreasing, with
\begin{align}
\|\mathcal{A}\mathcal{F}(\bar U)\|_{\mathcal{X}}&\leq Y_0, \label{eq:Y0}\\
\|I-\mathcal{A}\mathcal{A}^\dagger\|&\leq Z_0, \label{eq:Z0}\\
\|\mathcal{A}(D\mathcal{F}(\bar U)-\mathcal{A}^\dagger)\|&\leq Z_1, \label{eq:Z1}\\
\|\mathcal{A}(D\mathcal{F}(\bar U+h)-D\mathcal{F}(\bar U))\|&\leq Z_2(r)r\qquad\text{for } \|h\|_{\mathcal{X}}\leq r, \label{eq:Z2}
\end{align}
where the norms in \eqref{eq:Z0}, \eqref{eq:Z1}, and \eqref{eq:Z2} are operator norms on $\mathcal{X}$. Assume that $r_0>0$ satisfies
\begin{equation}\label{eq:radii}
p(r_0):=Z_2(r_0)r_0^2-(1-Z_0-Z_1)r_0+Y_0<0.
\end{equation}
Then $\mathcal{F}$ has a unique zero $U_*$ in the closed ball of radius $r_0$ about $\bar U$ in $\mathcal{X}$. Moreover $\mathcal{A}D\mathcal{F}(U_*)$ is invertible on $\mathcal{X}$, and in particular $D\mathcal{F}(U_*)$ is injective.
\end{theorem}

\begin{proof}
Let $\mathcal{T}(U)=U-\mathcal{A}\mathcal{F}(U)$. Then $\mathcal{T}\in C^1(\mathcal{X},\mathcal{X})$, and the fixed points of $\mathcal{T}$ are the zeros of $\mathcal{F}$ because $\mathcal{A}$ is injective. For $U$ in the closed ball $B$ of radius $r_0$ about $\bar U$ we write
\begin{align*}
D\mathcal{T}(U)=(I-\mathcal{A}\mathcal{A}^\dagger)-\mathcal{A}(D\mathcal{F}(\bar U)-\mathcal{A}^\dagger)-\mathcal{A}(D\mathcal{F}(U)-D\mathcal{F}(\bar U)),
\end{align*}
and \eqref{eq:Z0}, \eqref{eq:Z1}, \eqref{eq:Z2}, and \eqref{eq:radii} give
\begin{align*}
\|D\mathcal{T}(U)\|\leq Z_0+Z_1+Z_2(r_0)r_0=:q<1-\frac{Y_0}{r_0}\leq1.
\end{align*}
By the mean value inequality $\mathcal{T}$ is a contraction on $B$ with constant $q$, and for $U\in B$
\begin{align*}
\|\mathcal{T}(U)-\bar U\|_{\mathcal{X}}\leq\|\mathcal{T}(U)-\mathcal{T}(\bar U)\|_{\mathcal{X}}+\|\mathcal{A}\mathcal{F}(\bar U)\|_{\mathcal{X}}\leq qr_0+Y_0<r_0.
\end{align*}
Hence $\mathcal{T}$ maps $B$ into itself and has a unique fixed point $U_*$ in $B$. Since $\|I-\mathcal{A}D\mathcal{F}(U_*)\|\leq q<1$, the operator $\mathcal{A}D\mathcal{F}(U_*)$ is invertible on $\mathcal{X}$ by a Neumann series.
\end{proof}

\subsection{The approximate inverse}

Let $\bar U=(\bar\sigma,\bar X)$, where $\bar X\in C_G$ is a trigonometric polynomial and $\bar\sigma\geq1$. We fix $\Theta>0$ such that $\bar X$ is a combination of the fields $e_i$ with $\theta_i\leq\Theta$, and we put
\begin{align*}
K=\{i\in\mathcal{I} : \theta_i\leq\Theta\},\qquad \Theta_+=\min\{\theta_i : i\in\mathcal{I}\setminus K\}.
\end{align*}
Let $\Pi$ be the projection of $\mathcal{X}$ onto the span of the coordinate $\sigma$ and of the fields $e_i$, $i\in K$, and let $\Pi_\infty=I-\Pi$. We use the same notation in $\mathcal{Y}$, where $\Pi$ projects onto the first component and the fields $e_i$, $i\in K$. By \eqref{eq:N} and \eqref{eq:N12} the derivative of $\mathcal{F}$ at $\bar U$ is
\begin{equation}\label{eq:DF}
D\mathcal{F}(\bar U)(s,\xi)=\Big(2\langle \bar X_x\cdot\xi_x\rangle-2\bar\sigma\langle \bar X_t\cdot\xi_t\rangle-s\langle |\bar X_t|^2\rangle\ ,\ \ \bar\Delta\xi-\mathcal{D}\xi+s\big(\bar X_{tt}-\mathcal{N}_2(\bar X)\big)\Big),
\end{equation}
where $\mathcal{D}\xi=M_0\xi+M_1\xi_x+M_2\xi_t$ with the matrices
\begin{equation}\label{eq:M}
\begin{aligned}
M_0^{ab}&=\frac12(\bar X_{a,x}\bar X_{b,x}+\bar\sigma\bar X_{a,t}\bar X_{b,t})-\frac14(|\bar X_x|^2+\bar\sigma|\bar X_t|^2)\delta_{ab},\\
M_1^{ab}&=\frac12(\bar X_{a,x}\bar X_b-\bar X_a\bar X_{b,x})+\frac12(\bar X\cdot\bar X_x)\delta_{ab},\\
M_2^{ab}&=\bar\sigma\Big[\frac12(\bar X_{a,t}\bar X_b-\bar X_a\bar X_{b,t})+\frac12(\bar X\cdot\bar X_t)\delta_{ab}\Big].
\end{aligned}
\end{equation}
The entries of $M_0$, $M_1$, and $M_2$ are trigonometric polynomials supported on frequencies $l$ with $\theta_l\leq2\Theta$. We denote by $\mathcal{D}_{i'i}$ the coefficients of $\mathcal{D} e_i=\sum_{i'\in\mathcal{I}}\mathcal{D}_{i'i}e_{i'}$.

Let $D\mathcal{F}_K=\Pi D\mathcal{F}(\bar U)\Pi$, regarded as a $(|K|+1)\times(|K|+1)$ square matrix, and let $\mathcal{A}_K$ be a real matrix of the same order, which will be a numerical approximation of the inverse of $D\mathcal{F}_K$. We define
\begin{equation}\label{eq:A}
\mathcal{A}=\mathcal{A}_K\Pi+\bar\Delta^{-1}\Pi_\infty,\qquad \mathcal{A}^\dagger=D\mathcal{F}_K\Pi+\bar\Delta\Pi_\infty.
\end{equation}
Thus $\mathcal{A}$ maps the first component and the coefficients with index in $K$ of an element of $\mathcal{Y}$ to the corresponding coordinates of $\mathcal{X}$ through $\mathcal{A}_K$, and it acts diagonally on the remaining coefficients. We introduce the constants
\begin{equation}\label{eq:constants}
\begin{aligned}
\alpha_\eta&=\|\mathcal{A}_K(1,0)\|_{\mathcal{X}}, &\alpha_i&=\frac{\|\mathcal{A}_K(0,e_i)\|_{\mathcal{X}}}{|e_i|_{1}}\quad(i\in K),\\
\alpha_{\max}&=\max\Big(\max_{i\in K}\alpha_i\,,\ \frac{1}{\Theta_+}\Big), &\alpha_{tt}&=\max\Big(\max_{i\in K}\frac{k_2(i)^2\alpha_i}{\theta_i}\,,\ \frac1{\bar\sigma}\Big),
\end{aligned}
\end{equation}
and $\beta_i=\alpha_i|e_i|_{1}$ for $i\in K$, $\beta_i=|e_i|_{1}/\theta_i$ for $i\in\mathcal{I}\setminus K$.

\begin{lem}\label{lem:A}
The operator $\mathcal{A}\colon\mathcal{Y}\to\mathcal{X}$ is bounded. For $(y,Y)\in\mathcal{Y}$ with $Y=\sum_{i}y_ie_i$,
\begin{equation}\label{eq:A-bound}
\|\mathcal{A}(y,Y)\|_{\mathcal{X}}\leq\alpha_\eta|y|+\sum_{i\in\mathcal{I}}\beta_i|y_i|.
\end{equation}
If $|Y|_{1}<\infty$, then $\|\mathcal{A}(0,Y)\|_{\mathcal{X}}\leq\alpha_{\max}|Y|_{1}$. If $\xi\in C_G$ and $\|\xi\|_{\theta}<\infty$, then $\|\mathcal{A}(0,\xi_{tt})\|_{\mathcal{X}}\leq\alpha_{tt}\|\xi\|_{\theta}$. If $\|I-\mathcal{A}_KD\mathcal{F}_K\|<1$, then $\mathcal{A}$ is injective.
\end{lem}

\begin{proof}
For $i\notin K$ we have $\|\bar\Delta^{-1}e_i\|_{\theta}=\omega_i/\theta_i^2=\beta_i$, and for $i\in K$ we have $\|\mathcal{A}_K(0,e_i)\|_{\mathcal{X}}=\beta_i$. This gives \eqref{eq:A-bound} by the triangle inequality. The second assertion follows from $\beta_i\leq\alpha_{\max}|e_i|_{1}$ and \eqref{eq:norms-basis}. For the third assertion let $\xi=\sum_ix_ie_i$, so that $\xi_{tt}=-\sum_ik_2(i)^2x_ie_i$ and
\begin{align*}
\|\mathcal{A}(0,\xi_{tt})\|_{\mathcal{X}}\leq\sum_{i\in K}\frac{k_2(i)^2\alpha_i}{\theta_i}\omega_i|x_i|+\sum_{i\notin K}\frac{k_2(i)^2}{\theta_i^2}\omega_i|x_i|\leq\alpha_{tt}\|\xi\|_{\theta},
\end{align*}
because $k_2(i)^2\leq\theta_i^2/\bar\sigma$. If $\|I-\mathcal{A}_KD\mathcal{F}_K\|<1$, then $\mathcal{A}_KD\mathcal{F}_K$ is invertible, hence the square matrix $\mathcal{A}_K$ is invertible, and $\mathcal{A}$ is injective.
\end{proof}

\subsection{The bounds}

We now define the constants $Y_0$, $Z_0$, $Z_1$ and the function $Z_2$ and verify \eqref{eq:Y0}--\eqref{eq:Z2} for the map $\mathcal{F}$ of \eqref{eq:F} and the operators \eqref{eq:A}. We write $\mathcal{N}(\bar U)=\sum_i\mathcal{N}(\bar U)_ie_i$ and use the same notation for the coefficients of other elements of $C_G$.

The field $\mathcal{N}(\bar U)$ is a trigonometric polynomial supported in $\theta_k\leq3\Theta$, and $\mathcal{N}(\bar U)_i=-(\mathcal{N}_1(\bar X)+\bar\sigma\mathcal{N}_2(\bar X))_i$ for $i\notin K$. We put
\begin{equation}\label{eq:Y0-def}
Y_0=\alpha_\eta|\eta(\bar U)|+\sum_{i\in\mathcal{I}}\beta_i|\mathcal{N}(\bar U)_i|,\qquad Z_0=\|I-\mathcal{A}_KD\mathcal{F}_K\|,
\end{equation}
where the last norm is the operator norm \eqref{eq:opnorm} on $\Pi\mathcal{X}$. The bound \eqref{eq:Y0} holds by \eqref{eq:A-bound}. Since $\bar\Delta$ is diagonal, $\mathcal{A}\mathcal{A}^\dagger=\mathcal{A}_KD\mathcal{F}_K\Pi+\Pi_\infty$, and \eqref{eq:Z0} holds.

For the bound $Z_1$ we fix $\Theta_0>0$, put $\Theta'=\Theta+\Theta_0$, and denote by $m_p$, $p\in\{0,1,2\}$, the norm $\|\cdot\|$ of the matrix obtained from $M_p$ by removing the Fourier modes $l$ with $\theta_l\leq\Theta_0$. We define the column sums
\begin{equation}\label{eq:columns}
c_\sigma=\frac1\vartheta\sum_{i'\notin K}\beta_{i'}|\mathcal{N}_2(\bar X)_{i'}|,\qquad
c_i=\frac{1}{\omega_i}\sum_{i'\notin K}\beta_{i'}|\mathcal{D}_{i'i}|\ \ (i\in K),\qquad
c_i=\frac1{\omega_i}\sum_{i'\in\mathcal{I}}\beta_{i'}|\mathcal{D}_{i'i}|\ \ (i\notin K),
\end{equation}
and
\begin{equation}\label{eq:Z1-def}
\begin{aligned}
Z_1^{\mathrm{far}}&=\frac{1}{\Theta}\Big[\|M_1\|+\tau\|M_2\|+\frac{\|M_0\|}{\Theta'}\Big]+\alpha_{\max}\Big[m_1+\tau m_2+\frac{m_0}{\Theta'}\Big],\\
Z_1&=\max\Big(c_\sigma\,,\ \max_{\theta_i\leq\Theta'}c_i\,,\ Z_1^{\mathrm{far}}\Big).
\end{aligned}
\end{equation}

\begin{lem}\label{lem:Z1}
With these definitions $\|\mathcal{A}(D\mathcal{F}(\bar U)-\mathcal{A}^\dagger)\|\leq Z_1$.
\end{lem}

\begin{proof}
Let $h=(s,\xi)\in\mathcal{X}$. Since $\bar X$ is supported in $K$, the first component of $D\mathcal{F}(\bar U)(0,\Pi_\infty\xi)$ in \eqref{eq:DF} vanishes by the orthogonality of the exponentials, and $\Pi_\infty\bar X_{tt}=0$. Therefore
\begin{equation}\label{eq:Z1-identity}
\big(D\mathcal{F}(\bar U)-\mathcal{A}^\dagger\big)h=\Big(0\ ,\ -\mathcal{D}\Pi_\infty\xi-\Pi_\infty\mathcal{D}\Pi\xi-s\Pi_\infty\mathcal{N}_2(\bar X)\Big).
\end{equation}
By \eqref{eq:opnorm} it suffices to bound the columns of $\mathcal{A}(D\mathcal{F}(\bar U)-\mathcal{A}^\dagger)$. By \eqref{eq:Z1-identity} and \eqref{eq:A-bound} the norm of the column of $\sigma$, divided by $\vartheta$, is at most $c_\sigma$, and the norm of the column of $e_i$, divided by $\omega_i$, is at most $c_i$. It remains to show that $c_i\leq Z_1^{\mathrm{far}}$ when $\theta_i>\Theta'$.

We split $M_p=M_p^{\leq}+M_p^{>}$, where $M_p^{\leq}$ contains the Fourier modes $l$ with $\theta_l\leq\Theta_0$, and correspondingly $\mathcal{D}=\mathcal{D}^{\leq}+\mathcal{D}^{>}$. The operator $\mathcal{D}$ maps $C_G$ into itself, because it is the derivative of an equivariant map at an equivariant point. The truncation $M\mapsto M^{\leq}$ acts on each entry by a Fourier multiplier which is invariant under $l\mapsto(\pm l_1,\pm l_2)$. It commutes with the translations and reflections of $\mathbb{T}^2$ and with the conjugation by constant matrices. Hence $\mathcal{D}^{\leq}$ and $\mathcal{D}^{>}$ also map $C_G$ into itself. Let $\theta_i>\Theta'$. Since $|\mathcal{D}_{i'i}|\leq|\mathcal{D}^{\leq}_{i'i}|+|\mathcal{D}^{>}_{i'i}|$, it suffices to bound the sum in \eqref{eq:columns} which defines $c_i$ with $\mathcal{D}$ replaced by $\mathcal{D}^{\leq}$ and by $\mathcal{D}^{>}$. Since $k\mapsto\theta_k$ is a norm on $\R^2$, the field $\mathcal{D}^{\leq}e_i$ is supported on frequencies $k$ with $\theta_k\geq\theta_i-\Theta_0>\Theta$, so that the first of these sums equals $\|\bar\Delta^{-1}\mathcal{D}^{\leq}e_i\|_{\theta}/\omega_i$. By Lemma~\ref{lem:algebra},
\begin{align*}
\frac{\|\bar\Delta^{-1}\mathcal{D}^{\leq}e_i\|_{\theta}}{\omega_i}\leq\frac{|\mathcal{D}^{\leq}e_i|_{1}}{(\theta_i-\Theta_0)\theta_i|e_i|_{1}}
&\leq\frac{\|M_0\|+|k_1(i)|\|M_1\|+|k_2(i)|\|M_2\|}{(\theta_i-\Theta_0)\theta_i}\\
&\leq\frac{1}{\Theta}\Big[\frac{\|M_0\|}{\Theta'}+\|M_1\|+\tau\|M_2\|\Big],
\end{align*}
where we used $|k_1(i)|\leq\theta_i$ and $|k_2(i)|\leq\tau\theta_i$. Since $\beta_{i'}\leq\alpha_{\max}|e_{i'}|_{1}$ for all $i'\in\mathcal{I}$, as in the proof of Lemma~\ref{lem:A}, the second sum satisfies
\begin{align*}
\frac{1}{\omega_i}\sum_{i'\in\mathcal{I}}\beta_{i'}|\mathcal{D}^{>}_{i'i}|\leq\frac{\alpha_{\max}|\mathcal{D}^{>}e_i|_{1}}{\theta_i|e_i|_{1}}\leq\alpha_{\max}\frac{m_0+|k_1(i)|m_1+|k_2(i)|m_2}{\theta_i}\leq\alpha_{\max}\Big[\frac{m_0}{\Theta'}+m_1+\tau m_2\Big].
\end{align*}
The sum of the two bounds is $Z_1^{\mathrm{far}}$.
\end{proof}

For the bound $Z_2$ we put
\begin{align*}
b_0=|\bar X|_{1},\qquad b_1=|\bar X_x|_{1},\qquad b_2=|\bar X_t|_{1},
\end{align*}
and, for $r\geq0$,
\begin{equation}\label{eq:Z2-def}
\begin{aligned}
\Gamma_2(r)&=\frac34\Big[(b_2+\tau r)^2+2\tau(b_0+r)(b_2+\tau r)\Big],\\
Z_2(r)&=\alpha_{\max}\Big[3b_0+3b_1+3\bar\sigma^{1/2}b_2+\frac92r\Big]+\frac{\alpha_{\max}\Gamma_2(r)+\alpha_{tt}}{\vartheta}+\alpha_\eta\Big[4+\frac{2\tau}{\vartheta}(b_2+\tau r)\Big].
\end{aligned}
\end{equation}

\begin{lem}\label{lem:Z2}
Let $r>0$ and $h\in\mathcal{X}$ with $\|h\|_{\mathcal{X}}\leq r$. Then $\|\mathcal{A}(D\mathcal{F}(\bar U+h)-D\mathcal{F}(\bar U))\|\leq Z_2(r)r$.
\end{lem}

\begin{proof}
Let $h=(s,\xi)$ and $h'=(s',\xi')$ with $\|h'\|_{\mathcal{X}}\leq1$. By \eqref{eq:N2} we have $|\xi|_{1}\leq r$, $|\xi_x|_{1}\leq r$, $|\xi_t|_{1}\leq\tau r$, $|s|\leq r/\vartheta$, and $\|\xi'\|_{\theta}+\vartheta|s'|\leq1$. A direct computation from \eqref{eq:F} and \eqref{eq:N12} shows that the first component of $(D\mathcal{F}(\bar U+h)-D\mathcal{F}(\bar U))h'$ is
\begin{equation}\label{eq:Z2-eta}
2\langle \xi_x\cdot\xi'_x\rangle-2\bar\sigma\langle \xi_t\cdot\xi'_t\rangle-2s\langle (\bar X_t+\xi_t)\cdot\xi'_t\rangle-s'\langle 2\bar X_t\cdot\xi_t+|\xi_t|^2\rangle,
\end{equation}
and that the second component is
\begin{equation}\label{eq:Z2-N}
\begin{aligned}
&-\big[D\mathcal{N}_1(\bar X+\xi)-D\mathcal{N}_1(\bar X)\big]\xi'-\bar\sigma\big[D\mathcal{N}_2(\bar X+\xi)-D\mathcal{N}_2(\bar X)\big]\xi'\\
&\qquad-sD\mathcal{N}_2(\bar X+\xi)\xi'-s'\big[\mathcal{N}_2(\bar X+\xi)-\mathcal{N}_2(\bar X)\big]+s\xi'_{tt}+s'\xi_{tt}.
\end{aligned}
\end{equation}
We bound the terms of \eqref{eq:Z2-N} except the last two in the norm $|\cdot|_{1}$. Since
\begin{align*}
    D\mathcal{N}_1(X)\xi'=\mathcal{C}(\xi',X_x,X_x)+\mathcal{C}(X,\xi'_x,X_x)+\mathcal{C}(X,X_x,\xi'_x),
\end{align*}
the difference $[D\mathcal{N}_1(\bar X+\xi)-D\mathcal{N}_1(\bar X)]\xi'$ is a sum of nine values of $\mathcal{C}$ in which at least one argument is $\xi$ or $\xi_x$. By \eqref{eq:cubic} and \eqref{eq:N2},
\begin{align*}
|\big[D\mathcal{N}_1(\bar X+\xi)-D\mathcal{N}_1(\bar X)\big]\xi'|_{1}\leq\frac34r\,(4b_1+2b_0+3r)\|\xi'\|_{\theta}.
\end{align*}
In the same way, using $\bar\sigma\tau=\bar\sigma^{1/2}$ and $\bar\sigma\tau^2=1$,
\begin{align*}
\bar\sigma|\big[D\mathcal{N}_2(\bar X+\xi)-D\mathcal{N}_2(\bar X)\big]\xi'|_{1}\leq\frac34r\,(4\bar\sigma^{1/2}b_2+2b_0+3r)\|\xi'\|_{\theta}.
\end{align*}
Furthermore, we have
\begin{align*}
    |D\mathcal{N}_2(\bar X+\xi'')\xi'|_{1}\leq\Gamma_2(r)\|\xi'\|_{\theta}
\end{align*}
for every $\xi''$ with $\|\xi''\|_{\theta}\leq r$. This bounds the third term of \eqref{eq:Z2-N} by $(r/\vartheta)\Gamma_2(r)\|\xi'\|_{\theta}$ and, after integration along the segment from $\bar X$ to $\bar X+\xi$, the fourth term by $|s'|\Gamma_2(r)r$. The sum of these four bounds is at most
\begin{align*}
&r\Big[3b_0+3b_1+3\bar\sigma^{1/2}b_2+\frac{9r}2\Big]\|\xi'\|_{\theta}+\frac{r\Gamma_2(r)}{\vartheta}(\|\xi'\|_{\theta}+\vartheta|s'|)\\
&\qquad\leq r\Big[3b_0+3b_1+3\bar\sigma^{1/2}b_2+\frac{9r}2+\frac{\Gamma_2(r)}{\vartheta}\Big],
\end{align*}
and by Lemma~\ref{lem:A} the operator $\mathcal{A}$ maps the corresponding part of \eqref{eq:Z2-N} to an element of norm at most $\alpha_{\max}$ times this quantity. By Lemma~\ref{lem:A} again,
\begin{align*}
\|\mathcal{A}(0,s\xi'_{tt}+s'\xi_{tt})\|_{\mathcal{X}}\leq\alpha_{tt}(|s|\|\xi'\|_{\theta}+|s'|\|\xi\|_{\theta})\leq\frac{\alpha_{tt}r}{\vartheta}.
\end{align*}
Finally \eqref{eq:N1} and \eqref{eq:N2} bound the modulus of \eqref{eq:Z2-eta} by
\begin{align*}
r\Big[4+\frac{2\tau}{\vartheta}(b_2+\tau r)\Big]\|\xi'\|_{\theta}+|s'|\tau r(2b_2+\tau r)\leq r\Big[4+\frac{2\tau}{\vartheta}(b_2+\tau r)\Big],
\end{align*}
and $\mathcal{A}$ maps the first component to an element of norm at most $\alpha_\eta$ times its modulus. Adding the three contributions gives $Z_2(r)r$.
\end{proof}

The function $Z_2$ is non-decreasing. Lemmas~\ref{lem:C1}, \ref{lem:A}, \ref{lem:Z1}, and~\ref{lem:Z2} show that Theorem~\ref{thm:radii} applies to $\mathcal{F}$, $\bar U$, $\mathcal{A}$, and $\mathcal{A}^\dagger$ with the constants \eqref{eq:Y0-def}, \eqref{eq:Z1-def}, and \eqref{eq:Z2-def}, provided that $Z_0<1$. Upper bounds for these constants may be used in their place.

\section{The computation}\label{sec:computation}

We now evaluate the bounds of \S\ref{sec:functional} for an explicit approximate solution. The explicit inequalities used in the proof of Proposition~\ref{prop:existence} are collected in Proposition~\ref{prop:cap-existence}, and those used in \S\ref{sec:geometry} in Proposition~\ref{prop:cap-post}. Their proofs are computations in exact integer and rational arithmetic or in ball arithmetic, which we describe in this section; see also \S\ref{sec:assumptions}. \S\ref{sec:index-computation} and \S\ref{sec:rot-computation} describe the corresponding computations for the index and the rotational torus.

\begin{proposition}\label{prop:existence}
Let $\bar U=(\bar\sigma,\bar X)$ be the approximate solution described in \S\ref{sec:ubar}, let $\vartheta$ be as in \eqref{eq:vartheta}, and let
\begin{align*}
r_{\min}=8.1948870190792571\cdot10^{-9},\qquad r_{\max}=2.1897248302498472\cdot10^{-6}.
\end{align*}
Then the map $\mathcal{F}$ of \eqref{eq:F} has a zero $U_*=(\sigma_*,X_*)\in\mathcal{X}$ with $\|U_*-\bar U\|_{\mathcal{X}}\leq r_{\min}$. It is the only zero of $\mathcal{F}$ in the closed ball of radius $r_{\max}$ about $\bar U$ in $\mathcal{X}$, and $D\mathcal{F}(U_*)$ is injective on $\mathcal{X}$.
\end{proposition}

\subsection{The approximate solution}\label{sec:ubar}

The approximate solution is a pair $\bar U=(\bar\sigma,\bar X)$, where
\begin{align*}
\bar\sigma=1.1560357696393733
\end{align*}
is a binary64 floating-point number, written here in its shortest decimal form, and $\bar X$ is the element of $C_G$ whose coefficients in \eqref{eq:ansatz} are, for the indices in
\begin{align*}
K=\{i\in\mathcal{I} : \theta_i\leq\Theta\},\qquad \Theta=250,
\end{align*}
the binary64 numbers listed in the file \texttt{data/existence/ubar\_T250.json} of \cite{code}, and vanish for all other indices. The set $K$ contains $30516$ indices of the form $(W;m,j)$ and $15296$ indices of the form $(3;m,j)$, so that the truncated problem has $45813$ unknowns. The coefficients were obtained by Newton's method in floating-point arithmetic. They enter the proof only as explicit rational numbers. The weight of the coordinate $\sigma$ in the norm of $\mathcal{X}$ is the binary64 number
\begin{equation}\label{eq:vartheta}
\vartheta=42.730763029677775,
\end{equation}
which is an approximation of $\|\bar X\|_{\theta}/\bar\sigma$.

\begin{proposition}\label{prop:cap-existence}
Let $\bar U$ and $\vartheta$ be as above, and let $\mathcal{A}_K$ be the matrix of \S\ref{sec:finite-block}. Then
\begin{align*}
5.0974\cdot10^{-17}\leq\eta(\bar U)\leq5.0975\cdot10^{-17},
\end{align*}
the quantities of \S\ref{sec:functional} with $\Theta_0=200$ satisfy the upper bounds of Table~\ref{tab:bounds}, and $\Theta_+\geq250.00126$. Moreover, the polynomial $p$ of \eqref{eq:radii} satisfies
\begin{align*}
p(r_{\min})\leq-7.44720\cdot10^{-13},\qquad p(r_{\max})\leq-2.36634\cdot10^{-10}.
\end{align*}
\end{proposition}

\begin{proof}
The bounds are checked by the program \texttt{verify} of \cite{code}, in exact integer and rational arithmetic, by the computation described in \S\S\ref{sec:exact}--\ref{sec:radii-computation}.
\end{proof}

\subsection{Exact quantities}\label{sec:exact}

All coefficients of $\bar X$ and the number $\bar\sigma$ are dyadic rationals. The products of trigonometric polynomials with rational coefficients are computed exactly, by integer polynomial multiplication after a Kronecker substitution. In this way we obtain as exact rational numbers
\begin{itemize}
\item the number $\eta(\bar U)$ and the coefficients $\mathcal{N}(\bar U)_i$ for all $411361$ indices $i$ with $\theta_i\leq3\Theta$,
\item the Fourier coefficients of the $27$ entries of $M_0$, $M_1$, and $M_2$ in \eqref{eq:M}, and the coefficients of $\mathcal{N}_2(\bar X)$,
\item the first row of $D\mathcal{F}(\bar U)$ in \eqref{eq:DF},
\item the norms $b_0$, $b_1$, $b_2$, $|\bar X_{xx}|_{1}$, $|\bar X_{xt}|_{1}$, $\|M_p\|$, and $m_p$ for $\Theta_0=200$.
\end{itemize}
The index sets are determined by exact comparisons of the rational numbers $\theta_i^2$ and $\Theta^2$. Where the numbers $\theta_i$ themselves are needed, they are enclosed between dyadic rational numbers obtained from the integer square root. In particular, the exact value of $\eta(\bar U)$ gives the enclosure in Proposition~\ref{prop:cap-existence}.

\subsection{The finite block}\label{sec:finite-block}

The Fourier coefficients of the matrices \eqref{eq:M} are dyadic rational numbers, and we round each of them to the nearest integer multiple of $2^{-80}$. Every entry $\mathcal{D}_{i'i}$ is of the form $\frac1{16}\sum_lc_l\mu_l$ with at most $48$ terms, where the $\mu_l$ are Fourier coefficients of the matrices $M_p$ and the $c_l$ are integers. With the rounded coefficients in place of the $\mu_l$, the same sum is an integer multiple of $2^{-84}$, which we compute exactly in integer arithmetic and which differs from $\mathcal{D}_{i'i}$ by at most $2^{-85}\sum_l|c_l|$. The diagonal entries $-\theta_i^2$, the first row, and the column of $\sigma$ are computed exactly. We then round the entries of each column $c$ of the resulting matrix to integer multiples of $2^{-T_c}$, with $44\leq T_c\leq60$ chosen so that the corresponding integers are less than $2^{61}$ in modulus. This gives an integer matrix $\hat D$ and an explicit matrix $E$ with nonnegative entries such that
\begin{align*}
|D\mathcal{F}_K-\hat D\diag(2^{-T_c})|\leq E
\end{align*}
entrywise.

The matrix $\mathcal{A}_K$ is obtained by rounding a numerical inverse of $D\mathcal{F}_K$, computed in floating-point arithmetic, to dyadic rational numbers. The entries of the row $r$ are integer multiples of $2^{-S_r}$, with $57\leq S_r\leq75$ chosen so that the corresponding integers are at most $2^{60}$ in modulus. Thus $\mathcal{A}_K=\diag(2^{-S_r})\hat A$ with an integer matrix $\hat A$ of order $N=45813$. The integers $S_r$ and the entries of $\hat A$ are listed in the files \texttt{aprime\_S.npy} and \texttt{aprime\_Ahat.npy} of \cite{code}. We compute the product $\hat A\hat D$ exactly using FLINT’s integer matrix multiplication \cite{FL}, in $144$ blocks of order at most $4096$. Then
\begin{equation}\label{eq:Z0-enclosure}
|I-\mathcal{A}_KD\mathcal{F}_K|\leq|I-\mathcal{A}_K\hat D\diag(2^{-T_c})|+|\mathcal{A}_K|E,
\end{equation}
where the first matrix on the right is known exactly. In the weighted column sums \eqref{eq:opnorm} of the right side, we replace the entries of the first matrix by upper bounds which are integer multiples of powers of two, and the weights $\omega_i$ by dyadic rational upper bounds in the numerators and lower bounds in the denominators. The weighted column sums of $|\mathcal{A}_K|E$ are bounded by means of the row vector with the entries $\sum_{i'}\omega_{i'}|(\mathcal{A}_K)_{i'i}|$, which is computed exactly. This gives the bound $Z_0$ of Table~\ref{tab:bounds}.

Using dyadic bounds for the weights, we compute upper bounds for the weighted column sums in exact arithmetic. These give bounds for $\alpha_\eta$, $\alpha_i$, $\alpha_{\max}$, and $\alpha_{tt}$. The constant $Y_0$ of \eqref{eq:Y0-def} is evaluated in exact rational arithmetic from the exact residual. The first component and the indices in $K$ contribute at most $4.81377\cdot10^{-13}$, and the indices outside $K$ contribute at most $4.54219\cdot10^{-9}$.

\subsection{The bound \texorpdfstring{$Z_1$}{Z1}}

We take $\Theta_0=200$, so that $\Theta'=450$. The column sums $c_i$ of \eqref{eq:columns} are computed for the $45812$ indices $i\in K$ and for the $102393$ indices with $\Theta<\theta_i\leq\Theta'$. Since the matrices $M_p$ are supported in $\theta_l\leq2\Theta$, the rows $i'$ which occur satisfy $\theta_{i'}\leq\Theta'+2\Theta=950$, and there are $659856$ such rows. The entries $\mathcal{D}_{i'i}$ are computed in integer arithmetic from the rounded coefficients as in the finite block, and the error bound given there is added to their moduli. The column sums are accumulated exactly in integer arithmetic, after the weights $\beta_{i'}$ and the bounds for the moduli of the entries are rounded upward to dyadic grids. The constant $Z_1^{\mathrm{far}}$ of \eqref{eq:Z1-def} is evaluated from the exact norms $\|M_p\|$ and $m_p$ and from upper bounds for $\tau$ and $\alpha_{\max}$. The results are given in Table~\ref{tab:bounds}.

\subsection{The radii polynomial}\label{sec:radii-computation}

The function $Z_2$ of \eqref{eq:Z2-def} and the polynomial $p$ of \eqref{eq:radii} are evaluated in exact rational arithmetic from the certified upper bounds, whose decimal values rounded upward are listed in Table~\ref{tab:bounds}, with $\tau$ and $\bar\sigma^{1/2}$ replaced by rational upper bounds which are verified by squaring. This gives the bounds for $p(r_{\min})$ and $p(r_{\max})$ in Proposition~\ref{prop:cap-existence}.

\begin{table}[ht]
\centering
\small
\begin{tabular}{ll}
\toprule
quantity & certified upper bound\\
\midrule
$b_0$, $b_1$, $b_2$ & $11.3240$,\ \ $26.3869$,\ \ $32.8301$\\
$\|M_0\|$, $\|M_1\|$, $\|M_2\|$ & $116.331$,\ \ $43.7350$,\ \ $71.5742$\\
$m_0$ & $2.47859\cdot10^{-5}$\\
$m_1$ & $2.19918\cdot10^{-7}$\\
$m_2$ & $1.04079\cdot10^{-6}$\\
$\tau$, $\bar\sigma^{1/2}$ & $0.930068$,\ \ $1.07520$\\
$Y_0$ & $4.54267\cdot10^{-9}$\\
$Z_0$ & $3.35977\cdot10^{-11}$\\
$\alpha_\eta$ & $86.9276$\\
$\alpha_{\max}$, $\alpha_{tt}$ & $1010.29$,\ \ $696.160$\\
$c_\sigma$ & $8.73537\cdot10^{-11}$\\
$\max\{c_i : i\in K\}$ & $0.126576$\\
$\max\{c_i : \Theta<\theta_i\leq\Theta'\}$ & $0.285778$\\
$Z_1^{\mathrm{far}}$ & $0.443505$\\
$Z_1$ & $0.443505$\\
$Z_2(r_{\min})$, $Z_2(r_{\max})$ & $253142.49$,\ \ $253142.50$\\
\bottomrule
\end{tabular}
\medskip
\caption{Certified bounds for $\Theta=250$ and $\Theta_0=200$. In addition $\Theta_+\geq250.00126$.}
\label{tab:bounds}
\end{table}

\begin{proof}[Proof of Proposition~\ref{prop:existence}]
By Proposition~\ref{prop:cap-existence} and Table~\ref{tab:bounds} we have $Z_0<1$, so that $\mathcal{A}$ is injective by Lemma~\ref{lem:A}. By \S\ref{sec:functional} the hypotheses \eqref{eq:Y0}--\eqref{eq:Z2} of Theorem~\ref{thm:radii} hold, and by Proposition~\ref{prop:cap-existence} the inequality \eqref{eq:radii} holds for $r_0=r_{\min}$ and for $r_0=r_{\max}$. The theorem with $r_0=r_{\max}$ gives a unique zero $U_*$ in the ball of radius $r_{\max}$, and the theorem with $r_0=r_{\min}$ shows that $U_*$ lies in the ball of radius $r_{\min}$.
\end{proof}

\subsection{Quantities for the geometric conclusions}\label{sec:post}

The following quantities are used in \S\ref{sec:geometry}. The functions $|\bar X|$ and $|\bar X_x|$ are $G$-invariant, and $[0,\frac\pi3]\times[0,\pi]$ is a fundamental domain of the action of $G$ on $\mathbb{T}^2$. The values of $\bar X$ and $\bar X_x$ are enclosed on a grid of $719\times1712$ points of this domain, with spacings $h_x$ and $h_t$. Every point of the domain lies within $h_x/2$ and $h_t/2$ of a grid point in the two coordinates, and therefore
\begin{align*}
\min_{\mathbb{T}^2}|\bar X_x|\ \geq\ \min_{\text{grid}}|\bar X_x|-\frac12h_x|\bar X_{xx}|_{1}-\frac12h_t|\bar X_{xt}|_{1},
\end{align*}
with the analogous upper bound for $\sup|\bar X_x|$, and with $|\bar X_x|_1$ and $|\bar X_t|_1$ in place of $|\bar X_{xx}|_1$ and $|\bar X_{xt}|_1$ for $\sup|\bar X|$. The minimum over the grid is at least $0.409215$, and the sum of the last two terms is at most $0.299853$. The values at the grid points are enclosed in ball arithmetic, as products of matrices of balls which contain the cosines and sines at the grid points.

For $\sigma>0$ and $X\in C_G$ we put
\begin{equation}\label{eq:FJ}
\begin{aligned}
\mathcal{E}(\sigma,X)&=\frac{1}{4\pi\sqrt\sigma}\int_{\mathbb{T}^2}e^{-|X|^2/4}|X_x|^2dxdt,\\
\mathcal{J}_3(\sigma,X)&=\frac{1}{\sqrt\sigma}\int_{\mathbb{T}^2}\re\big((X_1+iX_2)^3\big)e^{-|X|^2/4}|X_x|^2dxdt.
\end{aligned}
\end{equation}
The integrals $\mathcal{E}(\bar U)$ and $\mathcal{J}_3(\bar U)$ are enclosed by the trapezoidal rule with $2304\times2048$ nodes. The integrands are $G$-invariant, and the sum is taken over the fundamental domain with the multiplicities of the orbits. The values of the integrands at the nodes are enclosed in ball arithmetic. The integrands extend to entire functions of $(x,t)$, where we write $\re((X_1+iX_2)^3)=X_1^3-3X_1X_2^2$. Let $f$ be one of them, let $a=\frac1{20}$, and let $M$ be an upper bound for $|f|$ on the set of points $(x,t)\in\mathbb{C}^2$ such that one of the two variables is real and the other has imaginary part at most $a$ in modulus. The estimate of \cite{TW} for periodic functions analytic in a strip, applied in each variable, shows that the trapezoidal rule with $N_x\times N_t$ nodes approximates $\int_{\mathbb{T}^2}f\,dxdt$ with an error at most
\begin{align*}
8\pi^2M\Big(\frac{1}{e^{aN_x}-1}+\frac{1}{e^{aN_t}-1}\Big).
\end{align*}
Let $B_0$ and $B_1$ be the norms $|\cdot|_1$ of $\bar X$ and $\bar X_x$ in which each Fourier coefficient of frequency $k$ is multiplied by $e^{a(|k_1|+|k_2|)}$. On the above set, the sums of the moduli of the components of $\bar X$ and $\bar X_x$ are at most $B_0$ and $B_1$. Since $|e^{-X\cdot X/4}|\leq e^{|X\cdot X|/4}$ and $|X_1^3-3X_1X_2^2|\leq(|X_1|+|X_2|)^3$, we may take $M=e^{B_0^2/4}B_1^2$ for $\mathcal{E}$ and $M=B_0^3e^{B_0^2/4}B_1^2$ for $\mathcal{J}_3$. We have $B_0\leq16.39954$ and $B_1\leq81.57345$, and with $N_x=2304$ and $N_t=2048$ the errors of the two integrals are at most $2.81205\cdot10^{-10}$ and $1.24027\cdot10^{-6}$. The results are given in Table~\ref{tab:post}.

\begin{table}[ht]
\centering
\begin{tabular}{ll}
\toprule
quantity & certified enclosure\\
\midrule
$\min_{\mathbb{T}^2}|\bar X_x|$ & $\geq0.109362$\\
$\sup_{\mathbb{T}^2}|\bar X|$ & $\leq4.63387$\\
$\sup_{\mathbb{T}^2}|\bar X_x|$ & $\leq10.06539$\\
$\mathcal{E}(\bar U)$ & $1.85257342255\pm2.2\cdot10^{-11}$\\
$\mathcal{J}_3(\bar U)$ & $-21.2233361\pm1.2\cdot10^{-6}$\\
\bottomrule
\end{tabular}
\medskip
\caption{Certified enclosures used in \S\ref{sec:geometry}.}
\label{tab:post}
\end{table}

\begin{proposition}\label{prop:cap-post}
The bounds of Table~\ref{tab:post} hold.
\end{proposition}

\begin{proof}
The bounds are checked by the program \texttt{verify} of \cite{code}, in ball arithmetic and in exact rational arithmetic, by the computation described above.
\end{proof}

\subsection{Rigorous numerics}\label{sec:assumptions}

All bounds are computed in exact integer and rational arithmetic or in ball arithmetic. Exact arithmetic uses the integers of Python and the integer and rational arithmetic of FLINT \cite{FL}, and ball arithmetic uses its component Arb \cite{Jo}, through python-flint version~0.9.0 with FLINT~3.6.0. Machine integers are used only for the integer approximations of the entries $\mathcal{D}_{i'i}$ in \S\ref{sec:finite-block} and in the computation of $Z_1$, and for the transfer of integer matrices to FLINT, and there the absence of overflow is verified by explicit bounds. Floating-point arithmetic is used only to compute the approximate solutions, the approximate inverses $\mathcal{A}_K$ before their rounding to dyadic rational numbers, the trigonometric polynomial $g$, the matrices $Y$, and the shifts $\nu_s$ of \S\ref{sec:index-computation}. These enter the proofs as explicit rational data, and no property of them other than their exact values is used.

The same arithmetic is used in \S\ref{sec:index-computation} and \S\ref{sec:rot-computation}.

\subsection{Code and data}\label{sec:code}

The programs and the data of the proofs, namely the approximate solutions $\bar U$ and $\bar U_{\mathrm{rot}}$, the matrices $\mathcal{A}_K$ and $\mathcal{A}_K^{\mathrm{rot}}$, and the trigonometric polynomial $g$, are contained in \cite{code}, together with the output of the computations. The program \texttt{verify} carries out all computations and checks every inequality stated in Propositions~\ref{prop:cap-existence}, \ref{prop:cap-post}, \ref{prop:cap-potential}, \ref{prop:cap-inertia}, and~\ref{prop:cap-rot}.

\section{Geometric conclusions}\label{sec:geometry}

To prove Theorem~\ref{thm:main}, we turn the zero of the conformal system into an embedded shrinking torus and identify its full symmetry group. Let $U_*=(\sigma_*,X_*)$ be the zero given by Proposition~\ref{prop:existence}, and put
\begin{align*}
r_0=r_{\min},\qquad \xi=X_*-\bar X,\qquad s=\sigma_*-\bar\sigma,\qquad \Sigma=X_*(\mathbb{T}^2).
\end{align*}
By Lemma~\ref{lem:algebra} and \eqref{eq:vartheta} we have, pointwise on $\mathbb{T}^2$ and for the Euclidean norm of $\R^3$,
\begin{equation}\label{eq:pointwise}
|\xi|\leq r_0,\qquad |\xi_x|\leq r_0,\qquad |s|\leq\frac{r_0}{\vartheta}\leq1.92\cdot10^{-10}.
\end{equation}

\subsection{Regularity and immersion}

We have $\sigma_*>0$ by \eqref{eq:pointwise}, and $X_*\in C_G$ is of class $C^1$ by Lemma~\ref{lem:algebra}. The equation $\mathcal{F}(U_*)=0$ in $\mathcal{Y}$ means that $\eta(X_*,\sigma_*)=0$ and that all coefficients in the nonlinearity $\mathcal{N}(X_*,\sigma_*)\in C_G$ vanish. By Proposition~\ref{prop:zeros} the map $X_*$ is real-analytic and weakly conformal. By Table~\ref{tab:post} and \eqref{eq:pointwise},
\begin{align*}
|X_{*,x}|\geq|\bar X_x|-|\xi_x|\geq0.1093-r_0>0
\end{align*}
on $\mathbb{T}^2$. Hence Proposition~\ref{prop:zeros} shows that $X_*$ is a real-analytic immersion of $\mathbb{T}^2$ which is conformal with respect to $dx^2+\sigma_*^{-1}dt^2$ and satisfies \eqref{eq:shrinker}. Its area element is $\sigma_*^{-1/2}|X_{*,x}|^2dxdt$. Consequently the Gaussian area of the immersion $X_*$, counted on the parameter domain, is $\mathcal{E}(\sigma_*,X_*)$ in the notation of \eqref{eq:FJ}.

\subsection{The Gaussian integrals}

\begin{lem}\label{lem:perturbation}
Let $U=(\sigma,X)$ and $\bar U=(\bar\sigma,\bar X)$ with $X,\bar X\in C_G$ of class $C^1$, and let $r_0>0$ and $0<\delta<\bar\sigma$ be such that $|X-\bar X|\leq r_0$ and $|X_x-\bar X_x|\leq r_0$ pointwise on $\mathbb{T}^2$ and $|\sigma-\bar\sigma|\leq\delta$. Let $a_0=\sup|\bar X|$ and $a_1=\sup|\bar X_x|$, and put
\begin{align*}
\kappa=\frac14(2a_0+r_0)(a_1+r_0)^2+2a_1+r_0.
\end{align*}
Then
\begin{align*}
|\mathcal E(U)-\mathcal{E}(\bar U)|&\leq\frac{\pi\kappa r_0}{(\bar\sigma-\delta)^{1/2}}+\frac{\delta\bar\sigma^{1/2}}{2(\bar\sigma-\delta)^{3/2}}\mathcal{E}(\bar U), \\
|\mathcal J_3(U)-\mathcal{J}_3(\bar U)|&\leq\frac{4\pi^2r_0}{(\bar\sigma-\delta)^{1/2}}\left(3(a_0+r_0)^2(a_1+r_0)^2+a_0^3\kappa\right)+\frac{\delta\bar\sigma^{1/2}}{2(\bar\sigma-\delta)^{3/2}}|\mathcal J_3(\bar U)|.
\end{align*}
\end{lem}

\begin{proof}
Let $w=e^{-|X|^2/4}$ and $\bar w=e^{-|\bar X|^2/4}$. Since $s\mapsto e^{-s/4}$ is $\frac14$-Lipschitz on $[0,\infty)$, the hypotheses give
\begin{align*}
|w-\bar w|\leq\frac14(2|\bar X|+r_0)r_0,    
\end{align*}
and moreover
\begin{align*}
||X_x|^2-|\bar X_x|^2|\leq(2|\bar X_x|+r_0)r_0.    
\end{align*}
Since $\bar w\leq1$, we obtain the pointwise bound
\begin{align*}
\Big|w|X_x|^2-\bar w|\bar X_x|^2\Big|\leq|w-\bar w||X_x|^2+\bar w\big||X_x|^2-|\bar X_x|^2\big|\leq\kappa r_0.
\end{align*}
With $W=X_1+iX_2$ and $W_0=\bar X_1+i\bar X_2$, we also have
\begin{align*}
    |\re(W^3)-\re(W_0^3)|\leq3(a_0+r_0)^2r_0,\qquad |\re(W_0^3)|\leq a_0^3,
\end{align*}
so the integrands of $\mathcal{J}_3$ differ by at most $(3(a_0+r_0)^2(a_1+r_0)^2+a_0^3\kappa)\,r_0$. Combining with
\begin{align*}
\big|\sigma^{-1/2}-\bar\sigma^{-1/2}\big|\leq\frac12(\bar\sigma-\delta)^{-3/2}\delta,
\end{align*}
the two estimates follow by the triangle inequality.
\end{proof}

By \eqref{eq:pointwise}, Lemma~\ref{lem:perturbation} applies to $U_*$ and $\bar U$ with $\delta=r_0/\vartheta$. Inserting the values of Proposition~\ref{prop:existence}, of \eqref{eq:vartheta}, and of Table~\ref{tab:post}, we find $\kappa\leq254.865$ and
\begin{equation}
|\mathcal{E}(U_*)-\mathcal{E}(\bar U)|\leq6.1028\cdot10^{-6},\qquad |\mathcal{J}_3(U_*)-\mathcal{J}_3(\bar U)|\leq9.5944\cdot10^{-3}.
\end{equation}
Together with Table~\ref{tab:post} this gives
\begin{equation}\label{eq:FI}
|\mathcal{E}(U_*)-1.8525734|\leq6.2\cdot10^{-6},\qquad |\mathcal{J}_3(U_*)+21.2233|\leq9.7\cdot10^{-3}.
\end{equation}
More precisely, the same bounds and Table~\ref{tab:post} give
\begin{equation}\label{eq:ball-enclosures}
\mathcal{E}(U_*)\in[1.8525673197,\,1.8525795254],\qquad \mathcal{J}_3(U_*)\in[-21.2329317,\,-21.2137405].
\end{equation}

\subsection{Embeddedness and entropy}

Since we have not yet established injectivity of the shrinker, we compute Gaussian area by integration on the parameter domain, possibly with multiplicity. The following standard argument, which we include for closed immersions, shows that Gaussian area below $2$ excludes both self-intersections and multiple coverings. The entropy identity follows from Section 7.2 in \cite{CM}, and the embeddedness criterion is stated in \cite{CIMW}.

\begin{lem}\label{lem:embedded}
Let $M$ be a closed surface and $X\colon M\to\R^3$ an immersion satisfying \eqref{eq:shrinker}. For $x_0\in\R^3$ and $t_0>0$, let
\begin{align*}
F_{x_0,t_0}(X)=\frac{1}{4\pi t_0}\int_Me^{-|X-x_0|^2/(4t_0)}\,dA,
\end{align*}
where $dA$ is the area element induced by $X$ on $M$. Then
\begin{align*}
F_{x_0,t_0}(X)\leq F_{0,1}(X)\qquad\text{for all }x_0\in\mathbb R^3,\,t_0>0.
\end{align*}
If additionally $F_{0,1}(X)<2$, then $X$ is an embedding.
\end{lem}

\begin{proof}
By \eqref{eq:shrinker}, the family of maps $X_t=(-t)^{1/2}X$ have normal velocity equal to their mean curvature vector, so they solve mean curvature flow up to a tangential term. Composing with a suitable time-dependent family of diffeomorphisms of $M$, we may remove the tangential part to obtain a smooth mean curvature flow of immersions on $(-\infty,0)$. Integrals over $M$ with respect to the induced measure are unchanged by this reparameterization, and therefore Huisken's monotonicity formula applies to the family $X_t$. This gives that for every $x_0\in\R^3$ and $s_0\in\R$, the Huisken Gaussian density defined as
\begin{align*}
\mathcal G(t)=\int_M\frac{-t}{4\pi(s_0-t)}\exp\Big(-\frac{|(-t)^{1/2}X-x_0|^2}{4(s_0-t)}\Big)dA
\end{align*}
is non-increasing for $t<\min(s_0,0)$. Taking $s_0=t_0-1$, we see that
\begin{align*}
\mathcal G(-1)=F_{x_0,t_0}(X),\qquad \lim_{t\to-\infty}\mathcal G(t)=F_{0,1}(X),
\end{align*}
since the integrand converges uniformly on $M$. This proves the first assertion.

Next assume that $X(p)=X(p')=x_0$ for two distinct points $p$ and $p'$ of $M$. We choose disjoint neighborhoods of $p$ and $p'$ on which $X$ is an embedding. The contribution of each neighborhood to $F_{x_0,t_0}(X)$ tends to $1$ as $t_0\to0$, hence $F_{0,1}(X)\geq2$. We conclude that $F_{0,1}(X)<2$ implies that $X$ is injective, and an injective immersion of a compact surface is an embedding.
\end{proof}

By \eqref{eq:FI}, the Gaussian area of the immersion $X_*$ is less than $2$. By Lemma~\ref{lem:embedded} the map $X_*\colon\mathbb{T}^2\to\R^3$ is an embedding. Hence $\Sigma$ is a real-analytic closed embedded surface of genus one which satisfies \eqref{eq:shrinker}, $F(\Sigma)=\mathcal{E}(U_*)$, and the first assertion of Lemma~\ref{lem:embedded} gives $\lambda(\Sigma)=F(\Sigma)$. Now \eqref{eq:FI} is the second assertion of Theorem~\ref{thm:main}.

The map $X_*$ is a conformal diffeomorphism from $\mathbb{T}^2$ with the metric $dx^2+\sigma_*^{-1}dt^2$ onto $\Sigma$. Hence $\Sigma$ is conformally equivalent to $\R^2/(2\pi\Z\times2\pi c\,\Z)$ with $c=\sigma_*^{-1/2}$. By \eqref{eq:pointwise} and the value of $\bar\sigma$ we have
\begin{align*}
    |c-\bar\sigma^{-1/2}|\leq\frac12(\bar\sigma-\delta)^{-3/2}\delta\leq7.8\cdot10^{-11}
\end{align*}
with $\delta=r_0/\vartheta$, and $|\bar\sigma^{-1/2}-0.93006727552|\leq7\cdot10^{-13}$. This proves the third assertion of Theorem~\ref{thm:main}.

\subsection{The symmetry group}

Let $\Gamma$ be the group of Euclidean isometries of $\R^3$ which map $\Sigma$ onto itself. Since $X_*\in C_G$, the set $\Sigma$ is invariant under $\rho$, $\iota_2$, and $\iota_3$, so that $G\subset\Gamma$.

\begin{lem}\label{lem:translation}
Let $\Sigma$ be a compact embedded self-shrinker and let $\varphi(X)=RX+v$, with $R\in O(3)$ and $v\in\R^3$, be an isometry of $\R^3$ with $\varphi(\Sigma)=\Sigma$. Then $v=0$.
\end{lem}

\begin{proof}
The mean curvature vector $\vec H=-Hn$ satisfies $\vec H=-\frac12X^\perp$ on $\Sigma$ and $\vec H(\varphi(p))=R\vec H(p)$ for $p\in\Sigma$. The normal space at $\varphi(p)$ is the image under $R$ of the normal space at $p$, so that the normal component of $\varphi(p)=Rp+v$ at $\varphi(p)$ is $R(p^\perp)+v^\perp$. It follows that $v^\perp=0$ at every point of $\varphi(\Sigma)=\Sigma$. Thus the constant vector $v$ is tangent to $\Sigma$ at every point. At a point where the function $\langle X,v\rangle$ attains its maximum on $\Sigma$ the vector $v$ is normal to $\Sigma$. Hence $v=0$.
\end{proof}

\begin{lem}\label{lem:groups}
Let $\Gamma\subset O(3)$ be a closed subgroup which contains the group $G$ of \eqref{eq:G} and does not contain $SO(3)$. If $\Gamma\ne G$, then $\Gamma$ contains a rotation about the $X_3$-axis by an angle which is not a multiple of $\frac{2\pi}{3}$.
\end{lem}

\begin{proof}
Assume first that $\Gamma$ is infinite. Then its identity component is a connected closed subgroup of $SO(3)$ of positive dimension different from $SO(3)$, and hence it is the group of all rotations about a line through the origin. This group is normalized by $\rho$, so that $\rho$ maps the line to itself, and the line is the $X_3$-axis. Thus $\Gamma$ contains all rotations about the $X_3$-axis.

Assume now that $\Gamma$ is finite and let $\Gamma^+=\Gamma\cap SO(3)$. The group $\Gamma^+$ contains $\rho$ and the rotation $\iota_2\iota_3$ by $\pi$ about the $X_1$-axis, which generate a dihedral group of order $6$. The finite subgroups of $SO(3)$ are the cyclic groups, the dihedral groups, and the tetrahedral, octahedral, and icosahedral groups. The group $\Gamma^+$ is not abelian, and it is not the tetrahedral group, which has no subgroup of order $6$. Suppose that $\Gamma^+$ is a dihedral group $D_k$ of order $2k$. Its elements of order three are rotations about its principal axis. Hence the principal axis is the $X_3$-axis and $k$ is a multiple of $3$. If $k>3$, then $\Gamma$ contains the rotation by $\frac{2\pi}{k}$ about the $X_3$-axis. If $k=3$, then $|\Gamma|\leq2|\Gamma^+|=12=|G|$ and $\Gamma=G$. Suppose finally that $\Gamma^+$ is the octahedral or the icosahedral group. The map $\gamma\mapsto\det(\gamma)\gamma$ is a homomorphism from $\Gamma$ to $SO(3)$ whose image is a finite group containing $\Gamma^+$. The octahedral and icosahedral groups are maximal among the finite subgroups of $SO(3)$, so that the image is $\Gamma^+$. Therefore $-\iota_3\in\Gamma^+$, and $-\iota_3$ is the rotation by $\pi$ about the $X_3$-axis.
\end{proof}

We can now prove the first assertion of Theorem~\ref{thm:main}. By Lemma~\ref{lem:translation} we have $\Gamma\subset O(3)$. The group $\Gamma$ is closed because $\Sigma$ is compact, and it does not contain $SO(3)$ because a torus is not a union of spheres centered at the origin. Let $R_\alpha\in\Gamma$ be a rotation by the angle $\alpha$ about the $X_3$-axis, and let
\begin{align*}
J=\int_\Sigma(X_1+iX_2)^3e^{-|X|^2/4}\,dA.
\end{align*}
Since $R_\alpha$ is an isometry of $\Sigma$ which multiplies $X_1+iX_2$ by $e^{i\alpha}$ and preserves $|X|$, we have $J=e^{3i\alpha}J$. The reflection $\iota_2\in\Gamma$ replaces $X_1+iX_2$ by its complex conjugate, so that $J$ is real, and $J=\mathcal{J}_3(U_*)$ by the formula for the area element. By \eqref{eq:FI} we have $J\ne0$. Hence $3\alpha\in2\pi\Z$, and Lemma~\ref{lem:groups} gives $\Gamma=G$.

If $\Sigma$ were a surface of revolution about a line $\ell$, then $\Gamma$ would contain all rotations about $\ell$, contradicting the finiteness of $\Gamma$. This proves the first assertion and completes the proof of Theorem~\ref{thm:main}. Corollary~\ref{cor:angenent} follows because the Angenent torus is a surface of revolution.

\subsection{Non-degeneracy in the symmetry class}

The injectivity of $D\mathcal{F}(U_*)$ in Proposition~\ref{prop:existence} has the following consequence.

\begin{proposition}\label{prop:nondegenerate}
The only $G$-invariant function $f$ on $\Sigma$ with $Lf=0$ is $f=0$.
\end{proposition}

\begin{proof}
Let $f$ be a $G$-invariant function on $\Sigma$ with $Lf=0$. Then $f$ is smooth, and $f\circ X_*$ is a smooth $G$-invariant function on $\mathbb{T}^2$. By the second part of Proposition~\ref{prop:jacobi}, there are $s'\in\R$ and a smooth $\xi'\in C_G$ with $\xi'\cdot n=f\circ X_*$ and $D\mathcal{F}(U_*)(s',\xi')=0$. A smooth map has finite norm $\|\cdot\|_{\theta}$, so that $(s',\xi')\in\mathcal{X}$. By Proposition~\ref{prop:existence} the operator $D\mathcal{F}(U_*)$ is injective on $\mathcal{X}$, and we conclude $\xi'=0$ and $f=0$.
\end{proof}

Thus the existence argument already gives non-degeneracy in the $G$-invariant class. The full spectral computation in \S\ref{sec:index} proves this again as part of Theorem~\ref{thm:index}, independently of Proposition~\ref{prop:nondegenerate}.

\section{The Morse index and the nullity}\label{sec:index}

We prove Theorem~\ref{thm:index} by reducing the Jacobi form to a Schr\"odinger form on the flat torus. From there we obtain rigorous enclosures of its potential, and bound its negative spectrum in each symmetry class. The same bounds give Corollary~\ref{cor:gap}, and the description of the kernel gives Corollary~\ref{cor:rigidity}. Let $U_*,r_0,\xi,s$, and $\Sigma$ be as in \S\ref{sec:geometry}. For a map $X\colon\mathbb{T}^2\to\R^3$ we set $\ell=|X_x|^2$, as well as $\ell_*=|X_{*,x}|^2$ and $\bar\ell=|\bar X_x|^2$ for the exact and approximate shrinkers. Since $X_*$ is conformal, we have $\ell_*=\sigma_*|X_{*,t}|^2$, and the area element of $\Sigma$ is $\sigma_*^{-1/2}\ell_* dxdt$.

\subsection{Reduction to a Schr\"odinger operator on the flat torus}

Let $\Sigma$ be a closed immersed self-shrinker and $w=e^{-|X|^2/4}$. Since $\nabla_\Sigma w=-\frac12wX^T$, the operator \eqref{eq:jacobi} is
\begin{align*}
Lf=-w^{-1}\dv_\Sigma(w\nabla_\Sigma f)-(|A|^2+1/2)f.
\end{align*}
It is self-adjoint in $L^2(\Sigma,wdA)$ and has compact resolvent, and its quadratic form is
\begin{equation}
\mathcal{B}_L(f)=\int_\Sigma\Big(|\nabla_\Sigma f|^2-(|A|^2+1/2)f^2\Big)wdA,\qquad f\in H^1(\Sigma).
\end{equation}

\begin{lem}\label{lem:ground}
Let $\Sigma$ be a closed immersed self-shrinker, $f\in H^1(\Sigma)$, and $\varphi=w^{1/2}f$. Then
\begin{align*}
\mathcal{B}_L(f)=\int_\Sigma\Big(|\nabla_\Sigma\varphi|^2+\Big(\frac{|X|^2}{16}+\frac{H^2}{4}-|A|^2-1\Big)\varphi^2\Big)\,dA.
\end{align*}
\end{lem}

\begin{proof}
We have $\nabla_\Sigma f=w^{-1/2}(\nabla_\Sigma\varphi+\frac14\varphi X^T)$, and hence
\begin{align*}
w|\nabla_\Sigma f|^2=|\nabla_\Sigma\varphi|^2+\frac14\langle X^T,\nabla_\Sigma\varphi^2\rangle+\frac1{16}|X^T|^2\varphi^2.
\end{align*}
By the divergence theorem the integral of the middle term over $\Sigma$ equals $-\frac14\int_\Sigma\varphi^2\dv_\Sigma X^T\,dA$. We have $\dv_\Sigma X=2$ and $\dv_\Sigma(\langle X,n\rangle n)=H\langle X,n\rangle$, so that \eqref{eq:shrinker} gives $\dv_\Sigma X^T=2-2H^2$ and $|X^T|^2=|X|^2-4H^2$. The coefficient of $\varphi^2$ is therefore
\begin{align*}
\frac1{16}(|X|^2-4H^2)-\frac14(2-2H^2)-|A|^2-\frac12=\frac{|X|^2}{16}+\frac{H^2}{4}-|A|^2-1.
\end{align*}
\end{proof}

We now return to the surface $\Sigma=X_*(\mathbb{T}^2)$.

\begin{lem}\label{lem:flat-form}
For $f\in H^1(\Sigma)$ let $\varphi=(w^{1/2}f)\circ X_*$. Then $\varphi\in H^1(\mathbb{T}^2)$ and
\begin{align*}
\mathcal{B}_L(f)=\sigma_*^{-1/2}\mathcal{B}_*(\varphi),\qquad \int_\Sigma f^2wdA=\sigma_*^{-1/2}\int_{\mathbb{T}^2}\ell_*\,\varphi^2dxdt,
\end{align*}
where
\begin{equation}\label{eq:form-star}
\mathcal{B}_*(\varphi)=\int_{\mathbb{T}^2}(\varphi_x^2+\sigma_*\varphi_t^2+V_*\varphi^2) dxdt,\qquad V_*=\ell_*\,\Big(\frac{|X|^2}{16}+\frac{H^2}{4}-|A|^2-1\Big)\circ X_*.
\end{equation}
\end{lem}

\begin{proof}
The map $X_*$ is a conformal diffeomorphism of $\mathbb{T}^2$ with the metric $dx^2+dy^2$, $y=t/\sqrt{\sigma_*}$, onto $\Sigma$. The Dirichlet integral is conformally invariant in dimension two, so that
\begin{align*}
\int_\Sigma|\nabla_\Sigma\varphi|^2\,dA=\int(\varphi_x^2+\varphi_y^2) dx\,dy=\sigma_*^{-1/2}\int_{\mathbb{T}^2}(\varphi_x^2+\sigma_*\varphi_t^2) dxdt,
\end{align*}
where on the left $\varphi$ is regarded as a function on $\Sigma$. The remaining terms follow from $dA=\sigma_*^{-1/2}\ell_* dxdt$ and Lemma~\ref{lem:ground}.
\end{proof}

Let $\mathcal{S}_*=-\partial_x^2-\sigma_*\partial_t^2+V_*$ be the self-adjoint operator in $L^2(\mathbb{T}^2)$, with respect to the measure $dxdt$, whose quadratic form is $\mathcal{B}_*$. The potential $V_*$ is real-analytic. The eigenvalues of $\mathcal{S}_*$ are not those of $L$, because the two operators are self-adjoint with respect to different measures. We compare the numbers of their negative and zero eigenvalues by means of the quadratic forms. Let $\mathcal{B}$ be a closed symmetric form which is bounded below, with form domain $\mathcal{Q}$ in a real Hilbert space. We denote by $n_-(\mathcal{B})$ the supremum of the dimensions of the subspaces $\mathcal{K}\subset\mathcal{Q}$ with $\mathcal{B}(\varphi)<0$ for all $\varphi\in\mathcal{K}\setminus\{0\}$. If the operator $T$ of $\mathcal{B}$ has compact resolvent, then by the min--max principle \cite[Chapter~XIII]{RS} $n_-(\mathcal{B})$ is the number of negative eigenvalues of $T$ counted with multiplicity, and by the representation theorem for closed semibounded forms \cite[Chapter~VI, Section~2]{Ka} the kernel of $T$ is the set of all $\varphi\in\mathcal{Q}$ with $\mathcal{B}(\varphi,\psi)=0$ for every $\psi\in\mathcal{Q}$.
We also write $n_-(T)$ for $n_-(\mathcal{B})$. If two such forms on the same domain satisfy $\mathcal{B}\leq\mathcal{B}'$, then $n_-(\mathcal{B}')\leq n_-(\mathcal{B})$. For a bounded function $c$ we denote by $T-c$ the operator of the form $\mathcal{B}(\varphi)-\int c\varphi^2$.

\begin{proposition}\label{prop:inertia}
For every $\beta\in\R$ the number of eigenvalues of $L$ which are less than $\beta$, counted with multiplicity, equals $n_-(\mathcal{S}_*-\beta\ell_*)$, and
\begin{align*}
\dim\ker(L-\beta)=\dim\ker(\mathcal{S}_*-\beta\ell_*).
\end{align*}
In particular the number of negative eigenvalues of $L$ is $n_-(\mathcal{S}_*)$, and the map $f\mapsto(w^{1/2}f)\circ X_*$ is a bijection from the kernel of $L$ onto the kernel of $\mathcal{S}_*$.
\end{proposition}

\begin{proof}
The map $\Phi f=(w^{1/2}f)\circ X_*$ is a linear bijection from $H^1(\Sigma)$ onto $H^1(\mathbb{T}^2)$, because $X_*$ is a real-analytic diffeomorphism of $\mathbb{T}^2$ onto $\Sigma$ and $w$ is bounded from above and from below by positive constants on $\Sigma$. By Lemma~\ref{lem:flat-form},
\begin{align*}
\mathcal{B}_L(f)-\beta\int_\Sigma f^2wdA=\sigma_*^{-1/2}\Big(\mathcal{B}_*(\Phi f)-\beta\int_{\mathbb{T}^2}\ell_*\,(\Phi f)^2dxdt\Big).
\end{align*}
The left side is the quadratic form of $L-\beta$ in $L^2(\Sigma,wdA)$, and the bracket is the quadratic form of $\mathcal{S}_*-\beta\ell_*$ in $L^2(\mathbb{T}^2)$. Hence $\Phi$ maps the subspaces on which the first form is negative onto the subspaces on which the second form is negative, and it maps the kernel of $L-\beta$ onto the kernel of $\mathcal{S}_*-\beta\ell_*$.
\end{proof}

\subsection{The potential}

For a map $X\colon\mathbb{T}^2\to\R^3$ of class $C^2$ and $\sigma>0$ we define
\begin{equation}\label{eq:P-pieces}
\begin{aligned}
\mathcal{P}_0&=-\ell+\frac18\ell|X|^2-\frac1{16}\big((X\cdot X_x)^2+\sigma(X\cdot X_t)^2\big),\\
\mathcal{P}_2&=|X_{xx}|^2+2\sigma|X_{xt}|^2+\sigma^2|X_{tt}|^2,\\
p_1&=X_x\cdot X_{xx},\qquad p_2=X_x\cdot X_{xt}.
\end{aligned}
\end{equation}
These functions are polynomials in $\sigma$ and in the derivatives of $X$ up to the second order.

\begin{proposition}\label{prop:potential}
The potential $V_*$ of \eqref{eq:form-star} is given by
\begin{equation}\label{eq:V}
V_*=\mathcal{P}_0-\frac{\mathcal{P}_2}{\ell}+\frac{4(p_1^2+\sigma p_2^2)}{\ell^2},
\end{equation}
where the right side is evaluated at $(X,\sigma)=(X_*,\sigma_*)$.
\end{proposition}

\begin{proof}
We work in the coordinates $(x,y)$ with $y=t/\sqrt{\sigma_*}$ and write $X=X_*$, $\sigma=\sigma_*$. The vectors $X_x$ and $X_y$ are orthogonal and have length $\ell^{1/2}$. Hence
\begin{align*}
\ell|X^T|^2=(X\cdot X_x)^2+(X\cdot X_y)^2=(X\cdot X_x)^2+\sigma(X\cdot X_t)^2,
\end{align*}
and \eqref{eq:shrinker} gives
\begin{align*}
    \ell H^2=\frac14\ell\langle X,n\rangle^2=\frac14\big(\ell|X|^2-(X\cdot X_x)^2-\sigma(X\cdot X_t)^2\big).
\end{align*}
Therefore
\begin{align*}
\ell\Big(\frac{|X|^2}{16}+\frac{H^2}{4}-1\Big)=\mathcal{P}_0.
\end{align*}
It remains to show that
\begin{equation}\label{eq:A-conformal}
\ell|A|^2=\frac{|D^2X|^2}{\ell}-\frac{|\nabla\ell|^2}{\ell^2},
\end{equation}
where the Hessian and the gradient are taken in the coordinates $(x,y)$. Since $\hat g_{\alpha\beta}=\ell\delta_{\alpha\beta}$, we have $|A|^2=\ell^{-2}\sum_{\alpha,\beta}h_{\alpha\beta}^2$ and
\begin{align*}
h_{\alpha\beta}^2=|X_{\alpha\beta}|^2-\frac1\ell\sum_\gamma(X_{\alpha\beta}\cdot X_\gamma)^2,\qquad X_{\alpha\beta}\cdot X_\gamma=\frac12(\delta_{\beta\gamma}\partial_\alpha\ell+\delta_{\alpha\gamma}\partial_\beta\ell-\delta_{\alpha\beta}\partial_\gamma\ell).
\end{align*}
The sum of the squares of the last expression over $\alpha,\beta,\gamma\in\{1,2\}$ is $\frac14(6+2-2-2)|\nabla\ell|^2=|\nabla\ell|^2$. This proves \eqref{eq:A-conformal}, which holds for every conformal immersion. In the coordinates $(x,t)$ we have
\begin{align*}
    |D^2X|^2=\mathcal{P}_2,\qquad |\nabla\ell|^2=(\partial_x\ell)^2+\sigma(\partial_t\ell)^2=4(p_1^2+\sigma p_2^2).
\end{align*}
Together with \eqref{eq:A-conformal} and the formula for $\mathcal{P}_0$, this gives \eqref{eq:V}.
\end{proof}

\subsection{Symmetry sectors}

We decompose functions on $\Sigma$ into six symmetry classes as follows. For a function $f$ obeying
\begin{align*}
f\circ\rho=f,\qquad f\circ\iota_2=\epsilon_2f,\qquad f\circ\iota_3=\epsilon_3f\qquad\text{for some }\epsilon_2,\epsilon_3\in\{\pm\},
\end{align*}
we assign the type
\begin{equation}\label{eq:types}
A_1'\text{ if }(\epsilon_2,\epsilon_3)=(+,+),\qquad A_2'\text{ if }(-,+),\qquad A_2''\text{ if }(+,-),\qquad A_1''\text{ if }(-,-).
\end{equation}
Further, a function $f$ satisfying
\begin{align*}
    f+f\circ\rho+f\circ\rho^2=0
\end{align*}
has type $E'$ if $f\circ\iota_3=f$, and type $E''$ if $f\circ\iota_3=-f$. The $G$-invariant functions have type $A_1'$. Every function decomposes uniquely into these six types and, since $L$ commutes with the symmetries, each is preserved by $L$.

The group $G$ acts on the functions on $\mathbb{T}^2$ by composition with the maps $\rho$, $\iota_2$, and $\iota_3$ of \S\ref{sec:conformal}. We write $\rho^*\varphi=\varphi\circ\rho$, so that $(\rho^*\varphi)(x,t)=\varphi(x+\frac{2\pi}{3},t)$. For $k=(k_1,k_2)$ with integers $k_1,k_2\geq0$ we consider the four products
\begin{align*}
\cos(k_1x)\cos(k_2t),\qquad \sin(k_1x)\cos(k_2t),\qquad \cos(k_1x)\sin(k_2t),\qquad \sin(k_1x)\sin(k_2t),
\end{align*}
and we define eight closed subspaces of $L^2(\mathbb{T}^2)$, which we call sectors, as the closed linear spans of the functions listed in Table~\ref{tab:sector-def}. The sectors are mutually orthogonal, and their sum is $L^2(\mathbb{T}^2)$. The sectors $A_1'$, $A_2'$, $A_2''$, and $A_1''$ consist of the functions which are invariant under $\rho^*$ and whose parities in $x$ and in $t$ are those of \eqref{eq:types}. The space $E'=E'_{\mathrm{c}}\oplus E'_{\mathrm{s}}$ consists of the functions $\varphi$ with $\varphi+\rho^*\varphi+(\rho^*)^2\varphi=0$ which are even in $t$, and $E''=E''_{\mathrm{c}}\oplus E''_{\mathrm{s}}$ consists of those which are odd in $t$. Since $X_*$ intertwines the actions of $G$ on $\mathbb{T}^2$ and on $\Sigma$ and $w$ is $G$-invariant, a function $f$ on $\Sigma$ is of one of the types of \S\ref{sec:intro-index} if and only if $(w^{1/2}f)\circ X_*$ belongs to the sector of the same name.

\begin{table}[ht]
\centering
\begin{tabular}{llll}
\toprule
sector & basis functions & frequencies & eigenfunctions of $L$\\
\midrule
$A_1'$ & $\cos(k_1x)\cos(k_2t)$ & $3\mid k_1$,\ \ $k_1\geq0$,\ \ $k_2\geq0$ & $H$\\
$A_2'$ & $\sin(k_1x)\cos(k_2t)$ & $3\mid k_1$,\ \ $k_1\geq1$,\ \ $k_2\geq0$ & $\langle e_3\times X,n\rangle$\\
$A_2''$ & $\cos(k_1x)\sin(k_2t)$ & $3\mid k_1$,\ \ $k_1\geq0$,\ \ $k_2\geq1$ & $\langle e_3,n\rangle$\\
$A_1''$ & $\sin(k_1x)\sin(k_2t)$ & $3\mid k_1$,\ \ $k_1\geq1$,\ \ $k_2\geq1$ & \\
$E'_{\mathrm{c}}$ & $\cos(k_1x)\cos(k_2t)$ & $3\nmid k_1$,\ \ $k_1\geq1$,\ \ $k_2\geq0$ & $\langle e_1,n\rangle$\\
$E'_{\mathrm{s}}$ & $\sin(k_1x)\cos(k_2t)$ & $3\nmid k_1$,\ \ $k_1\geq1$,\ \ $k_2\geq0$ & $\langle e_2,n\rangle$\\
$E''_{\mathrm{c}}$ & $\cos(k_1x)\sin(k_2t)$ & $3\nmid k_1$,\ \ $k_1\geq1$,\ \ $k_2\geq1$ & $\langle e_2\times X,n\rangle$\\
$E''_{\mathrm{s}}$ & $\sin(k_1x)\sin(k_2t)$ & $3\nmid k_1$,\ \ $k_1\geq1$,\ \ $k_2\geq1$ & $\langle e_1\times X,n\rangle$\\
\bottomrule
\end{tabular}
\medskip
\caption{The eight sectors of $L^2(\mathbb{T}^2)$. The last column lists the functions $f$ among $H$, $\langle e_a,n\rangle$, and $\langle e_a\times X,n\rangle$ for which $(w^{1/2}f)\circ X_*$ belongs to the sector. Their eigenvalues are $-1$, $-\frac12$, and $0$.}
\label{tab:sector-def}
\end{table}

\begin{lem}\label{lem:sectors}
Let $\sigma>0$, let $V'$ be a bounded $G$-invariant function on $\mathbb{T}^2$, and let $T'=-\partial_x^2-\sigma\partial_t^2+V'$. Then each of the eight sectors is invariant under $T'$. The restrictions of $T'$ to $E'_{\mathrm{c}}$ and to $E'_{\mathrm{s}}$ are unitarily equivalent, and so are its restrictions to $E''_{\mathrm{c}}$ and to $E''_{\mathrm{s}}$. Consequently
\begin{equation}\label{eq:sector-sum}
n_-(T')=\sum_{s\in\{A_1',A_2',A_2'',A_1''\}}n_-(T'|_s)+2n_-(T'|_{E'_{\mathrm{c}}})+2n_-(T'|_{E''_{\mathrm{c}}}),
\end{equation}
and the same formula holds for the dimensions of the kernels.
\end{lem}

\begin{proof}
The operator $T'$ commutes with $\rho^*$ and with the reflections $x\mapsto-x$ and $t\mapsto-t$. The orthogonal projection onto each sector is a product of three commuting projections, namely $\frac13(1+\rho^*+(\rho^*)^2)$ or its complement, the projection onto the even or onto the odd functions of $x$, and the projection onto the even or onto the odd functions of $t$. Hence it commutes with $T'$, which is the first assertion. Let
\begin{align*}
\mathcal{U}=\frac{1}{\sqrt3}\big(\rho^*-(\rho^*)^{-1}\big).
\end{align*}
Then $\mathcal{U} e^{ik_1x}=i\chi(k_1)e^{ik_1x}$, where $\chi(k_1)$ equals $0$, $1$, or $-1$ when $k_1$ is congruent to $0$, $1$, or $2$ modulo $3$, and therefore
\begin{align*}
\mathcal{U}\cos(k_1x)=-\chi(k_1)\sin(k_1x),\qquad \mathcal{U}\sin(k_1x)=\chi(k_1)\cos(k_1x).
\end{align*}
The operator $\mathcal{U}$ is real, acts on the variable $x$ only, and commutes with $T'$. It maps the orthogonal basis of $E'_{\mathrm{c}}$ in Table~\ref{tab:sector-def} onto that of $E'_{\mathrm{s}}$ up to signs, and it preserves the norms of the basis functions. The same holds for $E''_{\mathrm{c}}$ and $E''_{\mathrm{s}}$. Thus $\mathcal{U}$ is a unitary map from $E'_{\mathrm{c}}$ onto $E'_{\mathrm{s}}$ and from $E''_{\mathrm{c}}$ onto $E''_{\mathrm{s}}$ which intertwines the restrictions of $T'$. Finally $n_-$ and the dimension of the kernel are additive over an orthogonal decomposition into reducing subspaces, which gives \eqref{eq:sector-sum}.
\end{proof}

The potential $V_*$ is $G$-invariant. Indeed, $\ell_*=\mu+\frac12P$ is $G$-invariant by \S\ref{sec:conformal}, and the functions $|X|^2$, $H^2$, and $|A|^2$ on $\Sigma$ are invariant under the isometries of $\R^3$ which preserve $\Sigma$. Hence Lemma~\ref{lem:sectors} applies to $\mathcal{S}_*$ and to $\mathcal{S}_*-\beta\ell_*$ for every $\beta\in\R$. Moreover Proposition~\ref{prop:inertia} holds, with the same proof, for the restriction of $L$ to the functions of a given type and the restriction of $\mathcal{S}_*-\beta\ell_*$ to the corresponding sector.

\begin{lem}\label{lem:rotfields}
For every $v\in\R^3$ the function $\langle v\times X,n\rangle$ on $\Sigma$ satisfies $L\langle v\times X,n\rangle=0$, and it does not vanish identically unless $v=0$. The functions
\begin{align*}
\varphi_a=(w^{1/2}\langle e_a\times X,n\rangle)\circ X_*,\qquad a\in\{1,2,3\},
\end{align*}
are linearly independent elements of the kernel of $\mathcal{S}_*$, and
\begin{align*}
\varphi_3\in A_2',\qquad \varphi_2\in E''_{\mathrm{c}},\qquad \varphi_1\in E''_{\mathrm{s}}.
\end{align*}
\end{lem}

\begin{proof}
Let $v\in\R^3$ and let $R_\varepsilon$ be the rotation $\exp(\varepsilon v\times\cdot\,)$. The maps $R_\varepsilon X_*$ parametrize the self-shrinkers $R_\varepsilon\Sigma$, so that $S(R_\varepsilon X_*)=0$ for all $\varepsilon$. As in the proof of Lemma~\ref{lem:linearization}, the derivative of $S$ at $X_*$ in a direction $\xi'$ is $L(\xi'\cdot n)$. With $\xi'=v\times X_*$ we obtain $L\langle v\times X,n\rangle=0$. Assume that $\langle v\times X,n\rangle$ vanishes identically on $\Sigma$ for some $v\ne0$. Then the vector field $v\times X$ is tangent to the closed embedded surface $\Sigma$, and its flow, which consists of the rotations about the line spanned by $v$, preserves $\Sigma$. This contradicts the first assertion of Theorem~\ref{thm:main}, by which the group of isometries of $\Sigma$ is finite. Hence the linear map $v\mapsto\langle v\times X,n\rangle$ is injective, and by Proposition~\ref{prop:inertia} the functions $\varphi_a$ are linearly independent elements of the kernel of $\mathcal{S}_*$.

It remains to determine the sectors. Let $W=X_{*,1}+iX_{*,2}$ and $N=n_1+in_2$. Since $X_*$ and $n$ belong to $C_G$, Lemma~\ref{lem:fourier} shows that $W$ and $N$ are sums of the functions $e^{i(3m+1)x}\cos(jt)$ with real coefficients, and that $X_{*,3}$ and $n_3$ are sums of the functions $\cos(3mx)\sin(jt)$. We have
\begin{align*}
\langle e_3\times X,n\rangle=X_1n_2-X_2n_1=\im\big(\overline{W}N\big),
\end{align*}
and $\overline WN$ is a sum of the functions $e^{3imx}\cos(jt)\cos(j't)$ with real coefficients. Its imaginary part is odd in $x$, even in $t$, and invariant under $\rho^*$, so that it belongs to $A_2'$. Next
\begin{align*}
\langle e_2\times X,n\rangle=X_3n_1-X_1n_3,\qquad \langle e_1\times X,n\rangle=X_2n_3-X_3n_2.
\end{align*}
The functions $X_1$ and $n_1$ are even in $x$, the functions $X_2$ and $n_2$ are odd in $x$, and all four are even in $t$ and contain only the frequencies $k_1$ which are not divisible by $3$. The functions $X_3$ and $n_3$ are even in $x$, odd in $t$, and contain only the frequencies $k_1$ which are divisible by $3$. Hence $\langle e_2\times X,n\rangle\circ X_*$ belongs to $E''_{\mathrm{c}}$ and $\langle e_1\times X,n\rangle\circ X_*$ belongs to $E''_{\mathrm{s}}$. The multiplication by the $G$-invariant function $e^{-|X_*|^2/8}$ preserves each sector.
\end{proof}

The same expansions show that $H=\frac12\langle X,n\rangle$ is $G$-invariant and that the functions $\langle e_a,n\rangle$ lie in the sectors given in Table~\ref{tab:sector-def}.

\subsection{Enclosure of the potential}\label{sec:enclosure}

Let $g$ be a $G$-invariant trigonometric polynomial on $\mathbb{T}^2$, which will be a numerical approximation of $1/\bar\ell$. We denote by $\bar{\mathcal{P}}_0$, $\bar{\mathcal{P}}_2$, $\bar p_1$, and $\bar p_2$ the functions \eqref{eq:P-pieces} for $(X,\sigma)=(\bar X,\bar\sigma)$. We fix $\Theta_V>0$ and put
\begin{equation}\label{eq:V0}
V_{\mathrm{pol}}=\bar{\mathcal{P}}_0-g\bar{\mathcal{P}}_2+4(g\bar p_1)^2+4\bar\sigma\,(g\bar p_2)^2,\qquad V_0=\sum_{\theta_k\leq\Theta_V}(V_{\mathrm{pol}})^{\wedge}_ke^{i(k_1x+k_2t)},
\end{equation}
where $(V_{\mathrm{pol}})^{\wedge}_k$ are the Fourier coefficients of $V_{\mathrm{pol}}$. Both are $G$-invariant trigonometric polynomials: the algebraic operations in \eqref{eq:V0} preserve the symmetry, as does the frequency truncation. We do not assume that $\bar X$ is conformal or satisfies the shrinker equation. Formula \eqref{eq:V0} simply evaluates the algebraic expression of \eqref{eq:V} at the approximate data, with $1/\bar\ell$ replaced by $g$.

For a vector field $Y$ on $\mathbb{T}^2$ we use, in addition to $|Y|_{1}$, the quantity $|Y|_{1,\infty}=\max_a|Y_a|_{1}$, for which $|Y\cdot Y'|_{1}\leq|Y|_{1,\infty}|Y'|_{1}$. By Proposition~\ref{prop:existence} and Lemma~\ref{lem:algebra},
\begin{equation}\label{eq:xi-bounds}
|\xi|_{1}\leq r_0,\qquad |\xi_x|_{1}\leq r_0,\qquad |\xi_t|_{1}\leq\tau r_0,\qquad |s|\leq\frac{r_0}{\vartheta},
\end{equation}
and we put $\sigma_\pm=\bar\sigma\pm r_0/\vartheta$. The second derivatives of $\xi$ are controlled by the equation.

\begin{lem}\label{lem:second}
Let
\begin{equation}
\begin{aligned}
m_{\mathcal{D}}&=\|M_0\|+\|M_1\|+\tau\|M_2\|,\qquad c_{\mathcal{D}}=\frac34(2b_0+2b_1+2\bar\sigma^{1/2}b_2+2r_0),\\
\varepsilon_\Delta&=|\mathcal{N}(\bar U)|_{1}+m_{\mathcal{D}}r_0+c_{\mathcal{D}}r_0^2+\frac{r_0}{\vartheta}\Big(|\mathcal{N}_2(\bar X)|_{1}+\Gamma_2(r_0)r_0+|\bar X_{tt}|_{1}\Big),
\end{aligned}
\end{equation}
with the notation of \S\ref{sec:functional}. Then
\begin{align*}
|\xi_{xx}|_{1}+\sigma_*|\xi_{tt}|_{1}\leq\varepsilon_\Delta,\qquad 2\sigma_*^{1/2}|\xi_{xt}|_{1}\leq\varepsilon_\Delta.
\end{align*}
\end{lem}

\begin{proof}
By \eqref{eq:N12} the equations $\mathcal{N}(X_*,\sigma_*)=0$ and $\mathcal{N}(\bar U)=\bar\Delta\bar X-\mathcal{N}_1(\bar X)-\bar\sigma\mathcal{N}_2(\bar X)$ give
\begin{align*}
\xi_{xx}+\sigma_*\xi_{tt}=-\mathcal{N}(\bar U)+\big[\mathcal{N}_1(X_*)-\mathcal{N}_1(\bar X)\big]+\bar\sigma\big[\mathcal{N}_2(X_*)-\mathcal{N}_2(\bar X)\big]+s\mathcal{N}_2(X_*)-s\bar X_{tt}.
\end{align*}
The sum of the two brackets equals $\mathcal{D}\xi$ plus a remainder, where $\mathcal{D}$ is the operator of \eqref{eq:DF}. By \eqref{eq:N1} and \eqref{eq:xi-bounds} we have $|\mathcal{D}\xi|_{1}\leq m_{\mathcal{D}}r_0$. The remainder consists of the values of the trilinear form $\mathcal{C}$ of \eqref{eq:cubic} in which at least two arguments are $\xi$ or its first derivatives, namely
\begin{align*}
\mathcal{C}(\xi,\xi_x,\bar X_x)+\mathcal{C}(\xi,\bar X_x,\xi_x)+\mathcal{C}(\bar X,\xi_x,\xi_x)+\mathcal{C}(\xi,\xi_x,\xi_x)
\end{align*}
and the four analogous terms with the derivatives in $t$ and the factor $\bar\sigma$. By \eqref{eq:cubic}, \eqref{eq:xi-bounds}, $\bar\sigma\tau=\bar\sigma^{1/2}$, and $\bar\sigma\tau^2=1$, its norm $|\cdot|_{1}$ is at most $c_{\mathcal{D}}r_0^2$. The proof of Lemma~\ref{lem:Z2} gives $|\mathcal{N}_2(X_*)|_{1}\leq|\mathcal{N}_2(\bar X)|_{1}+\Gamma_2(r_0)r_0$. Hence $|\xi_{xx}+\sigma_*\xi_{tt}|_{1}\leq\varepsilon_\Delta$. For each Cartesian component the Fourier coefficients of $\xi_{xx}+\sigma_*\xi_{tt}$ are $-(k_1^2+\sigma_*k_2^2)\hat\xi_k$, so that $|\xi_{xx}|_{1}+\sigma_*|\xi_{tt}|_{1}=|\xi_{xx}+\sigma_*\xi_{tt}|_{1}$, and the bound for $\xi_{xt}$ follows from $2\sigma_*^{1/2}|k_1k_2|\leq k_1^2+\sigma_*k_2^2$.
\end{proof}

\begin{lem}\label{lem:neumann}
Let
\begin{equation}\label{eq:rho}
\bar\varrho=|1-g\bar\ell|_{1},\qquad \varsigma_0=|g|_{1},\qquad \varsigma_1=|g\bar X_x|_{1,\infty},\qquad \varrho=\bar\varrho+2\varsigma_1r_0+\varsigma_0r_0^2,
\end{equation}
and assume that $\varrho<1$. Then there is a function $\upsilon$ on $\mathbb{T}^2$ with
\begin{align*}
\frac{1}{\ell_*}=g\,(1+\upsilon),\qquad |\upsilon|_{1}\leq\varrho':=\frac{\varrho}{1-\varrho}.
\end{align*}
\end{lem}

\begin{proof}
Since $\ell_*=\bar\ell+2\bar X_x\cdot\xi_x+|\xi_x|^2$, we have
\begin{align*}
1-g\ell_*=(1-g\bar\ell)-2(g\bar X_x)\cdot\xi_x-g|\xi_x|^2,
\end{align*}
and $|1-g\ell_*|_{1}\leq\varrho$ by \eqref{eq:xi-bounds}. Since $\varrho<1$, the series $\sum_{j\geq0}(1-g\ell_*)^j$ converges in the Wiener algebra, and its product with $g\ell_*$ equals $1$. Hence $1/\ell_*=g\sum_{j\geq0}(1-g\ell_*)^j$, and the assertion holds with $\upsilon=\sum_{j\geq1}(1-g\ell_*)^j$.
\end{proof}

The next estimate separates the error into the polynomial part of the potential, the second-derivative and squared-gradient terms, the reciprocal approximation, and the final Fourier truncation.

\begin{proposition}\label{prop:eps}
Assume that $\varrho<1$ in Lemma~\ref{lem:neumann}. Let
\begin{equation}\label{eq:exact-norms}
\begin{aligned}
\varsigma_{xx}&=|g\bar X_{xx}|_{1,\infty}, & \varsigma_{xt}&=|g\bar X_{xt}|_{1,\infty}, & \varsigma_{tt}&=|g\bar X_{tt}|_{1,\infty},\\
q_1&=|g\bar p_1|_{1}, & q_2&=|g\bar p_2|_{1}, & q_{\mathcal{P}}&=|g\bar{\mathcal{P}}_2|_{1},\\
q_{xt}&=|g|\bar X_{xt}|^2|_{1}, & q_{tt}&=|g|\bar X_{tt}|^2|_{1}, & t_V&=|V_{\mathrm{pol}}-V_0|_{1},
\end{aligned}
\end{equation}
and define
\begin{align*}
\mathcal{M}(y_0,y_1,y_2,y_3)&=y_1^2+\frac{3}{16}y_0^2y_1^2+\frac{1}{16}y_3y_0^2y_2^2,\\
\varepsilon_1&=\mathcal{M}(b_0+r_0,\,b_1+r_0,\,b_2+\tau r_0,\,\sigma_+)-\mathcal{M}(b_0,b_1,b_2,\bar\sigma),\\
\varepsilon_2&=2\varepsilon_\Delta\Big[\max(\varsigma_{xx},\sigma_+\varsigma_{tt})+\sigma_+^{1/2}\varsigma_{xt}\Big]+\frac32\varsigma_0\varepsilon_\Delta^2+\frac{r_0}{\vartheta}\Big[2q_{xt}+(\bar\sigma+\sigma_+)q_{tt}\Big],\\
\varkappa_1&=\varsigma_{xx}r_0+(\varsigma_1+\varsigma_0r_0)\varepsilon_\Delta,\qquad \varkappa_2=\varsigma_{xt}r_0+\frac{(\varsigma_1+\varsigma_0r_0)\varepsilon_\Delta}{2\sigma_-^{1/2}},\\
\varepsilon_3&=4\Big[2q_1\varkappa_1+\varkappa_1^2+\sigma_+(2q_2\varkappa_2+\varkappa_2^2)\Big]+\frac{4r_0}{\vartheta}q_2^2,\\
q_+&=4\Big[(q_1+\varkappa_1)^2+\sigma_+(q_2+\varkappa_2)^2\Big],\\
\varepsilon_4&=(q_{\mathcal{P}}+\varepsilon_2)\varrho'+q_+\,(2\varrho'+\varrho'^2).
\end{align*}
Then the function $V_*-V_0$ is $G$-invariant and satisfies
\begin{equation}\label{eq:epsV}
\|V_*-V_0\|_{L^\infty}\leq|V_*-V_0|_{1}\leq\varepsilon_V:=\varepsilon_1+\varepsilon_2+\varepsilon_3+\varepsilon_4+t_V.
\end{equation}
\end{proposition}

\begin{proof}
We denote by $\mathcal{P}_0^*$, $\mathcal{P}_2^*$, $p_1^*$, and $p_2^*$ the functions \eqref{eq:P-pieces} for $(X,\sigma)=(X_*,\sigma_*)$. By Proposition~\ref{prop:potential} and Lemma~\ref{lem:neumann},
\begin{align*}
V_*=\mathcal{P}_0^*-g\mathcal{P}_2^*\,(1+\upsilon)+4\Big[(gp_1^*)^2+\sigma_*(gp_2^*)^2\Big](1+\upsilon)^2,
\end{align*}
and therefore
\begin{equation}\label{eq:V-difference}
\begin{aligned}
V_*-V_{\mathrm{pol}}&=[\mathcal{P}_0^*-\bar{\mathcal{P}}_0]-g[\mathcal{P}_2^*-\bar{\mathcal{P}}_2]+4\big[(gp_1^*)^2-(g\bar p_1)^2\big]+4\big[\sigma_*(gp_2^*)^2-\bar\sigma(g\bar p_2)^2\big]\\
&\qquad-g\mathcal{P}_2^*\,\upsilon+4\Big[(gp_1^*)^2+\sigma_*(gp_2^*)^2\Big](2\upsilon+\upsilon^2).
\end{aligned}
\end{equation}
We bound the six terms in the norm $|\cdot|_{1}$.

The function $\mathcal{P}_0$ is a sum of four monomials in $X$, $X_x$, $X_t$, and $\sigma$. By \eqref{eq:N1} the norm of the increment of a product is at most the increment of the product of the norms of its factors, and by \eqref{eq:xi-bounds} the factors $X_*$, $X_{*,x}$, $X_{*,t}$, and $\sigma_*$ have norms at most $b_0+r_0$, $b_1+r_0$, $b_2+\tau r_0$, and $\sigma_+$. Hence $|\mathcal{P}_0^*-\bar{\mathcal{P}}_0|_{1}\leq\varepsilon_1$.

For the second term we write
\begin{align*}
\mathcal{P}_2^*-\bar{\mathcal{P}}_2&=2\bar X_{xx}\cdot\xi_{xx}+|\xi_{xx}|^2+2\sigma_*(2\bar X_{xt}\cdot\xi_{xt}+|\xi_{xt}|^2)+\sigma_*^2(2\bar X_{tt}\cdot\xi_{tt}+|\xi_{tt}|^2)\\
&\qquad+2s|\bar X_{xt}|^2+s\,(\sigma_*+\bar\sigma)|\bar X_{tt}|^2.
\end{align*}
After multiplication by $g$, the terms which are linear in the second derivatives of $\xi$ have norm at most
\begin{align*}
2\varsigma_{xx}|\xi_{xx}|_{1}+2\sigma_*\varsigma_{tt}\sigma_*|\xi_{tt}|_{1}+4\sigma_*\varsigma_{xt}|\xi_{xt}|_{1}\leq2\varepsilon_\Delta\Big[\max(\varsigma_{xx},\sigma_+\varsigma_{tt})+\sigma_+^{1/2}\varsigma_{xt}\Big]
\end{align*}
by Lemma~\ref{lem:second}. The quadratic terms have norm at most
\begin{align*}
\varsigma_0\Big[(|\xi_{xx}|_{1}+\sigma_*|\xi_{tt}|_{1})^2+2\sigma_*|\xi_{xt}|_{1}^2\Big]\leq\frac32\varsigma_0\varepsilon_\Delta^2,
\end{align*}
and the two terms with the factor $s$ have norm at most $(r_0/\vartheta)\,[2q_{xt}+(\bar\sigma+\sigma_+)q_{tt}]$. Hence $|g\,(\mathcal{P}_2^*-\bar{\mathcal{P}}_2)|_{1}\leq\varepsilon_2$.

For the third and the fourth term we note that
\begin{align*}
g\,(p_1^*-\bar p_1)=(g\bar X_{xx})\cdot\xi_x+(g\bar X_x)\cdot\xi_{xx}+g\xi_x\cdot\xi_{xx}
\end{align*}
has norm at most $\varkappa_1$, and that the analogous expression for $p_2$ with $\xi_{xt}$ in place of $\xi_{xx}$ has norm at most $\varkappa_2$, by Lemma~\ref{lem:second} and $\sigma_*\geq\sigma_-$. Since
\begin{align*}
(gp_1^*)^2-(g\bar p_1)^2=2(g\bar p_1)g\,(p_1^*-\bar p_1)+\big(g\,(p_1^*-\bar p_1)\big)^2
\end{align*}
and
\begin{align*}
\sigma_*(gp_2^*)^2-\bar\sigma(g\bar p_2)^2=\sigma_*\Big[(gp_2^*)^2-(g\bar p_2)^2\Big]+s\,(g\bar p_2)^2,
\end{align*}
the sum of the norms of the third and the fourth term is at most $\varepsilon_3$.

Finally $|g\mathcal{P}_2^*|_{1}\leq q_{\mathcal{P}}+\varepsilon_2$ and $|4(gp_1^*)^2+4\sigma_*(gp_2^*)^2|_{1}\leq q_+$, and $|\upsilon|_{1}\leq\varrho'$ by Lemma~\ref{lem:neumann}. Hence the last two terms of \eqref{eq:V-difference} have norm at most $\varepsilon_4$. Together with the definition of $t_V$ this proves \eqref{eq:epsV}, since the supremum of a function is bounded by its norm $|\cdot|_{1}$. The function $V_*-V_0$ is $G$-invariant because $V_*$ and $V_0$ are.
\end{proof}

\subsection{The counting lemma}

Let $V_0$ and $\varepsilon_V$ be as in \S\ref{sec:enclosure}. We compare $\mathcal{S}_*$ with the operator
\begin{align*}
\mathcal{S}_0=-\bar\Delta+V_0=-\partial_x^2-\bar\sigma\partial_t^2+V_0
\end{align*}
in $L^2(\mathbb{T}^2)$, whose coefficients are explicit. By \eqref{eq:xi-bounds} the number
\begin{equation}\label{eq:varpi}
\varpi=\frac{r_0}{\vartheta\bar\sigma}\qquad\text{satisfies}\qquad|\sigma_*-\bar\sigma|\leq\varpi\bar\sigma.
\end{equation}

\begin{lem}\label{lem:shift}
Assume that $\varpi<1$, and for $\nu\in\R$ put
\begin{equation}\label{eq:epshat}
\hat\varepsilon(\nu)=\frac{\varepsilon_V+\varpi\,(|V_0|_{1}+|\nu|)}{1-\varpi}.
\end{equation}
Then
\begin{align*}
n_-\big(\mathcal{S}_0-\nu+\hat\varepsilon(\nu)\big)\leq n_-(\mathcal{S}_*-\nu)\leq n_-\big(\mathcal{S}_0-\nu-\hat\varepsilon(\nu)\big),
\end{align*}
and the same inequalities hold for the restrictions of the three operators to each sector.
\end{lem}

\begin{proof}
In the sense of quadratic forms on $H^1(\mathbb{T}^2)$, we have
\begin{align*}
(1-\varpi)(-\bar\Delta)\leq-\partial_x^2-\sigma_*\partial_t^2\leq(1+\varpi)(-\bar\Delta)    
\end{align*}
by \eqref{eq:varpi}, and $|V_*-V_0|\leq\varepsilon_V$ by Proposition~\ref{prop:eps}. Hence
\begin{align*}
\mathcal{S}_*-\nu\geq(1-\varpi)(\mathcal{S}_0-\nu)+\varpi\,(V_0-\nu)-\varepsilon_V\geq(1-\varpi)\big[\mathcal{S}_0-\nu-\hat\varepsilon(\nu)\big],
\end{align*}
where we used $V_0-\nu\geq-|V_0|_{1}-|\nu|$. In the same way
\begin{align*}
\mathcal{S}_*-\nu\leq(1+\varpi)(\mathcal{S}_0-\nu)-\varpi\,(V_0-\nu)+\varepsilon_V\leq(1+\varpi)\big[\mathcal{S}_0-\nu+\hat\varepsilon(\nu)\big].
\end{align*}
The number $n_-$ is not changed by a positive factor and is monotone with respect to the order of the forms. This proves the two inequalities. The three operators are reduced by the sectors by Lemma~\ref{lem:sectors}, since $V_*$ and $V_0$ are $G$-invariant, and the argument applies to the restrictions.
\end{proof}

The next lemma extends the finite-dimensional Schur complement identity for inertia \cite{Ha} to a positive infinite-dimensional block. We give the form-domain argument because it also identifies the kernel.

\begin{lem}\label{lem:schur}
Let $\mathcal{H}$ be a real Hilbert space and $\mathcal{B}$ a densely defined closed symmetric form in $\mathcal{H}$ which is bounded below, with form domain $\mathcal{Q}$, and let $T$ be the self-adjoint operator associated with $\mathcal{B}$. Let $\Pi$ be an orthogonal projection of finite rank whose range is contained in the domain of $T$, and $\Pi^\perp=I-\Pi$. Assume further that there is $d>0$ with
\begin{align*}
\mathcal{B}(\varphi)\geq d\,\|\varphi\|^2\qquad\text{for all }\varphi\in\mathcal{Q}\text{ with }\Pi\varphi=0.
\end{align*}
Then there is a symmetric operator $Z$ on the range of $\Pi$ with
\begin{equation}\label{eq:schur}
0\leq Z\leq\frac1d\,(\Pi^\perp T\Pi)^*(\Pi^\perp T\Pi)\qquad\text{and}\qquad n_-(\mathcal{B})=n_-(\Pi T\Pi-Z).
\end{equation}
Moreover, the dimension of $\{\varphi\in\mathcal{Q}:\mathcal{B}(\varphi,\psi)=0\,\forall\psi\in\mathcal{Q}\}$ equals the dimension of $\ker(\Pi T\Pi-Z)$.
\end{lem}

\begin{proof}
Since the range of $\Pi$ is contained in $\mathcal{Q}$, the space $\mathcal{Q}_0=\{\varphi\in\mathcal{Q}:\Pi\varphi=0\}$ equals $\Pi^\perp\mathcal{Q}$. It is dense in $\Pi^\perp\mathcal{H}$ and closed with respect to the form norm. The restriction of $\mathcal{B}$ to $\mathcal{Q}_0$ is therefore a closed form in the Hilbert space $\Pi^\perp\mathcal{H}$ which is bounded below by $d$. Let $T_0$ be its operator. Then $T_0^{-1}$ is a bounded operator on $\Pi^\perp\mathcal{H}$ with $0<T_0^{-1}\leq d^{-1}$ and with range in $\mathcal{Q}_0$. For $p$ in the range of $\Pi$ we put $Bp=T_0^{-1}\Pi^\perp Tp\in\mathcal{Q}_0$. For every $\psi\in\mathcal{Q}_0$ we have
\begin{align*}
\mathcal{B}(\psi,Bp)=\langle \psi,\Pi^\perp Tp\rangle=\langle \psi,Tp\rangle=\mathcal{B}(\psi,p),
\end{align*}
that is, $\mathcal{B}(\psi,p-Bp)=0$. It follows that
\begin{align*}
\mathcal{B}(p-Bp+\psi)=\mathcal{B}(p-Bp)+\mathcal{B}(\psi)=\langle (\Pi T\Pi-Z)p,p\rangle+\mathcal{B}(\psi),\qquad Z=(\Pi^\perp T\Pi)^*\,T_0^{-1}(\Pi^\perp T\Pi),
\end{align*}
because $\mathcal{B}(Bp,p-Bp)=0$ and therefore $\mathcal{B}(p-Bp)=\mathcal{B}(p,p-Bp)=\langle Tp,p\rangle-\langle Tp,Bp\rangle$. The bounds on $T_0^{-1}$ give the first part of \eqref{eq:schur}. The map $(p,\psi)\mapsto p-Bp+\psi$ is a linear bijection from the direct sum of the range of $\Pi$ and $\mathcal{Q}_0$ onto $\mathcal{Q}$, with the inverse $\varphi\mapsto(\Pi\varphi,\Pi^\perp\varphi+B\Pi\varphi)$. In these coordinates $\mathcal{B}$ is the direct sum of the quadratic form of $\Pi T\Pi-Z$ and of a positive definite form. This gives the second part of \eqref{eq:schur} and the last assertion.
\end{proof}

We apply the lemma in a fixed sector. Let $\mathcal{H}_s\subset L^2(\mathbb{T}^2)$ be one of the eight sectors and let $\Omega>0$. We denote by $\Pi$ the orthogonal projection of $\mathcal{H}_s$ onto the span of the basis functions of the sector in Table~\ref{tab:sector-def} whose frequencies satisfy $\theta_k\leq\Omega$, and we put $\Pi^\perp=I-\Pi$ on $\mathcal{H}_s$. The multiplication by $V_0$ maps $\mathcal{H}_s$ into itself. For $\nu\in\R$ and $d>0$ we define the symmetric operators
\begin{equation}\label{eq:Rpm}
\mathcal{R}^+(\nu)=\Pi\,(\mathcal{S}_0-\nu)\Pi,\qquad \mathcal{V}=\Pi V_0\Pi^\perp V_0\Pi=\Pi V_0^2\Pi-(\Pi V_0\Pi)^2,\qquad \mathcal{R}^-_d(\nu)=\mathcal{R}^+(\nu)-\frac1d\mathcal{V}
\end{equation}
on the range of $\Pi$, which has finite dimension.

\begin{proposition}\label{prop:updown}
Let $\nu\in\R$ and let $0<d\leq\Omega^2-|V_0|_{1}-\nu$. Then
\begin{align*}
n_-\big(\mathcal{R}^+(\nu)\big)\leq n_-\big((\mathcal{S}_0-\nu)|_{\mathcal{H}_s}\big)\leq n_-\big(\mathcal{R}^-_d(\nu)\big).
\end{align*}
\end{proposition}

\begin{proof}
We apply Lemma~\ref{lem:schur} in $\mathcal{H}_s$ to the form of $\mathcal{S}_0-\nu$, whose domain is $H^1(\mathbb{T}^2)\cap\mathcal{H}_s$. The functions in the range of $\Pi^\perp$ contain only frequencies with $\theta_k>\Omega$. For such a function $\varphi$ in the form domain we have $\langle -\bar\Delta\varphi,\varphi\rangle\geq\Omega^2\|\varphi\|^2$, because $-\bar\Delta$ has the symbol $\theta_k^2$, and $V_0\geq-|V_0|_{1}$. Hence the hypothesis of the lemma holds with the given $d$. Since $\bar\Delta$ commutes with $\Pi$, we have $\Pi^\perp(\mathcal{S}_0-\nu)\Pi=\Pi^\perp V_0\Pi$ and $(\Pi^\perp V_0\Pi)^*(\Pi^\perp V_0\Pi)=\mathcal{V}$. Thus the operator $Z$ of the lemma satisfies $0\leq Z\leq d^{-1}\mathcal{V}$, and
\begin{align*}
\mathcal{R}^-_d(\nu)\leq\Pi\,(\mathcal{S}_0-\nu)\Pi-Z\leq\mathcal{R}^+(\nu).
\end{align*}
The assertion follows from \eqref{eq:schur} and the monotonicity of $n_-$.
\end{proof}

The first inequality is the Rayleigh--Ritz bound. To determine both the index and nullity, we apply these bounds at opposite shifts and use the known rotation fields to account for the eigenvalues between them. The next proposition gives the precise criterion.

\begin{proposition}\label{prop:count}
Let $\mathcal{H}_s$ be a sector, and let $m\geq0$ and $z\geq0$ be integers such that the kernel of the restriction of $\mathcal{S}_*$ to $\mathcal{H}_s$ has dimension at least $z$. Let $\delta>0$, let $\nu=\delta+\hat\varepsilon(\delta)$ with the function $\hat\varepsilon$ of \eqref{eq:epshat}, and let $0<d\leq\Omega^2-|V_0|_{1}-\nu$. Assume that
\begin{equation}\label{eq:tests}
n_-\big(\mathcal{R}^+(-\nu)\big)\geq m\qquad\text{and}\qquad n_-\big(\mathcal{R}^-_d(\nu)\big)\leq m+z.
\end{equation}
Then the restriction of $\mathcal{S}_*$ to $\mathcal{H}_s$ has exactly $m$ eigenvalues which are less than $-\delta$, counted with multiplicity. Moreover, its kernel has dimension $z$ with no non-zero eigenvalues in the interval $[-\delta,\delta)$. In particular, it has exactly $m$ negative eigenvalues.
\end{proposition}

\begin{proof}
All operators are restricted to $\mathcal{H}_s$. Since $\hat\varepsilon(-\delta)=\hat\varepsilon(\delta)$, Lemma~\ref{lem:shift} and Proposition~\ref{prop:updown} give
\begin{align*}
n_-(\mathcal{S}_*+\delta)\geq n_-(\mathcal{S}_0+\nu)\geq n_-\big(\mathcal{R}^+(-\nu)\big)\geq m
\end{align*}
and
\begin{align*}
n_-(\mathcal{S}_*-\delta)\leq n_-(\mathcal{S}_0-\nu)\leq n_-\big(\mathcal{R}^-_d(\nu)\big)\leq m+z.
\end{align*}
The left sides are the numbers of the eigenvalues of $\mathcal{S}_*$ which are less than $-\delta$ and less than $\delta$. The eigenvalue $0$ lies in $[-\delta,\delta)$ and has multiplicity at least $z$. Hence the number of eigenvalues in $[-\delta,\delta)$ is at least $z$ and at most $m+z-m=z$, and all the inequalities are equalities.
\end{proof}

\subsection{Certified eigenvalue counts}\label{sec:index-computation}

In this section we describe the rigorous evaluation of $\varepsilon_V$ and the verification of \eqref{eq:tests}.

Recall the approximation $g$ of $1/\bar\ell$ introduced in \S\ref{sec:enclosure}. We define $g$ as a trigonometric polynomial of the $9866$ functions
\begin{align*}
    \{\cos(k_1x)\cos(k_2t):k_1\in3\mathbb Z_{\geq0},\,k_2\geq0,\,\theta_k\leq200\},
\end{align*}
with binary64 coefficients obtained from a discrete Fourier transform of $1/\bar\ell$, which are listed in the file \texttt{data/index/g\_K200.json} of \cite{code}. The truncation in \eqref{eq:V0} is $\Theta_V=300$. With the exact polynomial arithmetic of \S\ref{sec:computation} we compute as exact rational numbers the norms $|\mathcal{N}(\bar U)|_{1}$, $|\mathcal{N}_2(\bar X)|_{1}$, and $|\bar X_{tt}|_{1}$, the quantities \eqref{eq:rho} and \eqref{eq:exact-norms}, all the Fourier coefficients of $V_{\mathrm{pol}}$ and of $V_0$, and the Fourier coefficients of $V_0^2$ which are needed below. The products are formed in full before the norms are taken, and no intermediate truncation is made. The polynomial $V_0$ has $87\,639$ nonzero coefficients. The quantities of Lemmas~\ref{lem:second} and~\ref{lem:neumann}, of Proposition~\ref{prop:eps}, and of \eqref{eq:varpi} are then evaluated in exact rational arithmetic from these norms, from Proposition~\ref{prop:existence}, from \eqref{eq:vartheta}, and from the upper bounds of Table~\ref{tab:bounds}, with rational bounds for the square roots which are verified by squaring. The results are given in Table~\ref{tab:eps}. In particular $\varrho<1$ and $\varpi<1$.

\begin{table}[ht]
\centering
\small
\begin{tabular}{ll}
\toprule
quantity & certified upper bound\\
\midrule
$|\mathcal{N}(\bar U)|_{1}$ & $1.18068\cdot10^{-6}$\\
$|\mathcal{N}_2(\bar X)|_{1}$, $|\bar X_{tt}|_{1}$ & $320.414$,\ \ $326.625$\\
$\bar\varrho$ & $1.64388\cdot10^{-10}$\\
$\varsigma_0$, $\varsigma_1$ & $6.05851$,\ \ $2.70759$\\
$\varsigma_{xx}$, $\varsigma_{xt}$, $\varsigma_{tt}$ & $4.15435$,\ \ $3.10536$,\ \ $4.85806$\\
$q_1$, $q_2$ & $5.51231$,\ \ $6.70840$\\
$q_{\mathcal{P}}$, $q_{xt}$, $q_{tt}$ & $261.629$,\ \ $20.4878$,\ \ $137.318$\\
$t_V$ & $1.31928\cdot10^{-9}$\\
$|V_0|_{1}$ & $106.451$\\
$\varepsilon_\Delta$ & $3.16201\cdot10^{-6}$\\
$\varrho$ & $4.45411\cdot10^{-8}$\\
$\varepsilon_1$ & $5.58071\cdot10^{-5}$\\
$\varepsilon_2$ & $5.66999\cdot10^{-5}$\\
$\varepsilon_3$ & $6.27666\cdot10^{-4}$\\
$\varepsilon_4$ & $4.10184\cdot10^{-5}$\\
$\varepsilon_V$ & $7.81192\cdot10^{-4}$\\
$\varpi$ & $1.65895\cdot10^{-10}$\\
\bottomrule
\end{tabular}
\medskip
\caption{Certified bounds for the enclosure of the potential, with $g$ supported in $\theta_k\leq200$ and $\Theta_V=300$.}
\label{tab:eps}
\end{table}

\begin{proposition}\label{prop:cap-potential}
For the trigonometric polynomial $g$ above and $\Theta_V=300$, the quantities \eqref{eq:rho} and \eqref{eq:exact-norms}, the norms $|V_0|_1$, $|\mathcal{N}(\bar U)|_1$, $|\mathcal{N}_2(\bar X)|_1$, and $|\bar X_{tt}|_1$, and the constants of Lemmas~\ref{lem:second} and~\ref{lem:neumann}, of Proposition~\ref{prop:eps}, and of \eqref{eq:varpi} satisfy the upper bounds of Table~\ref{tab:eps}. In particular $\varrho<1$ and $\varpi<1$.
\end{proposition}

\begin{proof}
The bounds are checked by the program \texttt{verify} of \cite{code}, in exact rational arithmetic, by the computation described above.
\end{proof}

In the sequel $\varepsilon_V$, $\varpi$, and $|V_0|_{1}$ denote the rational upper bounds computed by the program, of which Table~\ref{tab:eps} gives upward roundings, and with these numbers Lemma~\ref{lem:shift} and Propositions~\ref{prop:updown} and~\ref{prop:count} remain valid. We consider the six sectors $A_1'$, $A_2'$, $A_2''$, $A_1''$, $E'_{\mathrm{c}}$, and $E''_{\mathrm{c}}$, and we take
\begin{align*}
\Omega=48,\qquad d=2048
\end{align*}
in all of them. For each sector $s$ we choose a dyadic rational number $\nu_s>0$ and put
\begin{equation}
\delta_s=(1-\varpi)\nu_s-\varepsilon_V-\varpi|V_0|_{1},
\end{equation}
so that $\delta_s>0$ and $\nu_s=\delta_s+\hat\varepsilon(\delta_s)$ by \eqref{eq:epshat}. The numbers $\nu_s$ are given in Table~\ref{tab:sectors}. By Proposition~\ref{prop:cap-inertia} below, $\Omega^2-|V_0|_{1}-\nu_s>2193$ in all six sectors, so that the condition on $d$ in Propositions~\ref{prop:updown} and~\ref{prop:count} holds.

The operators \eqref{eq:Rpm} are represented by matrices as follows. Let $\psi_k(x,t)=c_x(k_1x)c_t(k_2t)$ be the basis functions of the sector with $\theta_k\leq\Omega$, where $c_x$ and $c_t$ are the cosine or the sine according to Table~\ref{tab:sector-def}, and let $\epsilon_x=1$ if $c_x$ is the cosine and $\epsilon_x=-1$ if $c_x$ is the sine, and likewise $\epsilon_t$. The index set is determined by exact comparisons of $k_1^2+\bar\sigma k_2^2$ with $\Omega^2$. The functions $\psi_k$ are orthogonal, and the numbers $\pi^{-2}\int_{\mathbb{T}^2}\psi_k^2dxdt$ are equal to $1$, $2$, or $4$ when none, one, or both of $k_1$ and $k_2$ vanish. We denote by $\Lambda$ the diagonal matrix with these entries. Let $\hat V_k$ be the Fourier coefficients of $V_0$, which are real, even in $k_1$ and in $k_2$, and vanish unless $3\mid k_1$. Then
\begin{equation}\label{eq:BV}
\begin{aligned}
\frac{1}{\pi^2}\int_{\mathbb{T}^2}\psi_kV_0\psi_{k'} dxdt&=\hat V_{(k_1-k_1',\,k_2-k_2')}+\epsilon_x\hat V_{(k_1+k_1',\,k_2-k_2')}\\
&\qquad+\epsilon_t\hat V_{(k_1-k_1',\,k_2+k_2')}+\epsilon_x\epsilon_t\hat V_{(k_1+k_1',\,k_2+k_2')}.
\end{aligned}
\end{equation}
We denote by $B_V$ the matrix with the entries \eqref{eq:BV}, and by $B_{V^2}$ the matrix which is obtained in the same way with the Fourier coefficients of $V_0^2$ in place of those of $V_0$. The matrices
\begin{equation}\label{eq:Apm}
A^+(\nu)=\Lambda\diag(k_1^2+\bar\sigma k_2^2-\nu)+B_V,\qquad A^-(\nu)=A^+(\nu)-\frac1d\,(B_{V^2}-B_V\Lambda^{-1}B_V)
\end{equation}
are the matrices of the quadratic forms of $\mathcal{R}^+(\nu)$ and of $\mathcal{R}^-_d(\nu)$ in the basis $(\psi_k)$, divided by $\pi^2$. By Sylvester's law of inertia, they have the same numbers of negative eigenvalues as these operators. Since $\bar\sigma$, $\nu_s$, and the Fourier coefficients of $V_0$ are dyadic rational numbers, so are all the entries of $A^+(\pm\nu_s)$, $B_V$, and $B_{V^2}$, and we compute them without rounding.

The inertia of the twelve matrices $A^+(-\nu_s)$ and $A^-(\nu_s)$ is determined by means of the following lemma.

\begin{lem}
Let $A$ be a real symmetric $N\times N$ matrix, $Y$ and $B$ real $N\times N$ matrices, and $F$ a matrix with nonnegative entries such that
\begin{align*}
    |(Y^TAY-B)_{ij}|\leq F_{ij}\qquad\text{for all }i,j.
\end{align*}
Assume further that
\begin{equation}\label{eq:dominance}
|B_{ii}|>\sum_{j\ne i}|B_{ij}|+\sum_{j}F_{ij}\qquad\text{for every } i.
\end{equation}
Then $A$ is non-singular, and the number of negative eigenvalues equals the number of indices $i$ with $B_{ii}<0$.
\end{lem}

\begin{proof}
Let $A'=Y^TAY$, which is symmetric since $A$ is. By \eqref{eq:dominance}, we have
\begin{align*}
    |A'_{ii}|\geq|B_{ii}|-F_{ii}>\sum_{j\ne i}|A'_{ij}|,
\end{align*}
and moreover $A'_{ii}$ has the sign of $B_{ii}$. The symmetric matrices which are obtained from $A'$ by multiplying its off-diagonal entries by a factor in $[0,1]$ are strictly diagonally dominant, hence non-singular, and their eigenvalues depend continuously on this factor. Therefore the number of negative eigenvalues of $A'$ is the number of its negative diagonal entries. Since $A'$ is non-singular, the matrices $Y$ and $A$ are non-singular, and $A$ and $A'$ have the same number of negative eigenvalues by Sylvester's law of inertia.
\end{proof}

We choose $Y$ from the generalized eigenvectors of the pencil $(A^{\mathrm{mid}},\Lambda)$, computed in floating-point arithmetic, where $A^{\mathrm{mid}}$ is a binary64 approximation of $A$, so that $Y^TAY$ is expected to be nearly diagonal. We apply the lemma with $B=Y^TAY$ and $F=0$. The exact rational matrices $A^+(-\nu_s)$, $B_V$, $B_V\Lambda^{-1}$, and $B_{V^2}$ are enclosed in matrices of balls with $200$ bits, the matrices $A^-(\nu_s)$ and the products $Y^TAY$ are evaluated in ball arithmetic, and \eqref{eq:dominance} is decided from the endpoints of the resulting balls. In every row of every matrix the left side of \eqref{eq:dominance} exceeds the right side by at least $0.132$, and the smallest value occurs in the sector $E'_{\mathrm{c}}$.

\begin{table}[ht]
\centering
\small
\begin{tabular}{lrllcccc}
\toprule
sector & order & $\nu_s$ & lower bound for $\delta_s$ & $n_-(A^+(-\nu_s))$ & $n_-(A^-(\nu_s))$ & $m_s$ & $z_s$\\
\midrule
$A_1'$ & $589$ & $33/64$ & $0.514843$ & $3$ & $3$ & $3$ & $0$\\
$A_2'$ & $544$ & $31/8$ & $3.87421$ & $0$ & $1$ & $0$ & $1$\\
$A_2''$ & $572$ & $31/64$ & $0.483593$ & $1$ & $1$ & $1$ & $0$\\
$A_1''$ & $528$ & $4$ & $3.99921$ & $0$ & $0$ & $0$ & $0$\\
$E'_{\mathrm{c}}$ & $1140$ & $7/64$ & $0.108593$ & $3$ & $3$ & $3$ & $0$\\
$E''_{\mathrm{c}}$ & $1108$ & $23/16$ & $1.43671$ & $0$ & $1$ & $0$ & $1$\\
\bottomrule
\end{tabular}
\medskip
\caption{Certified numbers of negative eigenvalues of the matrices \eqref{eq:Apm} for $\Omega=48$ and $d=2048$, and the resulting numbers $m_s$ of negative eigenvalues and dimensions $z_s$ of the kernel of $\mathcal{S}_*$ in each sector. The sectors $E'_{\mathrm{s}}$ and $E''_{\mathrm{s}}$ are equivalent to $E'_{\mathrm{c}}$ and $E''_{\mathrm{c}}$.}
\label{tab:sectors}
\end{table}

\begin{proposition}\label{prop:cap-inertia}
Let $\Omega=48$ and $d=2048$, and for each of the six sectors $s$ of Table~\ref{tab:sectors} let $\nu_s$ be the number given there. Then $\Omega^2-|V_0|_1-\nu_s>2193$, the numbers $\delta_s$ satisfy the lower bounds of Table~\ref{tab:sectors}, and the matrices $A^+(-\nu_s)$ and $A^-(\nu_s)$ of \eqref{eq:Apm} are non-singular and have the numbers of negative eigenvalues listed in Table~\ref{tab:sectors}.
\end{proposition}

\begin{proof}
The first two assertions are inequalities between explicit rational numbers, by Proposition~\ref{prop:cap-potential}. For the third, the matrices $B=Y^TAY$ computed as described above satisfy \eqref{eq:dominance} with $F=0$, and the numbers of their negative diagonal entries are those of Table~\ref{tab:sectors}. The assertion follows from the lemma above. The computations are carried out by the program \texttt{verify} of \cite{code}.
\end{proof}

\subsection{Proof of Theorem~\ref{thm:index} and of its corollaries}

The proof gives the following more precise statement.

\begin{proposition}\label{prop:types}
The negative spectral subspace has dimension 3 in type $A_1'$, 1 in type $A_2''$, and 6 in type $E'$. The functions $\langle e_a\times X,n\rangle$, $a\in\{1,2,3\}$, which span the kernel of $L$, are of type $E''$, $E''$, and $A_2'$, respectively.
\end{proposition}

The four negative directions corresponding to dilation and translation are of type $A_1'$, $A_2''$, and $E'$. The remaining six consist of two directions of type $A_1'$ and four of type $E'$.

\begin{proof}[Proof of Theorem~\ref{thm:index} and Proposition~\ref{prop:types}]
Let $s$ be one of the six sectors of Table~\ref{tab:sectors}, and let $m_s$ and $z_s$ be the integers given there. By Lemma~\ref{lem:rotfields} the kernel of the restriction of $\mathcal{S}_*$ to $A_2'$ contains $\varphi_3$, and the kernel of its restriction to $E''_{\mathrm{c}}$ contains $\varphi_2$. Hence the kernel of the restriction of $\mathcal{S}_*$ to $s$ has dimension at least $z_s$ in all six sectors. By Proposition~\ref{prop:cap-inertia} and the remarks which follow \eqref{eq:Apm}, the inequalities \eqref{eq:tests} hold with $m=m_s$, $z=z_s$, $\delta=\delta_s$, $\nu=\nu_s$, $\Omega=48$, and $d=2048$. Proposition~\ref{prop:count} shows that the restriction of $\mathcal{S}_*$ to $s$ has exactly $m_s$ negative eigenvalues and a kernel of dimension $z_s$. By Lemma~\ref{lem:sectors},
\begin{align*}
n_-(\mathcal{S}_*)=3+0+1+0+2\cdot3+2\cdot0=10,\qquad \dim\ker\mathcal{S}_*=0+1+0+0+2\cdot0+2\cdot1=3.
\end{align*}
By Proposition~\ref{prop:inertia} the operator $L$ has exactly $10$ negative eigenvalues, and its kernel has dimension three. The three functions $\langle e_a\times X,n\rangle$ are linearly independent elements of this kernel by Lemma~\ref{lem:rotfields}, and therefore they span it. The statements about the types follow from Proposition~\ref{prop:inertia} for the restrictions, from Lemma~\ref{lem:sectors}, and from Lemma~\ref{lem:rotfields}. The third assertion is the case of the sector $A_1'$.
\end{proof}

\begin{remark}
The proof uses the first assertion of Theorem~\ref{thm:main} through Lemma~\ref{lem:rotfields}. Without the lower bound for the kernel, the counts of Table~\ref{tab:sectors} in the sectors $A_2'$ and $E''_{\mathrm{c}}$ show only that the restriction of $\mathcal{S}_*$ to each of these sectors has at most one eigenvalue which is less than $\delta_s$.
\end{remark}

\begin{proof}[Proof of Corollary~\ref{cor:gap}]
By Table~\ref{tab:post} and \eqref{eq:pointwise}, we have
\begin{align*}
\ell_{\max}:=\sup_{\mathbb{T}^2}\ell_*\leq\big(\sup_{\mathbb{T}^2}|\bar X_x|+r_0\big)^2\leq101.3121.
\end{align*}
Let $s$ be one of the six sectors of Table~\ref{tab:sectors}, and let us restrict all operators to $s$ for the remainder of the proof. By Proposition~\ref{prop:count}, the operator $\mathcal{S}_*$ has exactly $m_s$ eigenvalues less than $-\delta_s$, and its only eigenvalue in the interval $[-\delta_s,\delta_s)$ is $0$, which has multiplicity $z_s$. Consequently, for every constant $c$ we have
\begin{align*}
n_-(\mathcal{S}_*-c)=
\begin{cases}
m_s,&-\delta_s\leq c\leq 0\\
m_s+z_s,&0<c\leq\delta_s.
\end{cases}
\end{align*}
Now let $\beta\neq0$ with $|\beta|\leq\delta_s/\ell_{\max}$. The function $\beta\ell_*$ takes values between the constants $\beta\inf\ell_*$ and $\beta\ell_{\max}$, which have the same sign as $\beta$ and modulus at most $\delta_s$, so the monotonicity of $n_-$ gives
\begin{align*}
n_-(\mathcal{S}_*-\beta\ell_*)=
\begin{cases}
m_s,&-\delta_s/\ell_{\max}\leq\beta<0\\
m_s+z_s,&0<\beta\leq\delta_s/\ell_{\max}.
\end{cases}
\end{align*}
By Proposition~\ref{prop:inertia} for the restrictions, the left-hand side is the number of eigenvalues of $L$ less than $\beta$ whose eigenfunctions are of the type of $s$. Since this number does not change as $\beta$ varies in either interval, $L$ has no eigenvalue other than $0$ in $[-\delta_s/\ell_{\max},\delta_s/\ell_{\max})$ on the functions of this type, and by Lemma~\ref{lem:sectors} the same holds for the types of the sectors $E'_{\mathrm{s}}$ and $E''_{\mathrm{s}}$. Finally, Table~\ref{tab:sectors} gives $\delta_s\geq0.108593$ in all six sectors, and hence every nonzero eigenvalue of $L$ has modulus at least
\begin{align*}
\frac{0.108593}{101.3121}>1.07\cdot10^{-3}.
\end{align*}
\end{proof}

\begin{proof}[Proof of Corollary~\ref{cor:rigidity}]
For $u$ in a neighborhood of $0$ in $C^{2,\alpha}(\Sigma)$, the map $p\mapsto p+u(p)n(p)$ is an embedding of $\Sigma$, and we denote by $\Psi(u)\in C^{0,\alpha}(\Sigma)$ the function $H-\frac12\langle X,n\rangle$ of its image, as a function of $p$. The map $\Psi$ is continuously differentiable with $\Psi(0)=0$ and $D\Psi(0)=L$.

Let $\mathcal{K}$ be the kernel of $L$, which has dimension three by Theorem~\ref{thm:index}, and let $\mathcal{K}^\perp$ be its orthogonal complement in $L^2(\Sigma,wdA)$. Since $L$ is elliptic and symmetric with respect to this inner product, it is an isomorphism from $C^{2,\alpha}(\Sigma)\cap\mathcal{K}^\perp$ onto $C^{0,\alpha}(\Sigma)\cap\mathcal{K}^\perp$. Let $\Pi_{\mathcal{K}}^\perp$ be the projection of $C^{0,\alpha}(\Sigma)$ onto $C^{0,\alpha}(\Sigma)\cap\mathcal{K}^\perp$ along $\mathcal{K}$. By the implicit function theorem applied to the map $(k,v)\mapsto\Pi_{\mathcal{K}}^\perp\Psi(k+v)$ for $k\in\mathcal{K}$ and $v\in C^{2,\alpha}(\Sigma)\cap\mathcal{K}^\perp$, there is a $C^1$ map $k\mapsto v(k)$ in a neighborhood of $0$ in $C^{2,\alpha}(\Sigma)$ with $v(0)=0$ such that every zero $u$ of $\Psi$ in this neighborhood has the form $u=k+v(k)$ with $k\in\mathcal{K}$.

For any rotation $R\in SO(3)$ close to the identity, the surface $R\Sigma$ is a self-shrinker and can be written as a graph $u_R$ over $\Sigma$ in the above sense. Since $\Sigma$ is smooth, the map $R\mapsto u_R$ is $C^1$ with values in $C^{2,\alpha}(\Sigma)$ with derivative at the identity given by $v\mapsto\langle v\times X,n\rangle$ for $v\in\R^3$. Let $k_R\in\mathcal{K}$ be the orthogonal projection of $u_R$. By Lemma~\ref{lem:rotfields} and Theorem~\ref{thm:index}, the derivative of $R\mapsto k_R$ at the identity is an isomorphism from $\R^3$ onto $\mathcal{K}$. By the inverse function theorem, the set of functions $k_R$ contains a neighborhood $B$ of $0$ in $\mathcal{K}$. Since $\Psi(u_R)=0$, we have $u_R=k_R+v(k_R)$. We conclude that every zero of $\Psi$ in the form $u=k+v(k)$ with $k\in B$ equals $u_R$ for some rotation $R$, and the corresponding surface is $R\Sigma$.
\end{proof}

\section{A rotationally symmetric torus}\label{sec:rot}

To prove Theorem~\ref{thm:rot}, we apply the existence argument of \S\ref{sec:functional}--\ref{sec:computation} within the subspace of surfaces of revolution. This produces a compact embedded torus whose entropy can be enclosed by the same Gaussian-area estimates, yielding Corollary~\ref{cor:least}.

\subsection{The rotational subspace}\label{sec:rot-subspace}

For $\alpha\in\R$, let $R_\alpha$ denote rotation by the angle $\alpha$ about the $X_3$-axis, and define the space
\begin{align*}
C_{\mathrm{rot}}=\big\{X\in C_G : X(x+\alpha,t)=R_\alpha X(x,t)\text{ for all }\alpha\in\R\big\}.
\end{align*}
A map $X\in C_G$ belongs to $C_{\mathrm{rot}}$ if and only if
\begin{equation}\label{eq:rot-ansatz}
X_1+iX_2=e^{ix}a(t),\qquad X_3=b(t),
\end{equation}
where $a$ is a real even function and $b$ is a real odd function of period $2\pi$. In the expansion \eqref{eq:ansatz}, these are the maps whose coefficients $z_{m,j}$ and $\zeta_{m,j}$ vanish for $m\ne0$. Let $\mathcal{I}_{\mathrm{rot}}\subset\mathcal{I}$ be the set of indices of the form $(W;0,j)$ with $j\geq0$ and $(3;0,j)$ with $j\geq1$, whose frequencies are $(1,j)$ and $(0,j)$. We denote by $\mathcal{X}_{\mathrm{rot}}$ and $\mathcal{Y}_{\mathrm{rot}}$ the closed subspaces of $\mathcal{X}$ and $\mathcal{Y}$ which are spanned by the first coordinate and by the fields $e_i$ with $i\in\mathcal{I}_{\mathrm{rot}}$, so that $\mathcal{X}_{\mathrm{rot}}=\{(\sigma,X)\in\mathcal{X}:X\in C_{\mathrm{rot}}\}$. If $X\in C_{\mathrm{rot}}$ is of the form \eqref{eq:rot-ansatz} and $a$ has no zero, then $X(\mathbb{T}^2)$ is the surface of revolution about the $X_3$-axis which is generated by the closed curve $t\mapsto(|a(t)|,b(t))$ in the half-plane with the coordinates $r=(X_1^2+X_2^2)^{1/2}$ and $z=X_3$. Since $a$ is even and $b$ is odd, this curve and therefore the surface are symmetric under $\iota_3$.

\begin{lem}\label{lem:rot-system}
Let $X\in C_{\mathrm{rot}}$ be of class $C^2$ and of the form \eqref{eq:rot-ansatz}, and let $\sigma\in\R$. Then $\mathcal{N}(X,\sigma)=(N_a\cos x,\,N_a\sin x,\,N_b)$ with
\begin{align*}
N_a&=\sigma a''-a-\frac12\sigma\,(aa'+bb')a'+\frac14\big(a^2+\sigma(a'^2+b'^2)\big)a,\\
N_b&=\sigma b''-\frac12\sigma\,(aa'+bb')b'+\frac14\big(a^2+\sigma(a'^2+b'^2)\big)b,
\end{align*}
and
\begin{align*}
\eta(X,\sigma)=\frac1{2\pi}\int_0^{2\pi}(a^2-\sigma(a'^2+b'^2))dt.
\end{align*}
Moreover, $\mathcal{F}$ maps $\mathcal{X}_{\mathrm{rot}}$ into $\mathcal{Y}_{\mathrm{rot}}$, and hence so does $D\mathcal{F}(U)$ for every $U\in\mathcal{X}_{\mathrm{rot}}$.
\end{lem}

\begin{proof}
From \eqref{eq:rot-ansatz}, one easily computes
\begin{align*}
X\cdot X_x=0,\quad X\cdot X_t=aa'+bb',\quad |X_x|^2=a^2,\quad |X_t|^2=a'^2+b'^2.
\end{align*}
Inserting these into \eqref{eq:N} and \eqref{eq:F} gives the claimed formulas. The operator $\mathcal{N}(\cdot,\sigma)$ maps $C_G$ into itself and commutes with rotations and with the translations of $\mathbb{T}^2$. Hence $\mathcal{N}(X,\sigma)\in C_{\mathrm{rot}}$ for $X\in C_{\mathrm{rot}}$, which proves that $\mathcal{F}:\mathcal X_{\mathrm{rot}}\to \mathcal Y_{\mathrm{rot}}$. The claim for $D\mathcal{F}(U)$ follows by differentiation, since $\mathcal{Y}_{\mathrm{rot}}$ is closed.
\end{proof}

\subsection{The Newton--Kantorovich argument on the subspace}\label{sec:rot-computation}

We consider pairs $\bar U_{\mathrm{rot}}=(\bar\sigma_{\mathrm{rot}},\bar X_{\mathrm{rot}})$ with the binary64 number
\begin{align*}
\bar\sigma_{\mathrm{rot}}=1.1742314783959917,
\end{align*}
and with $\bar X_{\mathrm{rot}}\in C_{\mathrm{rot}}$ a trigonometric polynomial whose coefficients are dyadic rational numbers for the indices in
\begin{align*}
K_{\mathrm{rot}}=\{i\in\mathcal{I}_{\mathrm{rot}} : \theta_i\leq\Theta_{\mathrm{rot}}\},\qquad \Theta_{\mathrm{rot}}=500,
\end{align*}
and vanish for all other indices.
Throughout this section, we use the definitions of \S\ref{sec:functional} with $\bar U_{\mathrm{rot}}$, $\Theta_{\mathrm{rot}}$, and $\vartheta_{\mathrm{rot}}$ in place of $\bar U$, $\Theta$, and $\vartheta$, where we take
\begin{align*}
\vartheta_{\mathrm{rot}}=13.22662925925403.
\end{align*}
In particular, we have $\theta_k=(k_1^2+\bar\sigma_{\mathrm{rot}}k_2^2)^{1/2}$. Note that $\bar\sigma_{\mathrm{rot}}\geq1$, as required in \S\ref{sec:functional}.

The set $K_{\mathrm{rot}}$ consists of the indices $(W;0,j)$ for $0\leq j\leq461$ along with $(3;0,j)$ for $1\leq j\leq461$, so the truncated problem has $924$ unknowns.

We apply Theorem~\ref{thm:radii} to the restriction of $\mathcal{F}$ to $\mathcal{X}_{\mathrm{rot}}$, with the operators
\begin{align*}
\mathcal{A}_{\mathrm{rot}}=\mathcal{A}_K^{\mathrm{rot}}\Pi+\bar\Delta^{-1}\Pi_\infty,\qquad \mathcal{A}_{\mathrm{rot}}^\dagger=D\mathcal{F}_K^{\mathrm{rot}}\Pi+\bar\Delta\Pi_\infty
\end{align*}
from $\mathcal{Y}_{\mathrm{rot}}$ to $\mathcal{X}_{\mathrm{rot}}$ and from $\mathcal{X}_{\mathrm{rot}}$ to $\mathcal{Y}_{\mathrm{rot}}$. Here $D\mathcal{F}_K^{\mathrm{rot}}$ is the matrix of $\Pi D\mathcal{F}(\bar U_{\mathrm{rot}})\Pi$ with respect to the coordinate $\sigma$ and the fields $e_i$, $i\in K_{\mathrm{rot}}$, which is a $924\times924$ square matrix by Lemma~\ref{lem:rot-system}, and $\mathcal{A}_K^{\mathrm{rot}}$ is a binary64 matrix which approximates its inverse. The bounds of \S\ref{sec:functional} are valid in this setting for the following reasons.

Note that the space $\mathcal{X}_{\mathrm{rot}}$ is the weighted $\ell^1$ space over the index set $\{\sigma\}\cup\mathcal{I}_{\mathrm{rot}}$ with the same weights as $\mathcal{X}$. Hence \eqref{eq:norms-basis} holds on $C_{\mathrm{rot}}$, and the operator norm on $\mathcal{X}_{\mathrm{rot}}$ is given by \eqref{eq:opnorm} with $i$ and $i'$ in $\{\sigma\}\cup\mathcal{I}_{\mathrm{rot}}$. Lemmas~\ref{lem:algebra} and~\ref{lem:C1} thus hold as well on $C_{\mathrm{rot}}$. Lemma~\ref{lem:A} and the bound \eqref{eq:Y0-def} hold with all sums and maxima restricted to $\mathcal{I}_{\mathrm{rot}}$ and $K_{\mathrm{rot}}$, and the sum in \eqref{eq:Y0-def} is complete because $\mathcal{N}(\bar U_{\mathrm{rot}})\in C_{\mathrm{rot}}$. The definition of $Z_0$ and the proof of \eqref{eq:Z0} carry over with $K$ replaced by $K_{\mathrm{rot}}$. The proof of Lemma~\ref{lem:Z2} uses only the inequalities of Lemmas~\ref{lem:algebra} and~\ref{lem:A}, and it applies to $h$ and $h'$ in $\mathcal{X}_{\mathrm{rot}}$.

For the bound $Z_1$, we note that the identity \eqref{eq:Z1-identity} holds in $\mathcal{Y}_{\mathrm{rot}}$, and that $\mathcal{D} e_i\in C_{\mathrm{rot}}$ for $i\in\mathcal{I}_{\mathrm{rot}}$ by Lemma~\ref{lem:rot-system}. Thus we have $\mathcal{D}_{i'i}=0$ for $i\in\mathcal{I}_{\mathrm{rot}}$ and $i'\notin\mathcal{I}_{\mathrm{rot}}$. Hence the sums in \eqref{eq:columns} for $c_\sigma$ and $c_i$, $i\in\mathcal{I}_{\mathrm{rot}}$ involve only indices $i'\in\mathcal{I}_{\mathrm{rot}}$. By \eqref{eq:opnorm} on $\mathcal{X}_{\mathrm{rot}}$, it follows that the bound of Lemma~\ref{lem:Z1} holds with the supremum taken over $\mathcal{I}_{\mathrm{rot}}$. The proof of Lemma~\ref{lem:Z1} for the columns with $\theta_i>\Theta'$ uses the splitting $\mathcal{D}=\mathcal{D}^{\leq}+\mathcal{D}^{>}$, and it requires that both parts map $C_{\mathrm{rot}}$ into itself. Each term of the matrices \eqref{eq:M} is either a tensor product of two of the fields $\bar X_{\mathrm{rot}}$, $\partial_x\bar X_{\mathrm{rot}}$, and $\partial_t\bar X_{\mathrm{rot}}$, or the identity multiplied by a function independent of $x$. Since these fields satisfy $Y(x+\alpha,t)=R_\alpha Y(x,t)$, it follows that
\begin{align*}
M_p(x+\alpha,t)=R_\alpha M_p(x,t)R_\alpha^{-1},\qquad p\in\{0,1,2\}.
\end{align*}
The same holds for the truncations $M_p^{\leq}$ and $M_p^{>}$, and therefore the operators $\mathcal{D}^{\leq}$ and $\mathcal{D}^{>}$ commute with the maps $X\mapsto R_\alpha X(\cdot-\alpha,\cdot)$. By the proof of Lemma~\ref{lem:Z1} they also map $C_G$ into itself, and hence they preserve $C_{\mathrm{rot}}$. The rest of the proof carries over to the rotational case unchanged. Thus Theorem~\ref{thm:radii} applies to $\bar U_{\mathrm{rot}}$, with the constants \eqref{eq:Y0-def}, \eqref{eq:Z1-def}, and \eqref{eq:Z2-def} computed over $\mathcal{I}_{\mathrm{rot}}$.

We evaluate these constants as in \S\ref{sec:computation}, with the index sets restricted to $\mathcal{I}_{\mathrm{rot}}$ and with the truncations $\Theta_{\mathrm{rot}}=500$ and $\Theta_0=400$, so that $\Theta'=900$. In this setting all quantities are computed in exact rational arithmetic, including the matrix $\mathcal{A}^{\mathrm{rot}}_K$, whose binary64 entries are used without rounding.
The column sums are computed for the $923$ indices in $K_{\mathrm{rot}}$ and for the $738$ indices in $\mathcal{I}_{\mathrm{rot}}$ with $\Theta_{\mathrm{rot}}<\theta_i\leq\Theta'$, and the rows which occur are the $3507$ indices in $\mathcal{I}_{\mathrm{rot}}$ with $\theta_{i'}\leq1900$. The certified bounds are given in Table~\ref{tab:rot-bounds}. We put
\begin{align*}
r^{\mathrm{rot}}_{\min}=5.5798218604081895\cdot10^{-14},\qquad r^{\mathrm{rot}}_{\max}=2.3485615439327366\cdot10^{-4}.
\end{align*}

\begin{table}[ht]
\centering
\small
\begin{tabular}{ll}
\toprule
quantity & certified upper bound\\
\midrule
$b_0$, $b_1$, $b_2$ & $8.09152$,\ \ $6.62942$,\ \ $10.8375$\\
$\|M_0\|$, $\|M_1\|$, $\|M_2\|$ & $11.5576$,\ \ $7.21986$,\ \ $17.1543$\\
$m_0$ & $2.64802\cdot10^{-17}$\\
$m_1$ & $6.34327\cdot10^{-20}$\\
$m_2$ & $2.69708\cdot10^{-17}$\\
$\tau$, $\bar\sigma_{\mathrm{rot}}^{1/2}$ & $0.922834$,\ \ $1.08362$\\
$Y_0$ & $5.31985\cdot10^{-14}$\\
$Z_0$ & $6.38882\cdot10^{-14}$\\
$\alpha_\eta$ & $10.4201$\\
$\alpha_{\max}$, $\alpha_{tt}$ & $42.0212$,\ \ $21.5162$\\
$c_\sigma$ & $4.09649\cdot10^{-22}$\\
$\max\{c_i : i\in K_{\mathrm{rot}}\}$ & $0.0142888$\\
$\max\{c_i : \Theta_{\mathrm{rot}}<\theta_i\leq\Theta'\}$ & $0.0318485$\\
$Z_1^{\mathrm{far}}$ & $0.0461264$\\
$Z_1$ & $0.0461264$\\
$Z_2(r^{\mathrm{rot}}_{\min})$, $Z_2(r^{\mathrm{rot}}_{\max})$ & $4060.79$,\ \ $4060.87$\\
\bottomrule
\end{tabular}
\medskip
\caption{Certified bounds for the rotationally symmetric problem, with $\Theta_{\mathrm{rot}}=500$ and $\Theta_0=400$. In addition $\Theta_+\geq500.632$.}
\label{tab:rot-bounds}
\end{table}

\begin{proposition}\label{prop:cap-rot}
There exist a trigonometric polynomial $\bar X_{\mathrm{rot}}$ as above and a $924\times924$ matrix $\mathcal{A}_K^{\mathrm{rot}}$ with dyadic rational entries such that the quantities of \S\ref{sec:functional} for $\bar U_{\mathrm{rot}}$, $\Theta_{\mathrm{rot}}=500$, and $\Theta_0=400$ satisfy the upper bounds of Table~\ref{tab:rot-bounds}, $\Theta_+\geq500.632$, and
\begin{gather*}
2.1823\cdot10^{-16}\leq\eta(\bar U_{\mathrm{rot}})\leq2.1824\cdot10^{-16},\\
p(r^{\mathrm{rot}}_{\min})\leq-2.59758\cdot10^{-17},\qquad p(r^{\mathrm{rot}}_{\max})\leq-3.64610\cdot10^{-8}.
\end{gather*}
Moreover,
\begin{gather*}
\min_{\mathbb{T}^2}|\partial_x\bar X_{\mathrm{rot}}|\geq0.143486,\qquad \sup_{\mathbb{T}^2}|\bar X_{\mathrm{rot}}|\leq3.66695,\qquad \sup_{\mathbb{T}^2}|\partial_x\bar X_{\mathrm{rot}}|\leq3.60835,\\
1.8512166716821\leq\mathcal{E}(\bar U_{\mathrm{rot}})\leq1.8512166716822.
\end{gather*}
\end{proposition}

\begin{proof}
The coefficients of $\bar X_{\mathrm{rot}}$ and the entries of $\mathcal{A}_K^{\mathrm{rot}}$ are listed in the files \texttt{ubar\_T500.json} and \texttt{s3\_AK.npy} in the directory \texttt{data/rot} of \cite{code}. They were obtained by Newton's method for the system of Lemma~\ref{lem:rot-system} and by a numerical inversion, both in floating-point arithmetic. The bounds are checked by the program \texttt{verify} of \cite{code}, in exact rational arithmetic and in ball arithmetic, by the computation described in \S\ref{sec:computation} and above. For the last four bounds we use the method of \S\ref{sec:post} with a grid of $25\times83$ points of the fundamental domain, on which the minimum of $|\partial_x\bar X_{\mathrm{rot}}|$ is at least $0.437123$ and the remainder is at most $0.293638$, and the trapezoidal rule with $2304\times2048$ nodes.
\end{proof}

In the rest of this section, $\bar X_{\mathrm{rot}}$ and $\mathcal{A}_K^{\mathrm{rot}}$ are as in Proposition~\ref{prop:cap-rot}.

\begin{proposition}\label{prop:rot-existence}
The map $\mathcal{F}$ has a zero $U_{\mathrm{rot}}=(\sigma_{\mathrm{rot}},X_{\mathrm{rot}})\in\mathcal{X}_{\mathrm{rot}}$ with $\|U_{\mathrm{rot}}-\bar U_{\mathrm{rot}}\|_{\mathcal{X}}\leq r^{\mathrm{rot}}_{\min}$. It is the only zero of $\mathcal{F}$ in the closed ball of radius $r^{\mathrm{rot}}_{\max}$ about $\bar U_{\mathrm{rot}}$ in $\mathcal{X}_{\mathrm{rot}}$.
\end{proposition}

\begin{proof}
By Proposition~\ref{prop:cap-rot} we have $Z_0<1$, so that $\mathcal{A}_{\mathrm{rot}}$ is injective by Lemma~\ref{lem:A}. The hypotheses of Theorem~\ref{thm:radii} hold for the restriction of $\mathcal{F}$ to $\mathcal{X}_{\mathrm{rot}}$ by the discussion above, and by Proposition~\ref{prop:cap-rot} the inequality \eqref{eq:radii} holds for $r_0=r^{\mathrm{rot}}_{\min}$ and for $r_0=r^{\mathrm{rot}}_{\max}$. The assertion follows as in the proof of Proposition~\ref{prop:existence}.
\end{proof}

The uniqueness in Proposition~\ref{prop:rot-existence} holds in the ball of radius $r^{\mathrm{rot}}_{\max}$ about $\bar U_{\mathrm{rot}}$ in $\mathcal{X}_{\mathrm{rot}}$, that is, among the surfaces of revolution about the $X_3$-axis which are symmetric under $\iota_3$ and are parametrized conformally as in \eqref{eq:rot-ansatz}.

\subsection{Proof of Theorem~\ref{thm:rot}}

Let $U_{\mathrm{rot}}=(\sigma_{\mathrm{rot}},X_{\mathrm{rot}})$ be the zero of Proposition~\ref{prop:rot-existence}, with $X_{\mathrm{rot}}$ of the form \eqref{eq:rot-ansatz} with functions $a_{\mathrm{rot}}$ and $b_{\mathrm{rot}}$, and put $\Sigma_{\mathrm{rot}}=X_{\mathrm{rot}}(\mathbb{T}^2)$. We write $r_0=r^{\mathrm{rot}}_{\min}$. As in \eqref{eq:pointwise} we have $|X_{\mathrm{rot}}-\bar X_{\mathrm{rot}}|\leq r_0$ and $|\partial_xX_{\mathrm{rot}}-\partial_x\bar X_{\mathrm{rot}}|\leq r_0$ pointwise on $\mathbb{T}^2$, and $|\sigma_{\mathrm{rot}}-\bar\sigma_{\mathrm{rot}}|\leq r_0/\vartheta_{\mathrm{rot}}$. In particular $\sigma_{\mathrm{rot}}>0$.

The pair $U_{\mathrm{rot}}$ is a zero of $\mathcal{F}$ with $X_{\mathrm{rot}}\in C_G$ of class $C^1$. By Proposition~\ref{prop:zeros} the map $X_{\mathrm{rot}}$ is real-analytic and weakly conformal. By Proposition~\ref{prop:cap-rot} we have $|\partial_xX_{\mathrm{rot}}|\geq0.143486-r_0>0$. Since $|\partial_xX_{\mathrm{rot}}|=|a_{\mathrm{rot}}|$, the function $a_{\mathrm{rot}}$ does not vanish. Thus $X_{\mathrm{rot}}$ is a conformal immersion satisfying \eqref{eq:shrinker} by Proposition~\ref{prop:zeros}, and $\Sigma_{\mathrm{rot}}$ is a surface of revolution about the $X_3$-axis symmetric under $\iota_3$ by \S\ref{sec:rot-subspace}.

For maps in $C_{\mathrm{rot}}$, the Gaussian area \eqref{eq:FJ} is
\begin{align*}
\mathcal{E}(\sigma,X)=\frac{1}{2\sqrt\sigma}\int_0^{2\pi}e^{-(a^2+b^2)/4}a^2\,dt.
\end{align*}
Lemma~\ref{lem:perturbation} applies to $U_{\mathrm{rot}}$ and $\bar U_{\mathrm{rot}}$, with the bounds of Proposition~\ref{prop:cap-rot} for $\sup|\bar X_{\mathrm{rot}}|$ and $\sup|\partial_x\bar X_{\mathrm{rot}}|$ in place of $a_0$ and $a_1$, and with $\delta=r_0/\vartheta_{\mathrm{rot}}$. We find $\kappa\leq31.09$ and $|\mathcal{E}(U_{\mathrm{rot}})-\mathcal{E}(\bar U_{\mathrm{rot}})|\leq5.04\cdot10^{-12}$. Together with the enclosure of $\mathcal{E}(\bar U_{\mathrm{rot}})$ in Proposition~\ref{prop:cap-rot}, this gives
\begin{equation}\label{eq:rot-entropy}
1.85121667167\leq\mathcal{E}(U_{\mathrm{rot}})\leq1.85121667169.
\end{equation}
In particular the Gaussian area of the immersion $X_{\mathrm{rot}}$ is less than $2$. By Lemma~\ref{lem:embedded} the map $X_{\mathrm{rot}}$ is an embedding and $\lambda(\Sigma_{\mathrm{rot}})=F(\Sigma_{\mathrm{rot}})=\mathcal{E}(U_{\mathrm{rot}})$. Thus $\Sigma_{\mathrm{rot}}$ is a real-analytic closed embedded surface of genus one, and \eqref{eq:rot-entropy} is the enclosure of the entropy in Theorem~\ref{thm:rot}.

By \eqref{eq:ball-enclosures} we have $\lambda(\Sigma)=\mathcal{E}(U_*)\geq1.8525673197$. Together with \eqref{eq:rot-entropy} this gives
\begin{align*}
\lambda(\Sigma)-\lambda(\Sigma_{\mathrm{rot}})\geq1.8525673197-1.85121667169>1.35064\cdot10^{-3}.
\end{align*}
This completes the proof of Theorem~\ref{thm:rot}. Corollary~\ref{cor:least} follows, since $\Sigma_{\mathrm{rot}}$ is a closed embedded self-shrinker of genus one with smaller entropy.

The upper bounds in \eqref{eq:ball-enclosures} and the lower bound in \eqref{eq:rot-entropy} give in the same way
\begin{align*}
\lambda(\Sigma)-\lambda(\Sigma_{\mathrm{rot}})\leq1.36286\cdot10^{-3}.
\end{align*}

\begin{remark}
As in \S\ref{sec:geometry}, the surface $\Sigma_{\mathrm{rot}}$ is conformally equivalent to the flat torus $\R^2/(2\pi\Z\times2\pi c_{\mathrm{rot}}\Z)$ with $c_{\mathrm{rot}}=\sigma_{\mathrm{rot}}^{-1/2}$. Since $|\sigma_{\mathrm{rot}}-\bar\sigma_{\mathrm{rot}}|\leq r_0/\vartheta_{\mathrm{rot}}$, we have
\begin{align*}
0.922833051919\leq c_{\mathrm{rot}}\leq0.922833051920.
\end{align*}
\end{remark}

\begin{remark}
Theorem~\ref{thm:radii} also shows that $\mathcal{A}_{\mathrm{rot}}D\mathcal{F}(U_{\mathrm{rot}})$ is invertible on $\mathcal{X}_{\mathrm{rot}}$, so that $D\mathcal{F}(U_{\mathrm{rot}})$ is injective on $\mathcal{X}_{\mathrm{rot}}$. The argument of Propositions~\ref{prop:jacobi} and~\ref{prop:nondegenerate} can be repeated in the rotational class. Let $f$ be a function on $\Sigma_{\mathrm{rot}}$ with $Lf=0$ which is invariant under the rotations about the $X_3$-axis and under $\iota_3$. Then $f\circ X_{\mathrm{rot}}$ depends only on $t$ and is even. In the proof of the second part of Proposition~\ref{prop:jacobi} the function $R_1$ depends only on $t$, the function $R_2$ vanishes because $h_{xt}=0$ for a map of the form \eqref{eq:rot-ansatz}, and consequently $\phi=0$, the function $\psi$ depends only on $t$, and $\xi=fn+\psi\partial_tX_{\mathrm{rot}}$ belongs to $C_{\mathrm{rot}}$. Hence $f=0$. Thus $\Sigma_{\mathrm{rot}}$ has no nonzero Jacobi field which is invariant under the rotations about its axis and under $\iota_3$.
\end{remark}

\begin{remark}
The properties of $\Sigma_{\mathrm{rot}}$ in Theorem~\ref{thm:rot} do not identify it with the torus of \cite{An} or with the torus of \cite{DN}, since no uniqueness theorem is available for embedded shrinking tori of revolution. For complete embedded rotationally symmetric self-shrinkers which are symmetric with respect to the plane orthogonal to the axis, the theorem of Ma, Muhammad, and M{\o}ller \cite[Theorem~1.3]{MMM} gives finiteness up to rigid motions.
\end{remark}

\appendix

\section{Uncertified computations}\label{sec:numerics}

We compare the rotational profile obtained from the conformal system with one obtained by shooting for Angenent’s equation. We also give numerical spectra for both surfaces. These computations are not certified and are not used in the proofs.

\subsection{The rotationally symmetric torus}\label{sec:numerics-rot}

We computed the rotational approximation in multiple-precision arithmetic using the system of Lemma~\ref{lem:rot-system}. Independently, we computed a rotational shrinking profile by shooting for Angenent's equation \cite{An}. The two computed profiles differ by at most $4.8\cdot10^{-13}$, and their Gaussian areas differ by $2.4\cdot10^{-15}$. The value $1.85121667168$ agrees with Berchenko-Kogan's numerical value $1.85122$ \cite{BK1}. We therefore expect that $\Sigma_{\mathrm{rot}}$ is the Angenent torus, although we do not prove that it arises from the same construction employed in \cite{An}.

\subsection{The spectrum of the Jacobi operator}

Table~\ref{tab:spectrum} lists the lowest computed eigenvalues of $L$ by symmetry type. We solved the generalized eigenvalue problem $\mathcal{S}\varphi=\beta\ell\varphi$ of Proposition~\ref{prop:inertia} by a Galerkin method in each sector, using the basis of Table~\ref{tab:sector-def} with cut-offs $\theta_k\leq140$ for $\Sigma$ and $\theta_k\leq120$ for the rotational torus.

\begin{table}[ht]
\centering
\begin{tabular}{lll}
\toprule
type & $\Sigma$ & $\Sigma_{\mathrm{rot}}$\\
\midrule
$A_1'$ & $-4.1514$,\ \ $-1$,\ \ $-0.1507$,\ \ $1.0065$,\ \ $1.6825$ & $-3.7398$,\ \ $-1$,\ \ $0.1129$,\ \ $0.9920$\\
$A_2'$ & $0$,\ \ $1.8820$,\ \ $3.7530$ & $0.1129$,\ \ $3.0493$\\
$A_2''$ & $-\frac12$,\ \ $0.9142$,\ \ $2.0700$ & $-\frac12$,\ \ $1.2977$,\ \ $1.8626$\\
$A_1''$ & $1.7523$,\ \ $4.0973$ & $1.8626$\\
$E'$ & $-1.0392$,\ \ $-\frac12$,\ \ $-0.4991$,\ \ $0.4445$,\ \ $0.9879$ & $-1$,\ \ $-\frac12$,\ \ $-0.4876$,\ \ $0.9037$,\ \ $1.7270$\\
$E''$ & $0$,\ \ $0.7826$,\ \ $1.6737$ & $0$,\ \ $0.8640$,\ \ $2.0403$\\
\bottomrule
\end{tabular}
\medskip
\caption{Numerical values of the lowest eigenvalues of the Jacobi operator, computed in floating-point arithmetic.}
\label{tab:spectrum}
\end{table}

For $\Sigma$, the numerical counts agree with the certified index and nullity reported in Theorem~\ref{thm:index}. The computation also resolves an eigenvalue corresponding to class $E'$ near $-0.4991$, close to the translation eigenvalue $-\frac12$.

For $\Sigma_{\mathrm{rot}}$, the computed index $9$ and nullity $2$ agree with those found by Berchenko-Kogan \cite{BK2}. The lowest eigenvalue of $\Sigma_{\mathrm{rot}}$ with angular frequency $3$ is approximately $0.1129$, with eigenfunctions $a(t)\cos(3x)$ of type $A_1'$ and $a(t)\sin(3x)$ of type $A_2'$. On $\Sigma$, these types instead contain the eigenvalue near $-0.1507$ and the rotational zero mode $\langle e_3\times X,n\rangle$, which vanishes identically on $\Sigma_{\mathrm{rot}}$ and, to first order in the amplitude of the threefold deformation, is a multiple of $a(t)\sin(3x)$. This suggests that along the path of \S\ref{sec:discovery} the pair of eigenvalues near $0.1129$ splits into the additional negative eigenvalue and the rotational zero mode of $\Sigma$.

\bibliographystyle{plain}
\bibliography{references}

\end{document}